\documentclass[a4paper,11pt,reqno]{amsart}
\usepackage[utf8x]{inputenc}
\usepackage{amsmath}
\usepackage{amssymb}
\usepackage{amsthm}
\usepackage{amsfonts}
\usepackage{mathtools}
\usepackage{array}
\usepackage{cite}
\usepackage{enumerate}
\usepackage{xcolor}
\usepackage{bbm}
\usepackage{esint}
\usepackage{graphicx}
\usepackage{float}
\usepackage{setspace}
\usepackage[hidelinks]{hyperref}
\usepackage{esint}
\usepackage[margin=3cm]{geometry}
\usepackage{mathrsfs}
\usepackage{stmaryrd}
\allowdisplaybreaks

\newtheorem{theorem}{Theorem}[section]
\newtheorem{definition}[theorem]{Definition}
\newtheorem*{definition*}{Definition}
\newtheorem*{theorem*}{Theorem}
\newtheorem{lemma}[theorem]{Lemma}
\newtheorem*{lemma*}{Lemma}
\newtheorem{corollary}[theorem]{Corollary}
\newtheorem{example}[theorem]{Example}

\newtheorem{remark}[theorem]{Remark}
\newtheorem*{remark*}{Remark}

\newtheorem{proposition}[theorem]{Proposition}

\newtheorem{alphtheorem}{Theorem}

\newcommand{\sbullet}{\,\begin{picture}(1,1)(-0.5,-2)\circle*{2}\end{picture}\,}
\newcommand{\wsc}{\overset{*}{\rightharpoondown}}
\newcommand{\ip}[2]{\left<#1,#2\right>}
\newcommand{\restrict}{	\begin{picture}(10,8)\put(2,0){\line(0,1){7}}\put(1.8,
0){\line(1,0){7}}\end{picture}}

\newcommand{\pd}{\partial}
\newcommand{\dd}{\mathrm{d}}

\newcommand{\mbR}{\mathbb{R}}
\newcommand{\mbN}{\mathbb{N}}
\newcommand{\mbfM}{\mathbf{M}}

\newcommand{\mbfR}{\mathbf{R}}

\newcommand{\Lp}{\mathrm{L}}
\newcommand{\W}{\mathrm{W}}
\newcommand{\C}{\mathrm{C}} 
\newcommand{\bB}{\mathbb{B}}

\newcommand{\lL}{\mathcal{L}^d}

\newcommand{\Lrm}{\mathrm{L}}

\newcommand{\Acal}{\mathcal{A}}
\newcommand{\Bcal}{\mathcal{B}}
\newcommand{\Ccal}{\mathcal{C}}
\newcommand{\Dcal}{\mathcal{D}}
\newcommand{\Ecal}{\mathcal{E}}
\newcommand{\Fcal}{\mathcal{F}}
\newcommand{\Gcal}{\mathcal{G}}
\newcommand{\Hcal}{\mathcal{H}}

\newcommand{\Jcal}{\mathcal{J}}

\newcommand{\Lcal}{\mathcal{L}}
\newcommand{\Mcal}{\mathcal{M}}
\newcommand{\Ncal}{\mathcal{N}}

\newcommand{\Hfrak}{\mathfrak{H}}

\newcommand{\Mbf}{\mathbf{M}}

\newcommand{\Rbf}{\mathbf{R}}

\DeclareMathOperator{\id}{id}

\DeclareMathOperator{\dist}{dist}

\DeclareMathOperator{\tr}{tr}

\DeclareMathOperator{\supp}{supp}

\DeclareMathOperator{\gr}{gr}

\newcommand{\setsmall}[2]{\{\, #1 \, \textup{\textbf{:}}\, #2 \,\}}

\newcommand{\setBB}[2]{\biggl\{\, #1 \ \ \textup{\textbf{:}}\ \ #2 \,\biggr\}}

\newcommand{\norm}[1]{\|#1\|}

\newcommand{\abs}[1]{|#1|}

\newcommand{\floor}[1]{\lfloor #1 \rfloor}

\newcommand{\floorb}[1]{\bigl\lfloor #1\bigr\rfloor}

\newcommand{\N}{\mathbb{N}}
\newcommand{\R}{\mathbb{R}}

\newcommand{\todown}{\downarrow}

\newcommand{\frarg}{\,\sbullet\,}
\newcommand{\BV}{\mathrm{BV}}

\DeclarePairedDelimiter{\asc}{[}{]}

\DeclareMathOperator*{\wstarlim}{w*-lim}

\newcommand{\RBVL}{\mathbf{R}^1}
\newcommand{\Fcalrw}{\mathcal{F}^{w*}_{**}}
\newcommand{\Fcalro}{\mathcal{F}^{\,1}_{**}}

\def\Xint#1{\mathchoice
{\XXint\displaystyle\textstyle{#1}}%
{\XXint\textstyle\scriptstyle{#1}}%
{\XXint\scriptstyle\scriptscriptstyle{#1}}%
{\XXint\scriptscriptstyle\scriptscriptstyle{#1}}%
\!\int}
\def\XXint#1#2#3{{\setbox0=\hbox{$#1{#2#3}{\int}$ }
\vcenter{\hbox{$#2#3$ }}\kern-.6\wd0}}

\def\dashint{\Xint-}

\newcommand{\RED}[1]{{\color{black}#1}}

\title[Relaxation for partially coercive integral functionals]{Relaxation for partially coercive integral functionals with linear growth}

\author{Filip Rindler}
\address[Filip Rindler]{Mathematics Institute, University of Warwick, Coventry CV4 7AL, UK, and The Alan Turing Institute, British Library, 96 Euston Road, London NW1 2DB, UK}
\email{F.Rindler@warwick.ac.uk}

\author{Giles Shaw}
\address[Giles Shaw]{University of Reading, Department of Mathematics and Statistics, Whiteknights, PO Box 220, Reading RG6 6AX, UK}
\email{giles.shaw@gmail.com}

\begin{document}		

\begin{abstract}
We prove an integral representation theorem for the $\Lp^1$-relaxation of the functional
\[
\Fcal\colon u\mapsto\int_\Omega f(x,u(x),\nabla u(x))\;\dd x,\quad u\in\W^{1,1}(\Omega;\R^m),
\]
where $\Omega\subset\R^d$ ($d \geq 2$) is a bounded Lipschitz domain, to the space $\BV(\Omega;\R^m)$ under very general assumptions: we require principally that $f$ is Carath\'eodory, that the partial coercivity and linear growth bound
\[
  g(x,y)|A|\leq f(x,y,A)\leq Cg(x,y)(1+|A|),
\]
holds, where $g\colon\overline{\Omega}\times\R^m\to[0,\infty)$ is a continuous function satisfying a weak monotonicity condition, and that $f$ is quasiconvex in the final variable. Our result is the first that applies to integrands which are unbounded in the $u$-variable and \RED{therefore allows for the treatment of} many problems from applications. Such functionals are out of reach of the classical blow-up approach introduced by Fonseca \& M\"uller [\emph{Arch.\ Ration.\ Mech.\ Anal.} 123 (1993), 1--49]. Our proof relies on an intricate truncation construction (in the $x$ and $u$ arguments simultaneously) made possible by the theory of liftings developed in a previous paper by the authors~[\emph{Arch.\ Ration.\ Mech.\ Anal.}, 232 (2019), 1227--1328], and features techniques which could be of use for other problems involving $u$-dependent integrands. 
\end{abstract}

\maketitle
\author



\section{Introduction}\label{secintro}

Inspired by problems in the theory of phase transitions, this paper is concerned with the identification of the integral representation for the relaxation \RED{$\Fcalro$} to $\BV(\Omega;\R^m)$ with respect to the strong $\Lp^1(\Omega;\R^m)$-topology of the functional
\begin{equation}\label{eqoriginalfunctional}
\Fcal[u]:=\int_\Omega f(x,u(x),\nabla u(x))\;\dd x,\quad u\in\W^{1,1}(\Omega;\R^m),
\end{equation}
\RED{where $\Omega\subset\R^d$ ($d \geq 2$) is a bounded Lipschitz domain} and $f \colon \Omega \times \R^m \times \R^{m \times d} \to [0,\infty)$ is quasiconvex and has linear growth in the final variable. That is, we aim to compute
\begin{equation} \label{eq:firstFcalro}
\Fcalro[u]:=\inf\left\{\liminf_{j\to\infty}\Fcal[u_j]\colon\;(u_j)_j\subset\W^{1,1}(\Omega;\R^m)\text{ and }u_j\to u\in\BV(\Omega;\R^m)\text{ in $\Lp^1$}\right\}
\end{equation}
at $u \in \BV(\Omega;\R^m)$.

To motivate this problem and to explain the main difficulties in identifying $\Fcalro$, consider the task of computing the $\Gamma$-limit of a sequence of singularly perturbed functionals of the form
\begin{equation}\label{eqperturbfamily}
\Ecal_\varepsilon[u]:=\varepsilon^{-1}\int_\Omega \left[g(x,u(x))\right]^2\;\dd x+\varepsilon\int_\Omega \left[h(x,u(x),\nabla u(x))\right]^2\;\dd x.
\end{equation}
Such questions occur in a variety of applied contexts where $g$ typically models the energy of different phase mixtures and $h$ is a bulk energetic cost for movement between phases: in applications arising from fluid phase transition problems, $g$ is often of the form $|a-y|^2|y-b|^2$ for vectors $a\in\R^m$, $b\in\R^m$ representing preferred phases for a vector of fluid densities $u$ and $h(x,y,A)=|A|$ penalises variations away from a constant phase~\cite{Modi87TGPT,Gurt87SRCG,Bald90MICP}. In the study of phase transitions within elastic solids, $u$ satisfies the constraint $\operatorname{curl}u=0$ (and is therefore a gradient) and $g$ vanishes on the rotation orbits $\{RA\colon R\in\operatorname{SO}(3)\}$ and $\{RB\colon R\in\operatorname{SO}(3)\}$ of two rank-one connected matrices $A,B\in\R^{3\times 3}$~\cite{Fonse89PTEM,CoFoGi02AGCR}. From the perspective of smectic liquid crystals, $u$ is again curl-free, but now $g(x,y)=||y|-1|$ and $h(x,y,A)=|\tr(A)|$,~\cite{AviGig87AMPR}. In the theory of harmonic maps and also applications from reaction diffusion processes for chemical reactions, $g=\dist(\frarg,\Ncal)^2$ for some closed target Riemannian manifold $\Ncal\subset\R^m$ and $h^2(\nabla u)=g^{\alpha,\beta}\ip{\pd_\alpha u}{\pd_\beta u}$ is the Dirichlet integral associated to a $d$-dimensional domain Riemannian manifold $\Mcal$,~\cite{RuStKe89RDPE,LiPaWa12PTPH,CheStru89EPRR}.

If we set $f(x,y,A):=g(x,y)h(x,y,A)$, then the Cauchy--Schwarz inequality immediately implies that computing $\Fcalro$ provides a (usually optimal) lower bound for the $\Gamma$-limit of $(\Ecal_\varepsilon)_{\varepsilon>0}$ as $\varepsilon \to 0$. Thus, minimisers of $\Fcalro$ are ``physically reasonable'' solutions to the energy minimisation for $\int_\Omega g(x,u(x))\;\dd x$, taking into account the fact that transitions between phases should have an energetic cost.

\RED{In some phase-transition problems only (strong) compactness in $\Lp^1$ can be established rather than full $\BV$-weak* compactness, see for instance~\cite{FonTar89GTPT}. The need therefore arises to study~\eqref{eqperturbfamily} under just an $\Lp^1$-compactness assumption. At present it is not clear what the right limiting space is for many such problems, but we hope that a better understanding of the problem in the case where we a-priori assume that the limit lies in $\BV$, yields valuable new insights. In fact, our results are such that in a more general setup they would provide the required relaxation lower bounds on the functional ``where the limit map is of bounded variation'' (which may not be everywhere, of course).}

From a mathematical perspective, $\Lp^1(\Omega;\R^m)$-relaxation problems with linear growth can be very badly behaved. It is known (see~\cite{FonLeo01LSC,GorMar02AESL} for detailed discussions of this topic) that, if no coercivity is assumed on $f$, no sensible formula for $\Fcalro$ is possible without strong additional assumptions in the $(x,y)$ variables. Indeed, an example of Dal Maso~\cite{DMas79IRBV} shows that there exists a continuous ($u$-independent) integrand $f\colon\Omega\times\R^d\to[0,\infty)$ which is both convex and positively one-homogeneous in the final variable, but for which $\Fcal$ is not equal to $\Fcalro$ over $\W^{1,1}(\Omega;\R)$ (despite the convexity of $f$!). On the other hand, it was recently shown in~\cite{RindlerShaw19} that a satisfactory integral formula for the sequential weak* relaxation, $\Fcalrw$, of $\Fcal$ to $\BV(\Omega;\R^m)$ (i.e., $\Fcalrw[u]$ is defined analogously to $\Fcalro[u]$ but with $\Lp^1$-convergence replaced by weak* convergence in $\BV$) does always exist for essentially any Carath\'eodory integrand $f$ which is quasiconvex and of linear growth in the final variable. However, integrands arising from limits of the kind of models given by~\eqref{eqperturbfamily} are not coercive, so it need not be the case that every sequence $(u_j)_j\subset\W^{1,1}(\Omega;\R^m)$ with $\limsup_j\Fcal[u_j]<\infty$ can be assumed to be weakly* convergent. However, they can all be assumed to be \emph{partially coercive}: that is, there exist $g\in\C(\Omega\times\R^m;[0,\infty))$ and $C>0$ such that
\begin{equation}\label{eqpartialcoercivity}
g(x,y)|A|\leq f(x,y,A)\leq Cg(x,y)(1+|A|)\quad\text{ for all }\quad(x,y,A)\in\Omega\times\R^m\times\R^{m\times d}.
\end{equation}
It turns out that partial coercivity implies that $\Fcal$ is coercive in small cylinders $B^d(x,r)\times B^m(y,R)$ about every pair $(x,y)\subset\Omega\times\R^m$ which `matters from the perspective of computing $\Fcal$ and $\Fcalro$'. In order to derive a formula for $\Fcalro$ which makes minimal assumptions on $f$, our approach is to find a way of making the preceding statement rigorous in order to reduce the computation of $\Fcalro$ to an application of the weak* technology developed in~\cite{RindlerShaw19}.

An additional difficulty in computing $\Fcalro$ arises from the fact that, for the prototypical integrand $f(x,y,A)=g(x,y)h(x,y,A)$, the function $g$ could in principle be \emph{any} energy density on $\Omega\times\R^m$. Consequently, $f$ might exhibit arbitrary growth in the $y$-variable and hence $\Fcalro$ need not be finite over all of $\BV(\Omega;\R^m)$. 

The study of these functionals in the vector-valued case originates in~\cite{FonMul93RQFB} (see also~\cite{AmbDal92RBVQ} for the $u$-independent case), where the authors showed that, if $f$ is quasiconvex in the final variable, partially coercive in the sense of~\eqref{eqpartialcoercivity}, and additionally satisfies some strong \emph{localisation hypotheses} \RED{(see~\eqref{eqfonassumption} in Section~\ref{secrelaxationcomparisons} and the comments preceding it)}, then
\begin{align}
\begin{split}\label{eqfunctional}
\RED{\Fcalro}[u] &= \int_\Omega f(x,u(x),\nabla u(x))\;\dd x+\int_\Omega f^\#\left(x,u(x),\frac{\dd D^cu}{\dd|D^c u|}(x)\right)\;\dd|D^cu|(x) \\
  &\qquad +\int_{\mathcal{J}_u}K_f[u](x)\;\dd\mathcal{H}^{d-1}(x),  
\end{split}
\end{align}
where $f^\#$ is the \RED{generalised (limsup)} recession function of $f$, which is defined by $f^\#(A):=\limsup_{t\to\infty}t^{-1}f(tA)$ \RED{(see the discussion after Theorem~\ref{L1lscthm} for a comparison of different types of recession functions)}, and
\[
\begin{aligned}
K_f[u](x):= \inf\, \setBB{\frac{1}{\omega_{d-1}}\int_{\bB^d}f^\#(x,\varphi(y),\nabla\varphi(y))\;\dd y}{&\varphi\in\C^\infty(\bB^d;\mbR^m), \\
&{\varphi}|_{\pd\bB^d}=u^\pm(x)\text{ if }\ip{y}{n_u(x)} \gtrless 0}
\end{aligned}
\]
is a generalised surface energy density associated with $f$ \RED{(this definition of $K_f$ has a different formulation compared to \RED{the one} in~\cite{FonMul93RQFB} but is equivalent, see Remark~\ref{rem:K_f})}. Here, we have used the usual decomposition
\[
  Du = \nabla u \Lcal^d + D^s u,  \qquad D^s u = (u^+ - u^-) \otimes n_u \, \Hcal^{d-1} \restrict \Jcal_u + D^c u,
\]
for the derivative $Du$ of a function $u\in\BV(\Omega;\R^m)$ (see, for example~\cite{AmFuPa00FBVF}). 

The localisation hypotheses found in~\cite{FonMul93RQFB} and its descendants, see for instance~\cite{BoFoMa98GMR,FonLeo01LSC}, are required for technical reasons, but have the undesirable consequence of precluding the application of their result to several integrands of potential interest (such as $[\dist(y,K)]^{p(x)}|A|$, \RED{where $K\subset\R^m$ is compact and $p\colon\overline{\Omega}\to[0,\infty)$ is continuous}, for instance), which can arise naturally from the applications discussed above. All previously available results additionally require that $f$ is bounded in the $y$ variable, which is also incompatible with integrands arising from the applications listed above.

In this paper we build on the weak* relaxation results developed in~\cite{RindlerShaw19} to establish an improved relaxation theory for $\Fcal$ as defined by~\eqref{eqoriginalfunctional} under natural conditions on the integrand $f$. In particular, we allow for arbitrarily large growth in the $y$-variable (but still requiring that $f$ does not become degenerate as $|y|\to\infty$). Our main result is the first of its kind that applies to integrands which are unbounded in $y$ and reads as follows:

\begin{alphtheorem}\label{L1lscthm}
\RED{Let $\Omega\subset\R^d$ ($d \geq 2$) be a bounded Lipschitz domain} and let $f\colon\overline{\Omega}\times\R^m\times\R^{m\times d}\to[0,\infty)$ be such that
\begin{enumerate}[(i)]
\item\label{l1hypoth1} $f$ is a Carath\'eodory function whose (strong) recession function $f^\infty$ exists in the sense of Definition~\ref{defrecessionfunction}.
\item\label{l1hypoth2} there exists a continuous function $g\colon\overline{\Omega}\times\R^m\to[0,\infty)$ such that
	\begin{enumerate}[(a)]
	\item\label{eql1intro2} $f$ satisfies a growth bound of the form 
	\[	
	g(x,y)|A|\leq f(x,y,A)\leq Cg(x,y)(1+|A|)
	\]
	for some $C>0$ and for all $(x,y,A)\in\overline{\Omega}\times\R^m\times\R^{m\times d}$;
	\item\label{eql1intro1} there exist $R>0$ and $C>1$ for which
	\[
	g(x,y)\leq C g(x,ty)
	\]
	for all $x\in\Omega$, $|y|\geq R$ and $t\geq 1$;
	\item\label{eql1intro3}  for every compact $K\subset\R^m$ and $\varepsilon>0$, there exists $R_\varepsilon>0$ such that
	\[
	|(f-f^\infty)(x,y,A)|\leq\varepsilon g(x,y)(1+|A|)
	\]
	for all $(x,y)\in\overline{\Omega}\times K$ and $A\in\R^{m\times d}$ with $|A|\geq R_\varepsilon$.
	\end{enumerate}
	\item $f(x,y,\frarg)$ is quasiconvex for every $(x,y)\in\overline{\Omega}\times\R^m$.
\end{enumerate}
For a fixed $g\in\C(\overline{\Omega}\times\R^m\colon[0,\infty))$ verifying conditions~\eqref{eql1intro2},~\eqref{eql1intro1}, and~\eqref{eql1intro3}, define
\[
  \Gcal:=\biggl\{u\in\Lp^1(\Omega;\R^m)\colon\int_\Omega g(x,u(x))\;\dd x<\infty\biggr\}.
\]
Then, the \textcolor{black}{restricted} $\Lp^1(\Omega;\R^m)$-relaxation
\RED{\[
\Fcalro[u]:=\inf\left\{\liminf_{j\to\infty}\Fcal[u_j]\colon\;(u_j)_j\subset\W^{1,1}(\Omega;\R^m)\cap\Gcal\text{ and }u_j\to u\in\BV(\Omega;\R^m)\text{ in $\Lp^1$}\right\},
\]}
\RED{where $u \in \BV(\Omega;\R^m)\cap\Gcal$,} of $\Fcal$ from $\W^{1,1}(\Omega;\R^m)\cap\Gcal$ to $\BV(\Omega;\R^m)\cap\Gcal$ is given by
\begin{align*}
\Fcalro[u]&=\int_\Omega f(x,u(x),\nabla u(x))\;\dd x+\int_\Omega f^\infty\left(x,u(x),\frac{\dd D^c u}{\dd|D^c u|}(x)\right)\;\dd|D^cu|(x)\\
&\qquad+\int_{\Jcal_u} H_f[u](x)\;\dd\Hcal^{d-1}(x).
\end{align*}
\end{alphtheorem}
The energy density $H_f$ which features here is defined similarly to $K_f$ (see~\eqref{eq:Hdef} in Section~\ref{secpreliminaries}) but, in the absence of any $(x,y)$-localisation assumptions on $f$, is in general strictly less than $K_f$ (see Example~\ref{exdensitiesnotthesame}). Hypothesis~\eqref{eql1intro3} is a technical requirement, versions of which feature in the previous works~\cite{FonMul93RQFB,BoFoMa98GMR,FonLeo01LSC}, and which we make use of to compute a lower bound for $\Fcalro$ over $\Jcal_u$ and in the construction of recovery sequences for $\Fcalro$. It is satisfied by all integrands of the form $f(x,y,A)=g(x,y)h(x,y,A)$ where $h^\infty$ exists. We remark that \textcolor{black}{while the restricted $\Lp^1(\Omega;\R^m)$-relaxation considered here is in general different from the one defined in~\eqref{eq:firstFcalro},} the restriction of $\Fcalro$ to the class $\Gcal$ causes no issues in applications coming from~\eqref{eqperturbfamily}, where we expect $g(x,u(x)) < \infty$ almost everywhere, or if $|g(x,y)|\leq C(1+|y|^{d/(d-1)})$. For some further discussion of this requirement see Remark~\ref{rem:G}.

Theorem~\ref{L1lscthm} assumes that $f^\infty$ exists in a stronger sense than has been classically required in the literature (see Definition~\ref{defrecessionfunction}), where only the upper recession function $f^\#$ is used. In fact, the other properties required of $f$ in~\cite{FonLeo01LSC} imply that their $f^\#$ must exist in the sense of Definition~\ref{defrecessionfunction} at every point of continuity for $f^\#$, that $f^\#$ must be lower semicontinuous, and such that $f^\#(x,y,A)$ is continuous in $(x,A)$ for every $y\in\R^m$.

This paper is structured as follows: Necessary concepts and preliminaries are introduced in Section~\ref{secpreliminaries}. Section~\ref{secrelaxationcomparisons} contains a discussion of the relationship between Theorem~\ref{L1lscthm} and other results in the existing literature. The proof of Theorem~\ref{L1lscthm} is then carried out in two parts. We first establish the lower bound
\begin{align}
\begin{split}\label{eqlscstatement}
\RED{\Fcalro}[u] &\geq \int_\Omega f(x,u(x),\nabla u(x))\;\dd x+\int_\Omega f^\infty\left(x,u(x),\frac{\dd D^cu}{\dd|D^c u|}(x)\right)\;\dd|D^cu|(x) \\
  &\qquad +\int_{\mathcal{J}_u}\RED{H_f}[u](x)\;\dd\mathcal{H}^{d-1}(x)  
 \end{split}
\end{align}
by employing an improved version of the blow-up technique pioneered in~\cite{FonMul92QCIL} and~\cite{FonMul93RQFB}: Letting \RED{$(u_j)_j\subset(\C^\infty\cap\W^{1,1})(\Omega;\R^m)$} be such that
\[
 u_j\to u\text{ in }\Lp^1(\Omega;\R^m)\quad\text{ and }\quad\lim_{j\to\infty}\Fcal[u_j]=\Fcalro[u]
\]
(the existence of such a sequence $(u_j)_j$ follows from the definition of $\Fcalro$ together with a diagonal argument), we can pass to a non-relabelled subsequence in order to assume that there exists a Radon measure $\mu\in\mbfM(\overline{\Omega})$ such that
\begin{equation*}
\wstarlim_{j\to\infty}f(x,u_j(x),\nabla u_j(x))\lL\restrict\Omega(\dd x)=\mu\quad\text{ in }\mbfM(\overline{\Omega}).
\end{equation*}
Using the Radon--Nikod\'{y}m Differentiation Theorem, we can write $\mu$ as the sum of mutually singular measures,
\[
\mu=\frac{\dd\mu}{\dd\lL}\lL\restrict\Omega+\frac{\dd\mu}{\dd|D^c u|}|D^c u|+\frac{\dd\mu}{\dd\Hcal^{d-1}\restrict\Jcal_u}\Hcal^{d-1}\restrict\Jcal_u+\mu^s.
\]
To obtain~\eqref{eqlscstatement}, it therefore suffices to prove the three pointwise inequalities
\begin{align*}
&\frac{\dd\mu}{\dd\lL}(x)\geq f(x,u(x),\nabla u(x))\quad\text{ for }\lL\text{-almost every }x\in\Omega,\\
&\frac{\dd\mu}{\dd|D^c u|}(x)\geq f^\infty\left(x,u(x),\frac{\dd D^cu}{\dd|D^cu|}(x)\right)\quad\text{ for }|D^cu|\text{-almost every }x\in\Omega,\\
&\frac{\dd\mu}{\dd\Hcal^{d-1}\restrict\Jcal_u}(x)\geq H_f[u](x)\quad\text{ for }\Hcal^{d-1}\text{-almost every }x\in\Jcal_u.
\end{align*}
The first two of these inequalities are proved simultaneously in Section~\ref{chapl1lsc} via an intricate measure-theoretic truncation argument, which allows us to replace the $\Lp^1(\Omega;\R^m)$-convergent sequence $(u_j)_j$ with a weakly* convergent sequence $(\widetilde{u}_j)_j$ and to reduce the problem to an application of the weak* lower semicontinuity theory developed in~\cite{RindlerShaw19}. We note that the work in~\cite{RindlerShaw19} features in two distinct ways here: the main weak* relaxation result (quoted below as Theorem~\ref{thmwsclsc}) is used to obtain the final inequality above, but the theory of liftings together with an associated Besicovitch Differentiation Theorem is also used in an essential way to control the error term which arises when exchanging $(u_j)_j$ for $(\widetilde{u}_j)_j$. Whilst previous works have only been able to use measure-theoretic techniques to localise at points $x\in\Omega$, the key point here is that liftings allow us to use measure theory to localise at points $(x,u(x))\in\Omega\times\R^m$ with a single limit. Their use here is a new technique which we hope will also find applications elsewhere. As a consequence of the localisation made possible by liftings, all of the results proved in this section hold for all Carath\'eodory integrands which are partially coercive and such that $f^\infty$ exists, without needing to assume either hypothesis~\eqref{eql1intro2} or~\eqref{eql1intro3} of Theorem~\ref{L1lscthm}.

Section~\ref{secl1jumplsc} is devoted to obtaining the optimal lower bound for $\frac{\dd\mu}{\dd\Hcal^{d-1}\restrict\Jcal_u}$. We are not able to reduce the problem here to an application of the weak* theory and hence must assume more of the integrand $f$ to proceed. The key technical result is Lemma~\ref{propfphi}, which shows that \RED{we may approximate the functional $\Fcal$ with functionals whose integrands are better behaved}.

Finally, in Section~\ref{chaprecoveryseqs} we deal with the second component of the proof of Theorem~\ref{L1lscthm}: we prove the existence of \emph{recovery sequences} \RED{$(u_j)_j\subset(\C^\infty\cap\W^{1,1})(\Omega;\R^m)$} with the property that $u_j\to u$ in $\Lp^1(\Omega;\R^m)$ and
\begin{align*}
\lim_{j\to\infty}\Fcal[u_j]=&\int_\Omega f(x,u(x),\nabla u(x))\;\dd x+\int_\Omega f^\infty\left(x,u(x),\frac{\dd D^c u}{\dd|D^c u|}(x)\right)\;\dd|D^c u|(x)\\
&\qquad+\int_{\mathcal{J}_u}H_f[u](x)\;\dd\mathcal{H}^{d-1}(x).
\end{align*}
These are explicitly constructed using a version of the technique developed in~\cite{RindlerShaw19}.

\section*{Acknowledgements}
The authors would like to thank Irene Fonseca, Jan Kristensen and Neshan Wickramasekera for several helpful discussions. \RED{We are also grateful to the referees for their thorough reading of the manuscript, which led to many improvements.}

This project has received funding from the European Research Council (ERC) under the European Union's Horizon 2020 research and innovation programme, grant agreement No 757254 (SINGULARITY). F.~R.\ also acknowledges the support from an EPSRC Research Fellowship on Singularities in Nonlinear PDEs (EP/L018934/1).

The work presented in this paper forms part of G.S.'s PhD thesis, undertaken at the Cambridge Centre for Analysis at the University of Cambridge, and supported by the UK Engineering and Physical Sciences Research Council (EPSRC) grant EP/H023348/1 for the University of Cambridge Centre for Doctoral Training, the Cambridge Centre for Analysis, whose support G.S. gratefully acknowledges.

\section{Preliminaries}\label{secpreliminaries}
Throughout this work, $\Omega\subset\R^d$ will always be assumed to be a bounded open domain with compact Lipschitz boundary $\pd\Omega$ in dimension $d\geq 2$, and $\bB^k$, $\pd\bB^k$ will denote the open unit ball in $\R^k$ and its boundary (the unit sphere) respectively. The open ball of radius $r$ centred at $x\in\R^k$ is $B(x,r)$, although we will sometimes write $B^k(x,r)$ if the dimension of the ambient space needs to be emphasised for clarity. The volume of the unit ball in $\R^k$ will be denoted by $\omega_k:=\Lcal^k(\bB^k)$, where $\Lcal^k$ is the usual $k$-dimensional Lebesgue measure. We will write $\R^{m\times d}$ for the space of $m\times d$ real-valued matrices. The map $\pi\colon\Omega\times\R^m\to\Omega$ denotes the projection $\pi((x,y)):=x$, and $T^{(x_0,r)}\colon\R^d\to\R^d$, $T^{(x_0,r),(y_0,s)}\colon\R^d\times\R^m\to\R^d\times\R^m$ represent the homotheties $x\mapsto (x-x_0)/r$ and $(x,y)\mapsto((x-x_0)/r,(y-y_0)/s))$, respectively. Tensor products $a\otimes b\in\R^{m\times d}$ and $f\otimes g$ for vectors $a\in\R^m$, $b\in\R^d$, and real-valued functions $f$, $g$, are defined componentwise by $(a\otimes b)_{i,j}:=a_ib_j$ and $(f\otimes g)(x,y):=f(x)g(y)$ respectively.

The closed subspaces of $\BV(\Omega;\R^m)$ and $\C^\infty(\Omega;\R^m)$ consisting only of the functions satisfying $(u):=\dashint_\Omega u(x)\;\dd x=0$ are denoted by $\BV_\#(\Omega;\R^m)$ and $\C^\infty_\#(\Omega;\R^m)$ respectively. We shall use the notation $(u)_\Omega$ when the domain of integration might not be clear from context, as well as the abbreviation $(u)_{x,r}:=(u)_{B(x,r)}$.  We shall sometimes use subscripts for clarity when taking the gradient with respect to a partial set of variables: that is, if $f=f(x,y)\in\C^1(\Omega\times\R^m)$ then $\nabla_x f:=(\pd_{x_1}f,\pd_{x_2}f,\ldots,\pd_{x_d}f)$ and $\nabla_yf:=(\pd_{y_1}f,\pd_{y_2}f,\ldots,\pd_{y_m}f)$.

\subsection{Measure theory}\label{subsecmeasuretheory}
For a separable locally convex metric space $X$ (in our case, $X$ will usually be either $\Omega$ or $\Omega\times\R^m$), the space of vector-valued Radon measures on $X$ taking values in a finite-dimensional normed vector space $V$ (usually $\R^{m\times d})$ will be written as $\mbfM(X;V)$ or just $\mbfM(X)$ if $V=\R$. The cone of positive Radon measures on $X$ is $\mbfM^+(X)$, and the set of elements $\mu\in\mbfM(X;V)$ whose total variation $|\mu|$ is a probability measure, is $\mbfM^1(X;V)$. The notation $\mu_j\wsc\mu$ will denote the usual weak* convergence of measures, i.e.\ $\ip{\varphi}{\mu_j}\to\ip{\varphi}{\mu}$ for all $\varphi\in\C_0(X;\RED{V^*})$, \RED{where $V^*$ denotes the dual space to $V$.} We recall that $\mu_j$ is said to converge to $\mu$ \textbf{strictly} if $\mu_j\wsc \mu$ and in addition $|\mu_j|(X)\to|\mu|(X)$. If $\mu_j\to\mu$ strictly, then we have that $\ip{\varphi}{\mu_j}\to\ip{\varphi}{\mu}$ for all $\varphi\in\C_b(X;V^*)$ rather than just for all $\varphi\in\C_0(X;V^*)$. Given a map $T$ from $X$ to another separable, locally convex metric space $Y$, the pushforward operator $T_\#\colon\Mbf(X;V)\to\mbfM(Y;V)$ is defined by
\[
\ip{\varphi}{T_\#\mu}:=\ip{\varphi\circ T}{\mu},\qquad\varphi\in\C_0(Y;V^*).
\]
If $T$ is continuous and proper, then $T_\#$ is continuous when $\mbfM(X;V)$ and $\mbfM(Y;V)$ are equipped with their respective weak* or strict topologies. 

For $k\in[0,\infty)$, the $k$-dimensional Hausdorff (outer) measure on $\R^d$ is written as $\Hcal^{k}$ and, if $A\in\Bcal(\R^d)$ is a Borel set satisfying $\Hcal^k(A)<\infty$, its restriction $\Hcal^k\restrict A$ to $A$ defined by $[\Hcal^k\restrict A](\frarg):=\Hcal^k(\frarg\cap A)$ is a finite Radon measure. A set $A\subset\R^m$ is said to be \textbf{countably $\boldsymbol{\Hcal^k}$-rectifiable} if there exists a sequence of Lipschitz functions $f_i\colon\R^k\to\R^d$ ($i \in \N$) such that
\[
\Hcal^k\left(A\setminus\bigcup_{i=1}^\infty f_i(\R^k)\right)=0,
\]
and \textbf{$\boldsymbol{\Hcal^k}$-rectifiable} if in addition $\Hcal^k(A)<\infty$.
We say that $\mu\in\mbfM(\R^d;V)$ is a \textbf{$\boldsymbol k$-rectifiable measure} if there exists a countably $\Hcal^k$-rectifiable set $A\subset\R^d$ and a Borel function $f\colon A\to V$ such that $\mu=f\Hcal^k\restrict A$.

With $A$ assumed to be countably $\Hcal^k$-rectifiable, we can compute the Radon--Nikod\'{y}m derivative for any $\mu\in\mbfM(\R^d)$ with respect to $\Hcal^k\restrict A$, given for $\Hcal^k$-almost every $x\in A$, by
\begin{equation*}
\frac{\dd\mu}{\dd\Hcal^k\restrict A}(x):=\lim_{r\to 0}\frac{\mu(B(x,r))}{\omega_{k}r^k}.
\end{equation*}
The function $\frac{\dd\mu}{\dd\Hcal^k\restrict A}$ is a Radon--Nikod\'{y}m derivative in the sense that $\frac{\dd\mu}{\dd\Hcal^k\restrict A}\Hcal^k\restrict A$ is a $k$-rectifiable measure and that we can decompose
\[
\mu=\frac{\dd\mu}{\dd\Hcal^k\restrict A}\Hcal^k\restrict A+\mu^s,\quad\text{ where $\mu^s$ satisfies }\quad\frac{\dd\mu}{\dd\Hcal^k\restrict A}\Hcal^k\restrict A\perp \mu^s,
\]
in analogy with the usual Lebesgue--Radon--Nikod\'{y}m decomposition. \RED{Here, $\mu \perp \nu$ for measures $\mu, \nu \in \mbfM(\R^d;V)$ means that $\mu,\nu$ are mutually singular.}

\RED{For $k\in\mbN$, a} measure $\mu\in\mbfM(\R^d;V)$ is said to admit a ($k$-dimensional) \textbf{approximate tangent space} at $x_0$ if there exists an (unoriented) $k$-dimensional hyperplane $\tau\subset\R^d$ and $\theta\in V$ such that
\[
r^{-k} T_\#^{(x_0,r)}\mu\to\theta\Hcal^k\restrict(\overline{\bB^d}\cap\tau) \quad\text{ strictly in }\mbfM(\overline{\bB^d};V)\quad\text{ as }r\to 0.
\]
As a consequence of the Lebesgue Differentiation Theorem, we can see that $f\mu$ admits an approximate tangent space at $x\in\R^d$ whenever $x$ is a $\mu$-Lebesgue point of $f$ and a point at which $\mu$ admits an approximate tangent space.
The existence of approximate tangent spaces characterises the class of rectifiable measures in the sense that $\mu\in\mbfM(\R^d;V)$ possesses a $k$-dimensional approximate tangent space at $|\mu|$-almost every $x_0\in\R^d$ if and only if $\mu$ is $k$-rectifiable (see Theorem~2.83 in~\cite{AmFuPa00FBVF}). 

We shall make use of the Vitali--Besicovitch Covering and Besicovitch Differentiation Theorems for Radon measures in $\R^d$ (see Theorems~2.19 and~2.22 in~\cite{AmFuPa00FBVF}):

\begin{theorem}[Vitali--Besicovitch Covering Theorem] \label{thm:besi}
Let $A\subset\R^d$ be a bounded Borel set. A collection $\Fcal$ of closed balls in $\R^d$ is said to be a \textbf{fine cover} for $A$ if, for every $x\in A$ and $R>0$, there exists $r<R$ such that $\overline{B(x,r)}\in\Fcal$. 

If $\Fcal$ is a fine cover for $A$ then, for every $\mu\in\mbfM^+(\R^d)$ there exists a disjoint countable family $\Fcal'\subset\Fcal$ such that
\[
\mu\left(A\setminus\RED{\bigcup_{\overline{B}\in\Fcal'}\overline{B}}\right)=0.
\]
\end{theorem}

\begin{theorem}[Besicovitch Differentiation Theorem]
Let $\mu\in\mbfM^+(\Omega)$ and $\eta\in\mbfM(\Omega;V)$. Then, for $\mu$-almost every $x\in\Omega$, the limit
\[
f(x):=\lim_{r\to 0}\frac{\eta(B(x,r))}{\mu(B(x,r))}=\lim_{r\to 0}\frac{\eta(\overline{B(x,r))}}{\mu(\overline{B(x,r))}}
\]
exists and is equal to the Radon--Nikod\'{y}m derivative $\frac{\dd\eta}{\dd\mu}(x)$.
\end{theorem}

In addition to the usual version of the Besicovitch Differentiation Theorem, we shall make use of a new generalised version, first proved in~\cite{RindlerShaw19}, which applies to measures that behave like graphs. The following two results are taken from Theorem~5.1 and Lemma~5.2 in~\cite{RindlerShaw19}. \RED{We note that graph-like measures have been used before in variational problems, for example in~\cite{DMas79IRBV} and~\cite{AviGig91VIMB,AviGig92MCGR}}.

Given a function $u\colon\Omega\to\R^m$, we denote its associated \textbf{graph map} by $\gr^u\colon x\mapsto (x,u(x))$. If $u$ is measurable with respect to a measure $\mu$, then the pushforward measure $\gr^u_\#\mu$ is well-defined as a measure on $\Omega\times\R^m$. We say that $\eta\in\mbfM(\Omega\times\R^m)$ (or $\mbfM(\Omega\times\R^m;\R^{m\times d})$) is a $u$-\textbf{graphical measure} (or just a graphical measure) if it arises in this way. Note that if $\eta$ is a $u$-graphical measure then it must be the case that $\eta=\gr^u_\#(\pi_\#\eta)$.

We will use the following theorem from~\cite{RindlerShaw19}\footnote{\RED{The statement of this theorem in~\cite{RindlerShaw19} is incorrect insofar as the computation of $\frac{\dd\lambda}{\dd\eta}(x,u(x))$ when $\lambda$ is \emph{not} singular with respect to $\eta$ requires the additional assumption that $\lim_{r\to 0}\frac{\eta\left(\overline{B^d(x,r)\times B(y_r^x,c_r^x)}\right)}{\mu\left(\overline{B^d(x,r)}\right)} = 1$; this assertion for $\lambda$ \emph{not} singular to $\eta$, however, is used neither in~\cite{RindlerShaw19} or here. For the sake of clarity we provide a full proof.}}.

\RED{\begin{theorem}[Generalised Besicovitch Differentiation Theorem for graphical measures]\label{thmgeneralisedbesicovitch}
Given $u\colon\Omega\to\R^m$, let $\eta=\gr^u_\#\mu\in\mbfM^+(\Omega\times\R^m)$ be a $u$-graphical measure. For each $x\in\Omega$, let $(c_r^x)_{r > 0} \subset (0,\infty)$ satisfy $\lim_{r\downarrow 0}c_r^x=0$ and $(y_r^x)_{r > 0} \subset \R^m$ satisfy $\lim_{r\downarrow 0}y_r^x=u(x)$. Then:

\begin{enumerate}[(i)]
\item\label{graphbesicovitchsingular} If $\lambda$ is a (possibly vector-valued) measure on $\Omega\times\R^m$ satisfying $\lambda\perp\eta$, we have that 
\[
0=\frac{\dd\lambda}{\dd\eta}(x,u(x))=\lim_{r\to 0}\frac{\lambda\left(\overline{B(x,r)\times B(y_r^x,c_r^x)}\right)}{\pi_\#\eta(\overline{B(x,r)})}=\lim_{r\to 0}\frac{\lambda\left(\overline{B(x,r)\times B(y_r^x,c_r^x)}\right)}{\mu(\overline{B(x,r)})}
\]
for $\eta$-almost every $(x,u(x))\in\Omega\times\R^m$, where $\pi\colon\Omega\times\R^m\to\Omega$ is the projection map $\pi((x,y)):=x$.
\item\label{cylindricallebesgue} A cylindrical version of the Lebesgue Differentiation Theorem holds in the sense that, for any $f\in\Lp^1(\Omega\times\R^m,\eta)$ (i.e., \textcolor{black}{the $\Lp^1$-space with respect to $\eta=\gr^u_\#\mu$}),
\[
\lim_{r\to 0}\frac{1}{\mu(\overline{B(x,r)})}\int_{\overline{B(x,r)}\times \R^m}|f(\overline{x},y)-f(x,u(x))|\;\dd\eta(\overline{x},y)=0
\]
for $\mu$-almost every $x\in\Omega$.
\item\label{graphbesicovitchgeneral} If $\lambda\in\mbfM(\Omega\times\R^m)$, then
\[
\frac{\dd\lambda}{\dd\eta}(x,u(x))=\lim_{r\to 0}\frac{\lambda\left(\overline{B(x,r)\times B(y_r^x,c_r^x)}\right)}{\pi_\#\eta(\overline{B(x,r)})}=\lim_{r\to 0}\frac{\lambda\left(\overline{B(x,r)\times B(y_r^x,c_r^x)}\right)}{\mu(\overline{B(x,r)})}
\]
for $\eta$-almost every $(x,u(x))\in\Omega\times\R^m$ for which $(c_r^x)_r$ and $(y_r^x)_r$ are such that $\lim_{r\to 0}\frac{\eta\left(\overline{B(x,r)\times B(y_r^x,c_r^x)}\right)}{\mu\left(\overline{B(x,r)}\right)} = 1$.
\end{enumerate}
 \end{theorem}

\begin{proof}
We first establish~\eqref{graphbesicovitchsingular}: let $\lambda$ be a measure on $\Omega\times\R^m$ satisfying $\lambda\perp\eta$ and define for $x \in \Omega$ such that $\mu(\overline{B(x,r)}) > 0$ for all $r > 0$ the function
\[
F(x):=\limsup_{r\to 0}\frac{|\lambda|(\overline{B(x,r)\times  B(y_r^x,c_r^x)})}{\mu(\overline{B(x,r)})}.
\]
Let $Z \subset \Omega$ be such that $\mu(Z) = \mu(\Omega)$, $|\lambda|(\gr^u(Z))=0$ and such that at every $x \in Z$ it holds that $\mu(\overline{B(x,r)}) > 0$ for all $r > 0$. Such a set $Z$ exists since $\eta = \gr^u_\# \mu$ is singular to $\lambda$ and $\eta$ is \textcolor{black}{concentrated on} the set $\gr^u(\Omega)$.

Let $E\subset Z$ be a Borel set such that $F(x)>t>0$ for all $x\in E$ and for $\varepsilon\in(0,t)$ arbitrary let $A\supset \gr^u(E)$ be an open set. Define the family
\begin{align*}
\mathcal{F}:=&\Bigl\{\overline{B(x,r)\times B(y_r^x,c_r^x)}\colon\; x\in E,\text{ }\overline{B(x,r)\times B(y_r^x,{c_r^x})}\subset A\text{ and }\\
&\qquad|\lambda|(\overline{B(x,r)\times B(y_r^x,c_r^x)})\geq(t-\varepsilon)\mu(\overline{B(x,r)})\Bigr\}\subset\R^d\times\R^m.
\end{align*}
Since $c^x_r\downarrow 0$ and $y_r^x\to u(x)$ as $r\downarrow 0$ it follows that $\pi_\#\mathcal{F}$ is a fine cover for $E$ and so, by the Vitali--Besicovitch Covering Theorem~\ref{thm:besi}, there exists a countable disjoint subfamily of $\pi_\#\mathcal{F}$, which we write as $\pi_\#\mathcal{F}'$ for some disjoint subcollection $\mathcal{F}'\subset\mathcal{F}$, that covers $\mu$-almost all of $E$. We therefore have that
\[
(t-\varepsilon)\mu(E)\leq(t-\varepsilon)\sum_{B\in\mathcal{F}'}\mu(\pi_\# B)\leq\sum_{B\in\mathcal{F}'}|\lambda|(B)\leq|\lambda|(A).
\]
First letting $\varepsilon\downarrow 0$ and then using the outer regularity of $|\lambda|$ to approximate $\gr^u(E)$ with a sequence of open sets, we obtain
\[
t\mu(E)\leq|\lambda|(\gr^u(E)).
\]
Since $|\lambda|(\gr^u(Z))=0$, also $|\lambda|(\gr^u(E))=0$ and if $E$ was such that $\eta(\gr^u(E))=\mu(E)>0$, then $t=0$, a contradiction. It follows that $F(x)=0$ for $\mu$-almost every $x\in\Omega$ and hence that
\[
\lim_{r\to 0}\frac{\lambda(\overline{B(x,r)\times B(y_r^x,c_r^x)})}{\mu(\overline{B(x,r)})}=\frac{\dd\lambda}{\dd\eta}(x,u(x))=0 \quad \text{ for }\eta\text{-a.e.\ }(x,u(x))\in\Omega\times\R^m,
\]
as required. 

Next, we prove~\eqref{cylindricallebesgue}: for $f\in\Lp^1(\Omega\times\R^m,\eta)$, the definition of $\eta=\gr^u_\#\mu$ implies
\begin{align*}
&\lim_{r\to 0}\frac{1}{\mu(\overline{B(x,r)})}\int_{\overline{B(x,r)}\times \R^m}|f(\overline{x},y)-f(x,u(x))|\;\dd\eta(\overline{x},y)\\
&\qquad =\lim_{r\to 0}\frac{1}{\mu(\overline{B(x,r)})}\int_{\overline{B(x,r)}}|f(\overline{x},u(\overline{x}))-f(x,u(x))|\;\dd\mu(\overline{x}).
\end{align*}
It can be checked that the conditions $f\in\Lp^1(\Omega\times\R^m,\eta)$ and $f\circ \gr^u\in\Lp^1(\Omega,\mu)$ are equivalent, and so it follows from $f\in\Lp^1(\Omega\times\R^m,\eta)$ that $\mu$-almost every $x\in\Omega$ is a $\mu$-Lebesgue point for $x\mapsto f(x,u(x))$. We therefore deduce that
\begin{equation}\label{eqlebesguepoint}
\lim_{r\to 0}\frac{1}{\mu(\overline{B(x,r)})}\int_{\overline{B(x,r)}\times \R^m}|f(\overline{x},y)-f(x,u(x))|\;\dd\eta(\overline{x},y)=0
\end{equation}
for $\mu$-almost every $x\in\Omega$, as required.

Finally, to prove~\eqref{graphbesicovitchgeneral} we argue as follows: for $\lambda\in\mbfM(\Omega\times\R^m)$ let
\[
\lambda = \frac{\dd\lambda}{\dd\eta}\eta + \lambda^s, \qquad\lambda\perp\eta
\]
be the usual Radon--Nikod\'{y}m decomposition of $\lambda$ with respect to $\eta$. Noting that~\eqref{cylindricallebesgue} implies
\begin{align*}
&\lim_{r\to 0}\frac{1}{\mu(\overline{B(x,r)})}\int_{\overline{B(x,r)\times B(y_r^x,c_r^x)}}\left|\frac{\dd\lambda}{\dd\eta}(\overline{x},y)-\frac{\dd\lambda}{\dd\eta}(x,u(x))\right|\;\dd\eta(\overline{x},y)\\
&\qquad \leq \lim_{r\to 0}\frac{1}{\mu(\overline{B(x,r)})}\int_{\overline{B(x,r)}\times \R^m}\left|\frac{\dd\lambda}{\dd\eta}(\overline{x},y)-\frac{\dd\lambda}{\dd\eta}(x,u(x))\right|\;\dd\eta(\overline{x},y)\\
&\qquad =0,
\end{align*}
we can use~\eqref{graphbesicovitchsingular} and~\eqref{cylindricallebesgue} together to see that
\begin{align*}
\lim_{r\to 0}\frac{\lambda\left(\overline{B(x,r)\times B(y_r^x,c_r^x)}\right)}{\mu(\overline{B(x,r)})} & = \lim_{r\to 0}\frac{1}{\mu(\overline{B(x,r)})}\int_{\overline{B(x,r)\times B(y_r^x,c_r^x)}}\frac{\dd\lambda}{\dd\eta}(\overline{x},y)\;\dd\eta(\overline{x},y) \\
&\qquad+ \lim_{r\to 0}\frac{\lambda^s\left(\overline{B(x,r)\times B(y_r^x,c_r^x)}\right)}{\mu(\overline{B(x,r)})}\\
& = \frac{\dd\lambda}{\dd\eta}(x,u(x)) \cdot \lim_{r\to 0}\frac{\eta\left(\overline{B(x,r)\times B(y_r^x,c_r^x)}\right)}{\mu(\overline{B(x,r)})} + 0
\end{align*}
for $\eta$-almost every $(x,u(x))\in\Omega\times\R^m$, at which point the conclusion follows.
\end{proof}}

Using Theorem~\ref{thmgeneralisedbesicovitch}, it can be proved that the behaviour of graphical measures under general homotheties is stable under multiplication by integrable functions:

\begin{lemma}\label{lemgraphicallebesgueconverg}
Let $\eta$ be a $u$-graphical measure on $\Omega\times\R^m$ (that is, $\eta=\gr^u(\pi_\#\eta)$ for some $\pi_\#|\eta|$-measurable function $u\colon\Omega\to\R^m$) and let $x_0\in\Omega$, $r_n\downarrow 0,c_n\downarrow 0$, \RED{$(a_n)_n$} and $(y_n)_n\subset\R^m$ with $y_n\to u(x_0)$ be such that 
\[
a_nT_\#^{(x_0,r_n),(y_n,c_n)}\eta\wsc\eta^0\qquad\text{ in }\mbfM(\overline{\bB^d}\times\R^m).
\]
If $f\in\Lp^1(\Omega\times\R^m,\mu)$ and $x_0$ is an $\eta$-cylindrical Lebesgue point for $f$ in the sense of Theorem~\ref{thmgeneralisedbesicovitch}, then
\[
a_nT_\#^{(x_0,r_n),(y_n,c_n)}(f\eta)\wsc f(x_0,u(x_0))\eta^0.
\]
Moreover, if $a_nT_\#^{(x_0,r_n),(y_n,c_n)}\eta\to\eta^0$ strictly, then $a_nT_\#^{(x_0,r_n),(y_n,c_n)}(f\eta)\to f(x_0)\eta^0$ strictly as well.
\end{lemma}

\subsection{BV-functions}
Given a function $u\in\BV(\Omega;\R^m)$, we recall the mutually singular decomposition $Du=\nabla\lL\restrict\Omega +D^cu +D^ju$ of the derivative $Du$, where $|D^cu|\ll \Hcal^{d-1}$, $D^ju$ is absolutely continuous with respect to $\Hcal^{d-1}\restrict\Jcal_u$, and $\Jcal_u$ is the countably $\Hcal^{d-1}$-rectifiable jump set of $u$. Each $u\in\BV(\Omega;\R^m)$ admits a precise representative $\widetilde{u}\colon\Omega\to\R^m$ which is defined $\Hcal^{d-1}$-almost everywhere in $\Omega\setminus\Jcal_u$. The \textbf{jump interpolant} associated to $u$ is then the function $u^\theta\colon\Omega\times[0,1]\to\R^m$ defined, up to a choice of orientation $n_u$ for the jump set $\Jcal_u$ of $u$, for $\Hcal^{d-1}$-almost every $x\in\Omega$ by
\begin{equation}\label{eqdefjumpinterpolant}
u^\theta(x):=\begin{cases}
\theta u^-(x)+(1-\theta)u^+(x) &\text{if }x\in\Jcal_u,\\
\widetilde{u}(x) &\text{otherwise,}
\end{cases}
\end{equation}
\RED{where $u^-(x)$ and $u^+(x)$ are the one-sided jump limits of $u$ at $x$ in negative and positive $n_u(x)$-direction, respectively.} The need to fix a choice of orientation for $\Jcal_u$ in order to properly define $u^\theta$ is obviated by the fact that $u^\theta$ will only appear in expressions of the form
\[
\int_0^1 \varphi(u^\theta(x))\;\dd\theta,
\]
which are invariant of our choice of $n_u$.

Given a function $u\in\BV(\Omega;\R^m)$, if $\mu$ is a measure on $\Omega$ satisfying both $|\mu|\ll\Hcal^{d-1}$ and $|\mu|(\Jcal_u)=0$ (we will usually take $\mu=|Du|\restrict(\Omega\setminus\Jcal_u)$), then the precise representative of $u$ (which we also denote by $u$) is $\mu$-measurable and so the pushfoward $\gr^u_\#\mu$ is well-defined as a $u$-graphical measure on $\Omega\times\R^m$.

A sequence $(u_j)_j\subset\BV(\Omega;\R^m)$ is said to converge \textbf{strictly} to $u\in\BV(\Omega;\R^m)$ if $u_j\to u$ in $\Lp^1(\Omega;\R^m)$ and $Du_j\to Du$ strictly in $\mbfM(\Omega;\R^{m\times d})$ as $j\to\infty$. We say that $u_j$ converges \textbf{area-strictly} to $u$ if $u_j\to u$ in $\Lp^1(\Omega;\R^m)$ and in addition
\[
\int_\Omega\sqrt{1+|\nabla u_j(x)|^2}\;\dd x+|D^s u_j|(\Omega)\to\int_\Omega\sqrt{1+|\nabla u(x)|^2}\;\dd x+|D^s u|(\Omega)
\]
as $j\to\infty$. It is the case that area-strict convergence implies strict convergence in $\BV(\Omega;\R^m)$ and that strict convergence implies weak* convergence. That none of these notions of convergence coincide follows from considering the sequence $(u_j)_j\subset\BV((-1,1))$ given by $u_j(x):=x+(a/j)\sin(jx)$ for some $a\neq 0$ fixed. This sequence converges weakly* to the function $x\mapsto x$, strictly if and only if $|a|\leq 1$, but (since the function $z\mapsto\sqrt{1+|z|^2}$ is strictly convex away from $0$) never area-strictly. Smooth functions are area-strictly (and hence strictly) dense in $\BV(\Omega;\R^m)$: indeed, if $u\in\BV(\Omega;\R^m)$ and $(u_\rho)_{\rho>0}$ is a family of radially symmetric mollifications of $u$ then it holds that $u_\rho\to u$ area-strictly as $\rho \todown 0$ (see Lemma~1 in~\cite{KriRin10CGYM}).

If $\Omega\subset\R^d$ is such that $\pd\Omega$ is Lipschitz and compact, then the trace onto $\pd\Omega$ of a function $u\in\BV(\Omega;\R^m)$ is denoted by $u|_{\pd\Omega}\in\Lp^1(\pd\Omega;\R^m)$. The trace map $u\mapsto u|_{\pd\Omega}$ is norm-bounded from $\BV(\Omega;\R^m)$ to $\Lp^1(\pd\Omega;\R^m)$ and is continuous with respect to strict convergence (see Theorem~3.88 in~\cite{AmFuPa00FBVF}). If $u,v\in\BV(\Omega;\R^m)$ are such that $u|_{\pd\Omega}=v|_{\pd\Omega}$, then we shall sometimes simply say that ``$u=v$ on $\pd\Omega$''.

The following proposition, a proof for which can be found in the appendix of~\cite{KriRin10CGYM} (or Lemma~B.1 of~\cite{Bild03CVP} in the case of a Lipschitz domain $\Omega$), states that we can even require that smooth area-strictly convergent approximating sequences satisfy the trace equality ${u_j}|_{\pd\Omega}=u|_{\pd\Omega}$:
\begin{proposition}\label{propfixedbdary}
For every $u\in\BV(\Omega;\R^m)$, there exists a sequence $(u_j)_j\subset(\C^\infty\cap\W^{1,1})(\Omega;\R^m)$ with the property that 
\[
u_j\to u\qquad\text{ area-strictly in }\BV(\Omega;\R^m)\text{ as }j\to\infty \qquad\text{ and }\qquad {u_j}|_{\pd\Omega}=u|_{\pd\Omega}.
\]
Moreover, if $u|_{\pd\Omega}\in\Lp^\infty(\pd\Omega;\R^m)$ we can assume that $(u_j)_j\subset(\C^\infty\cap\W^{1,1}\cap\Lp^\infty)(\Omega;\R^m)$  and, if $u\in\Lp^\infty(\Omega;\R^m)$, then we can also require that $\sup_j\norm{u_j}_{\Lp^\infty}\leq\norm{u}_{\Lp^\infty}$.
\end{proposition}

Next, we characterise the behaviour of $\BV(\Omega;\R^m)$-functions under rescaling.
\begin{theorem}[Blowing-up BV-functions]\label{thmbvblowup}
Let $u\in\BV(\Omega;\R^m)$ and write $\Omega$ as the disjoint union
\[
\Omega=\Dcal_u\cup\mathcal{J}_u\cup \Ccal_u\cup\mathcal{N}_u
\]
and $Du$ as the mutually singular sum
\[
Du=\nabla u\lL+D^ju+D^cu,\qquad \nabla u\lL=Du\restrict\Dcal_u,\quad D^ju=Du\restrict\Jcal_u,\quad D^cu=Du\restrict\Ccal_u,
\]
where $\mathcal{D}_u$ denotes the set of points at which $u$ is approximately differentiable, $\mathcal{J}_u$ denotes the set of jump points of $u$, $\mathcal{C}_u$ denotes the set of points where $u$ is approximately continuous but not approximately differentiable, and $\mathcal{N}_u$ satisfies $|Du|(\Ncal_u)=\mathcal{H}^{d-1}(\mathcal{N}_u)=0$. For $r>0$ and $x\in\Dcal_u\cup\Ccal_u\cup\Jcal_u$, define $u^r\in\BV(\bB^d;\R^m)$ by
\[
u^r(z):=c_r\left(\frac{u(x+rz)-(u)_{x,r}}{r}\right),\quad c_r:=\begin{cases}1&\text{ if }x\in\Dcal_u,\\
r&\text{ if }x\in\Jcal_u,\\
\frac{r^d}{|Du|(B(x,r))}&\text{ if }x\in\Ccal_u.
\end{cases}
\]
Then the following trichotomy relative to $\lL\restrict\Omega+|Du|$ holds:
\begin{enumerate}[(i)]
\item For $\lL$-almost every $x\in\Omega$, 
\[
u^r\to \nabla u(x)\frarg \text{ strongly in }\BV(\bB^d;\R^m)\text{ as $r\todown 0$}.
\]
\item\label{bvblowupjump} For $\Hcal^{d-1}$-almost every $x\in\Jcal_u$,
\[
u^r\to \frac{1}{2}\begin{cases}u^+(x)-u^-(x)\quad\text{ if }\ip{z}{n_{u}(x)}\geq 0,\\
u^-(x)-u^+(x)\quad\text{ if }\ip{z}{n_{u}(x)}<0,
\end{cases}
\]
strictly in $\BV(\bB^d;\R^m)$ as $r\downarrow 0$.
\item For $|D^c u|$-almost every $x\in\Omega$ and for any sequence $r_n\downarrow 0$, the sequence $(u^{r_n})_n$ contains a subsequence \RED{converging} weakly* in $\BV(\bB^d;\R^m)$ to a non-constant limit function \RED{of the form
\begin{equation}\label{eqcantorblowup}
z\in\R^d\mapsto a(x)\gamma\left(\ip{z}{ b(x)}\right)
\end{equation}
where $\gamma\in\BV((-1,1);\R)$ is non-constant and increasing, and $a(x)\in\pd\bB^m$, $ b(x)\in\pd\bB^d$ satisfy
\[
\frac{\dd D^c u}{\dd|D^c u|}(x)=a(x)\otimes b(x).
\]}
Moreover, if $(u^{r_n})_{n}$ is a sequence converging weakly* in this fashion, then, for any $\varepsilon>0$, there exists $\tau\in(1-\varepsilon,1)$ such that the sequence $(u^{\tau r_n})_{n}$ converges strictly in $\BV(\bB^d;\R^m)$ to a limit of the form described by~\eqref{eqcantorblowup}.
\end{enumerate}
In all three situations, we denote $\lim_r u^r$ (or $\lim_n u^{r_{n}}$) by $u^0$. If the base (blow-up) point $x$ needs to be specified explicitly to avoid ambiguity, then we shall write $u^r_{x}$, $u^{r_n}_{x}$ and $u^0_{x}$.
\end{theorem}

We refer to Theorem~2.4 in~\cite{RindlerShaw19} for a proof (which only uses standard results in the theory of BV-functions and Alberti's Rank-One Theorem~\cite{Albe93ROPD}).

For $x\in\Jcal_u$, the function $u^0$ gives a 'vertically recentered' description of the behaviour of $u$ near $x$. It will be convenient to have a compact notation for also describing this behaviour when $u$ is not recentered.

\begin{definition}\label{defjumpblowup}
For $x\in\Jcal_u$, define $u^\pm\in\BV(\bB^d;\R^m)$ by
\[
u^\pm(z):=\begin{cases}
u^+(x) &\text{ if }\ip{z}{n_u(x)}\geq 0,\\
u^-(x) &\text{ if }\ip{z}{n_u(x)}<0.
\end{cases}
\]
If the choice of base point $x\in\Jcal_u$ needs to be emphasised for clarity, we shall write $u^\pm_x$.
\end{definition}
This definition is independent of the choice of orientation $(u^+,u^-,n_u)$ and, for $\Hcal^{d-1}$-almost every $x\in\Jcal_u$, the rescaled function $u(x+r\frarg)$ converges strictly to $u^\pm$ as $r\downarrow 0$.

\subsection{Liftings}
The proof of Theorem~\ref{L1lscthm} makes essential use of the theory of liftings as developed in Section~3 of~\cite{RindlerShaw19}, although we note that some of these ideas were first explored by Jung \& Jerrard in~\cite{JunJer04SCML}. Liftings are graph-like measures associated to $\BV(\Omega;\R^m)$ functions which have good continuity properties and are very useful for the study of $u$-dependent functionals over $\BV(\Omega;\R^m)$, especially for localisation and blow-up arguments. They are used extensively in~\cite{RindlerShaw19} and serve as the primary technical tool for the calculations performed therein, but here we will only need to use a result which concerns their behaviour under rescaling.

\begin{definition}[Elementary Liftings]\label{defelementaryliftings}
Given $u\in\BV_\#(\Omega;\R^m)$, the \textbf{elementary lifting} $\gamma\asc{u}\in\mbfM(\Omega\times\R^m;\R^{m\times d})$ associated to $u$ is defined by
\[
\gamma\asc{u}:=|Du|\otimes\left(\frac{\dd Du}{\dd|Du|}\int_0^1\delta_{u^\theta}\;\dd\theta\right),
\]
that is,
\[
\ip{\varphi}{\gamma\asc{u}}=\int_\Omega\int_0^1\varphi(x,u^\theta(x))\;\dd\theta\;\dd Du(x)\quad\text{for all }\varphi\in\C_0(\Omega\times\R^m),
\]
where $u^\theta$ is the jump interpolant defined in Section~\ref{secpreliminaries}.
\end{definition}
Following the definition given in Section~\ref{subsecmeasuretheory}, we note that, if $u\in\BV_\#(\Omega;\R^m)$, then $\gamma\asc{u}\restrict((\Omega\setminus\Jcal_u)\times\R^m)=\gr^u_\#(\nabla u\lL\restrict\Omega+D^cu)$ is a $u$-graphical measure.
\begin{theorem}[Tangent Liftings at diffuse points]\label{thmdiffusetangentlifting}
Let $x_0\in\Dcal_u\cup\Ccal_u$ and $u^r$, $u^0\in\BV(\bB^d;\R^m)$ be as defined in Theorem~\ref{thmbvblowup}. Then, for $\lL+|D^c u|$-almost every $x_0\in\Dcal_u\cup\Ccal_u$,
\begin{itemize}
\item if $x_0\in\Dcal_u$ then $u^r\to u^0$ and $\gamma\asc{u^r}\to\gamma\asc{u^0}$ strictly in $\mbfM(\Omega\times\R^m;\R^{m\times d})$ as $r\to 0$,
\item if $x_0\in\Ccal_u$ then for any sequence $r_n\downarrow 0$ and $\varepsilon>0$, there exists $\tau\in(1-\varepsilon,1)$ such that $u^{\tau r_n}\to u^0$ and $\gamma\asc{u^{\tau r_n}}\to\gamma\asc{u^0}$ strictly in $\mbfM(\Omega\times\R^m;\R^{m\times d})$ as $n\to\infty$.
\end{itemize}
\end{theorem}
\RED{A key property of liftings is that the elementary lifting $\gamma\asc{u^r}$ associated to a rescaled $\BV$ function can be written as a rescaling of the elementary lifting $\gamma\asc{u}$ associated to the original $\BV$ function:
\begin{equation}\label{eqrescaledliftings}
\gamma\asc{u^{r}}=\frac{c_r}{r^d}(T^{(x,r), ((u)_{x,r}, c_r^{-1}r)})_\# \bigl( \gamma\asc{u}\restrict (B(x,r)\times\R^m) \bigr).
\end{equation}
The identity~\eqref{eqrescaledliftings} }and Theorem~\ref{thmdiffusetangentlifting} are proved in~\cite{RindlerShaw19} as Lemma~3.16 and Theorem~3.17, respectively.

\subsection{Integrands}
\begin{definition}[Recession functions]\label{defrecessionfunction}
For $f\colon\overline{\Omega}\times\R^m\times\R^{m\times d}\to\R$, define the \textbf{(strong) recession function} $f^\infty\colon\overline{\Omega}\times\R^m\times\R^{m\times d}\to\R$ of $f$ by
\begin{equation*}
 f^\infty\left(x,y,A\right)=\lim_{\substack{(x_k,y_k,A_k)\to (x,y,A)\\t_k\to\infty}}\frac{f\left(x_k,y_k,t_kA_k\right)}{t_k},
\end{equation*}
whenever the right-hand side exists \textcolor{black}{in $\R$} for every $(x,y,A)\in\overline{\Omega}\times\R^m\times\R^{m\times d}$ independently of the order in which the limits of the individual sequences $(x_k,y_k,A_k)_k\subset\overline{\Omega}\times\R^m\times\R^{m\times d}$, $(t_k)_k\subset(0,\infty)$ are taken and of the sequences used. The definition of $f^\infty$ implies that, whenever it exists, $f^\infty$ must be continuous.
\end{definition}

\begin{definition}[Representation integrands]\label{defrepresentationf}
A function $f\colon\overline{\Omega}\times\R^m\times\R^{m\times d}\to\R$ is said to be a member of $\mbfR(\Omega\times\R^m)$ if
$f$ is Carath\'eodory and its recession function $f^\infty$ exists.
\end{definition}

\begin{definition}\label{defrepresentationfl1}
An integrand $f\in\mbfR(\Omega\times\R^m)$ is said to be a member of the set $\RBVL(\Omega\times\R^m)$ if there exists $g\in\C(\overline{\Omega}\times\R^m;[0,\infty))$ satisfying
	\begin{enumerate}[(a)]
	\item\label{prelimgexists} there exists $C>0$ such that
	\[	
	0\leq g(x,y)|A|\leq f(x,y,A)\leq Cg(x,y)(1+|A|);
	\]
	\item\label{prelimgbound} there exist $R>0$ and $C>1$ for which $|y|\geq R$ and $t\geq 1$ imply $g(x,y)\leq C g(x,ty)$;
	\item\label{prelimgcompatible} for every $K\Subset\R^m$ and every $\varepsilon>0$, there exists $R_\varepsilon>0$ such that $|A|\geq R_\varepsilon$ implies $|(f-f^\infty)(x,y,A)|\leq\varepsilon g(x,y)(1+|A|)$ for all $(x,y)\in\overline{\Omega}\times K$.
	\end{enumerate}
\end{definition}

\subsection{Functionals and surface energies}\label{subsecfunctionals}
For $f\in\Rbf(\Omega\times\R^m)$, we define the extended functional $\Fcal\colon\BV(\Omega;\R^m)\to\R \cup\{+\infty\}$ by
\begin{equation}\label{eqctsfunctional}
\Fcal[u]:=\int_\Omega f(x,u(x),\nabla u(x))\;\dd x+\int_\Omega\int_0^1 f^\infty\left(x,u^\theta(x),\frac{\dd D^su}{\dd|D^su|}(x)\right)\;\dd\theta\;\dd|D^su|(x),
\end{equation}
where $u^\theta$ is the jump interpolant defined \RED{previously} by~\eqref{eqdefjumpinterpolant}. Theorem~\ref{thmareastrictcontinuity} below states that $\Fcal$ as defined by~\eqref{eqctsfunctional} is the area-strictly continuous extension of $u\mapsto\int_\Omega f(x,u(x),\nabla u(x))\;\dd x$ from $\W^{1,1}(\Omega;\R^m)$ (or \RED{$(\C^\infty\cap\W^{1,1})(\Omega;\R^m)$}) to $\BV(\Omega;\R^m)$. Proposition~\ref{propfixedbdary} therefore implies that Theorem~\ref{L1lscthm} can equivalently be seen as identifying the weak* relaxation of this continuously extended $\Fcal$ from $\BV(\Omega;\R^m)$ to $\BV(\Omega;\R^m)$, which is the approach that we take in what follows.

\begin{theorem}\label{thmareastrictcontinuity}
 Let $\Omega\subset\mbR^d$ be a bounded domain with Lipschitz boundary and let $f\in \mbfR(\Omega\times\R^m)$ satisfy the growth bound
\begin{equation}\label{eq:growthbound}
|f(x,y,A)|\leq C(1+|y|^{d/(d-1)}+|A|) \quad\text{ for all } (x,y,A)\in\Omega\times\mbR^m\times\mbR^{m\times d}.
\end{equation}
Then the functional $\Fcal\colon\BV(\Omega;\R^m)\to\R$ is area-strictly continuous.
\end{theorem}
Theorem~\ref{thmareastrictcontinuity} is proved under slightly more general hypotheses in~\cite{RinSha15SCE} as Theorem~5.2.

Given $u\in\BV(\Omega;\R^m)$ and $x\in\Jcal_u$, define the class of functions $\Acal_u(x)$ by
\begin{equation*}
\Acal_u(x):=\left\{\varphi\in\left(\C^\infty\cap\W^{1,1}\cap\Lp^\infty\right)(\bB^d;\R^m)\colon\;\varphi= u^\pm_{x}\text{ on }\pd\bB^d\right\},
\end{equation*}
where $u^\pm_x$ is as given in Definition~\ref{defjumpblowup} \RED{and equality on $\pd\bB^d$ is to be understood in the sense of Gagliardo's Trace Theorem,~\cite{Gagl57CTFR}}. For $f\in\Rbf(\Omega\times\R^m)$ and $u\in\BV(\Omega;\R^m)$, the surface energy densities $K_f[u]$ and $H_f[u]$ at $x\in\Jcal_u$ are defined by
\begin{equation*}
K_f[u](x):=\inf\left\{\frac{1}{\omega_{d-1}}\int_{\bB^d}f^\infty(x,\varphi(z),\nabla\varphi(z))\;\dd z\colon\;\varphi\in\Acal_u(x)\right\},
\end{equation*}
and
\begin{equation} \label{eq:Hdef}
H_f[u](x):=\liminf_{r\to 0}H^r_f[u](x),
\end{equation}
respectively, where $H^r_f[u]$ is given for each $r>0$ by
\begin{equation*}
H_f^r[u](x):=\inf\left\{\frac{1}{\omega_{d-1}}\int_{\bB^d}f^\infty\left(x+rz,\varphi(z),\nabla\varphi(z)\right)\;\dd z\colon\;\varphi\in\Acal_u(x),\;\norm{\varphi}_{\Lp^{1}}\leq 2\norm{u_x^\pm}_{\Lp^{1}}\right\}.
\end{equation*}

Example~\ref{exdensitiesnotthesame} in Section~\ref{secrelaxationcomparisons} demonstrates that, for general integrands $f\in\Rbf(\Omega\times\R^m)$, it may occur that $K_f\neq H_f$. As a consequence, the weak* and $\Lp^1$ relaxations, $\Fcalrw$ and $\Fcalro$, of $\Fcal$ should be expected to differ in general. This phenomenon is not apparent in the earlier works~\cite{FonMul93RQFB,BoFoMa98GMR,FonLeo01LSC}, since the assumptions on $f^\infty$ which feature therein (usually that $f^\infty(x_0,y,A)\leq (1+\varepsilon)f^\infty(x,y,A)$ uniformly in $(y,A)$ for each $x_0$ with $|x-x_0|$ sufficiently small) are strong enough to enforce $K_f=H_f$.

\RED{\begin{remark}\label{rem:K_f}
Note that in~\cite{FonMul93RQFB} the set $\Acal_u$ is defined using a different class of test functions, which are both defined on cubes rather than balls and which are not required to be fixed on the entirety of the boundary. The restriction of test functions to those whose behaviour on the boundary is completely fixed as in the case of our $\Acal_u$ was later obtained in~\cite{FonLeo07MMCV,BoFoMa98EBER,BoFoMa98GMR}. It can also be seen that balls and cubes as domains are interchangeable from the perspective of defining $K_f$ by using a covering argument to tile the bisecting hyperplane (with the appropriate orientation) of a cube with balls and vice versa. 
\end{remark}}

It follows from the definition $H_f[u](x):=\liminf_{r\to 0}H^r_f[u](x)$ that for each $x\in\Jcal_u$ there exists a sequence $r_n\downarrow 0$ such that $H_f[u](x)=\lim_{n\to\infty}H_f^{r_n}[u](x)$. A priori, the sequence $(r_n)_n$ depends on $x$ and we do not know whether a sequence $r_n\downarrow 0$ exists with respect to which $\lim_{n\to\infty}H_f^{r_n}[u](x)=H_f[u](x)$ for all $x\in\Jcal_u$ simultaneously. Fortunately, Lemma~\ref{lemenergyapprox} below shows that it is the case that $\liminf_{n\to\infty} H_f^{1/n}[u](x)=H_f[u](x)$ for every $x\in\Jcal_u$, so that we can always restrict our attention to the fixed sequence $(\frac{1}{n})_{n\in\mbN}$ for the purposes of computing $H_f[u]$.

\begin{lemma}\label{lemenergyapprox}
If $f\in\Rbf(\Omega\times\R^m)$ then, for $\Hcal^{d-1}$-almost every $x\in\Jcal_u$, it holds that
\[
H_f[u](x)=\liminf_{n\to\infty}H^\frac{1}{n}_f[u](x).
\]
\end{lemma}

\begin{proof}
The inequality
\[
H_f[u]\leq\liminf_{n\to\infty}H^\frac{1}{n}_f[u]
\]
follows immediately from the definition of $H_f[u]$, and so it remains to show that $H_f[u]\geq\liminf_nH^\frac{1}{n}_f[u]$. To see this, let $r_n\downarrow 0$ and $(\varphi_n)_n\subset\Acal_u(x)$ be sequences such that $\norm{\varphi_n}_{\Lp^{1}}\leq 2\norm{u_x^\pm}_{\Lp^{1}}$ and
\[
\omega_{d-1}H_f[u](x)=\lim_{n\to\infty}\int_{\bB^d}f^\infty(x+r_nz,\varphi_n(z),\nabla\varphi_n(z))\;\dd z.
\]
Let $(k_n)_n\subset\mbN$ be an increasing sequence such that $\frac{1}{k_n+1}\leq r_n<\frac{1}{k_n}$ and define \RED{$w_n\in(\C^\infty\cap\W^{1,1})({k_nr_n}\bB^d;\R^m)$} by $w_n(z):=\varphi_n((k_nr_n)^{-1}z)$. After a change of variables, we then have that
\[
\omega_{d-1}H_f[u](x)=\lim_{n\to\infty}(k_nr_n)^{1-d}\int_{{k_nr_n}\bB^d}f^\infty\left(x+\frac{z}{k_n},w_n(z),\nabla w_n(z)\right)\;\dd z.
\]
Since $k_nr_n<1$ and $w_n$ satisfies $w_n(z)=u^\pm((k_nr_n)^{-1}z)$ for $z\in\pd(k_nr_n\bB^d)$, we can extend each $w_n$ to an element of $\BV(\bB^d;\R^m)$ by setting $w_n := u_x^\pm$ in $\bB^d\setminus k_nr_n \bB^d$. \textcolor{black}{For each $n\in\mbN$, we can use Proposition~\ref{propfixedbdary} to obtain a sequence $(v_n^i)_{i\in\mbN}\subset(\W^{1,1}\cap\C^\infty)(\bB^d;\R^m)$ with $\sup_i\norm{v_n^i}_\infty<\infty$, ${v_n^i}_{|\pd\bB^d}=u^\pm_x$, and $v_n^i\to w_n$ area strictly in $\BV(\bB^d;\R^m)$ as $j\to\infty$.  Applying Theorem~\ref{thmareastrictcontinuity} for each $n\in\mbN$ to an integrand of the form $(x,y,A)\mapsto g(|y|)f^\infty(x,y,A)$, where $g\in\C_c([0,\infty);[0,\infty))$ is such that $g(t)=1$ whenever $|t|\leq \sup_i\norm{v_n^i}_\infty + \norm{w_n}_\infty$, we can therefore always find $i_n\in\mbN$ such that
\[
\left|\int_{\bB^d} f^\infty\left(x + \frac{z}{k_n}, w_n(z), \nabla w_n(z)\right)\;\dd z - \int_{\bB^d}f^\infty\left(x + \frac{z}{k_n}, v^{i_n}_n(z), \nabla v^{i_n}_n(z)\right)\;\dd z\right| < \frac{1}{n}.
\]
Letting $\widetilde{w}_n:= v_n^{i_n}$ and recalling the extension of $w_n$ to $\bB^d$, we therefore have a sequence  $(\widetilde{w}_n)_n\subset\Acal_u(x)$ which is such that 
\begin{align*}
\left|\int_{{k_nr_n}\bB^d}f^\infty\left(x+\frac{z}{k_n},w_n(z),\nabla w_n(z)\right)\;\dd z-\int_{{k_nr_n}\bB^d}f^\infty\left(x+\frac{z}{k_n},\widetilde{w}_n(z),\nabla \widetilde{w}_n(z)\right)\;\dd z\right|\to 0
\end{align*}
as $n\to\infty$ and
\begin{align*}
&\lim_{n\to\infty}\int_{\bB^d\setminus{k_nr_n}\bB^d}f^\infty\left(x+\frac{z}{k_n},\widetilde{w}_n(z),\nabla \widetilde{w}_n(z)\right)\;\dd z\\
=&\lim_{n\to\infty}\int_{\bB^d\setminus{k_nr_n}\bB^d}f^\infty\left(x+\frac{z}{k_n},u^\pm_x(z),\frac{\dd Du^\pm_x}{\dd|D u^\pm_x|}(z)\right)\;\dd|Du^\pm_x|(z)=0.
\end{align*}
Since each $\varphi_n$ satisfies $\norm{\varphi_n}_{\Lp^1(\bB^d;\R^m)}\leq 2\norm{u^\pm_x}_{\Lp^1(\bB^d;\R^m)}$,  we have that $\norm{w_n}_{\Lp^1(\bB^d;\R^m)}<2\norm{u^\pm_x}_{\Lp^1(\bB^d;\R^m)}$ (note $\norm{w_n}_{\Lp^1(k_nr_n\bB^d;\R^m)}\leq 2(k_nr_n)^d\norm{u^\pm_x}_{\Lp^1(\bB^d;\R^m)}$  and $\norm{u^\pm_x}_{\Lp^1(\bB^d\setminus k_nr_n\bB^d;\R^m)}=(1-(k_nr_n)^d)\norm{u^\pm_x}_{\Lp^1(\bB^d;\R^m)}$), and can therefore guarantee that $\norm{\widetilde{w}_n}_{\Lp^1(\bB^d;\R^m)}\leq 2\norm{u^\pm_x}_{\Lp^1(\bB^d;\R^m)}$ for each $n$.}

Because $k_nr_n\to 1$, we therefore have that
\begin{align*}
H_f[u](x)&=\lim_{n\to\infty}\frac{1}{\omega_{d-1}}\int_{\bB^d}f^\infty\left(x+\frac{z}{k_n},\widetilde{w}_n(z),\nabla \widetilde{w}_n(z)\right)\;\dd z\\
&\geq\liminf_{k_n\to\infty}H_f^{\frac{1}{k_n}}[u](x)\\
&\geq\liminf_{n\to\infty}H_f^{\frac{1}{n}}[u](x),
\end{align*}
as required.
\end{proof}

It is not clear a priori that the densities $K_f[u]$ and $H_f[u]$ are $\Hcal^{d-1}$-measurable and hence whether the integrals
\[
\int_{\Jcal_u}K_f[u](x)\;\dd\Hcal^{d-1}(x),\qquad \int_{\Jcal_u}H_f[u](x)\;\dd\Hcal^{d-1}(x),
\]
are always well-defined for every $u\in\BV(\Omega;\R^m)$. Lemma~\ref{lemenergymeasurable} below shows that the density $H^r_f[u]$ is always $\Hcal^{d-1}$-measurable; the corresponding result for $K_f[u]$ was shown in Lemma~2.15 and Corollary~2.16 of~\cite{RindlerShaw19}. By Lemma~\ref{lemenergyapprox}, $H_f[u]$ is then the pointwise infimum of a sequence of $\Hcal^{d-1}$-measurable functions and is therefore itself measurable.

\begin{lemma}\label{lemenergymeasurable}
If $f\in\Rbf(\Omega\times\R^m)$ and $r>0$ is sufficiently small, then $H_f^r[u]$ is $\Hcal^{d-1}\restrict\Jcal_u$-measurable and equal $\Hcal^{d-1}\restrict\Jcal_u$-almost everywhere to an upper-semicontinuous function.
\end{lemma}

\begin{proof}
First, fix a triple $(u^+,u^-,n_u)\colon\Jcal_u\to\R^m\times\R^m\times\pd\bB^d$ such that $n_u$ orients $\Jcal_u$ and $u^+$, $u^-$ are the one-sided jump limits of $u$ with respect to $n_u$. Fix also $\varepsilon>0$. The triple $(u^+,u^-,n_u)$ is Borel and hence $|D^ju|$-measurable, and so Lusin's Theorem (see Theorem~1.45 in~\cite{AmFuPa00FBVF}) implies that there exists a compact set $K_\varepsilon\Subset\Jcal_u$ such that $|D^ju|(\Jcal_u\setminus K_\varepsilon)\leq\varepsilon$ and $(u^+,u^-,n_u)$ is continuous when restricted to $K_\varepsilon$.

Let $x\in K_\varepsilon$ and $(x_j)_j\subset K_\varepsilon$ be such that $x_j\to x$. For each $j\in\mbN$ let $R_j\colon\bB^d\to\bB^d$ be a rotation mapping $n_u(x)$ to $n_u(x_j)$ such that $R_j\to\id_{\R^d}$ as $j\to\infty$ and let $S_j\colon\R^m\to\R^m$ be a sequence of linear maps mapping $u^+(x)$ to $u^+(x_j)$ and $u^-(x)$ to $u^-(x_j)$ such that $S_j\to\id_{\R^m}$ as $j\to\infty$ (such a choice of $(R_j)_j$, $(S_j)_j$ is possible by the fact that $x,x_j\in K_\varepsilon$). Now for $\delta>0$, let $\varphi\in\Acal_u(x)$ with $\norm{\varphi}_{\Lp^{1}}\leq 2\norm{u_x^\pm}_{\Lp^{1}}$ be such that
\[
\frac{1}{\omega_{d-1}}\int_{\bB^d}f^\infty(x+rz,\varphi(z),\nabla \varphi(z))\;\dd z\leq H_f^r[u](x)+\delta.
\]
Define $\varphi_j\in\C^\infty(\bB^d;\R^m)$ by $\varphi_j(z):=S_j(\varphi(R_jz))$ and note that $\varphi_j\in\Acal_u(x_j)$. By the convergence properties assumed of $(R_j)_j$ and $(S_j)_j$, we have that $\varphi_j\to\varphi$ strictly in $\BV(\bB^d;\R^m)$ as $j\to\infty$. For $r> 0$, define $\mu_j^r\in\mbfM(\Omega\times\R^m;\R^{m\times d})$ by
\[	
\mu_j^r:=\left(T^{(x_j,r),(-(\varphi_j)_{\bB^d},1)}\right)^{-1}_\#\gamma\asc{\varphi_j-(\varphi_j)_{\bB^d}},
\]
where $\gamma\asc{\varphi_j-(\varphi_j)_{\bB^d}}$ is \RED{the elementary lifting associated to $\varphi_j-(\varphi_j)_{\bB^d}$} as given by Definition~\ref{defelementaryliftings}, and let $\mu^r\in\mbfM(\Omega\times\R^m;\R^{m\times d})$ be given by
\[
\mu^r:=\left(T^{(x,r),(-(\varphi)_{\bB^d},1)}\right)^{-1}_\#\gamma\asc{\varphi-(\varphi)_{\bB^d}},
\]
It can easily be seen that $\mu^r_j$ converges strictly in $\mbfM(\bB^d\times\R^m;\R^{m\times d})$ to $\mu^r$ as $j\to\infty$. Using Reshetnyak's Continuity Theorem (see Theorem~2.39 in~\cite{AmFuPa00FBVF}) and the positive one-homogeneity of $f^\infty$, we therefore deduce
\begin{align*}
\int_{\bB^d}f^\infty(x_j+rz,\varphi_j(z),\nabla\varphi_j(z))\;\dd z&=\int_{\bB^d\times\R^m}f^\infty\left(z,y,\frac{\dd\mu^r_j}{\dd|\mu^r_j|}(z,y)\right)\;\dd|\mu^r_j|(z,y)\\
&\to\int_{\bB^d\times\R^m}f^\infty\left(z,y,\frac{\dd\mu^r}{\dd|\mu^r|}(z,y)\right)\;\dd|\mu^r|(z,y)\\
&=\int_{\bB^d}f^\infty(x+rz,\varphi(z),\nabla\varphi(z))\;\dd z
\end{align*}
as $j\to\infty$ for any $r>0$. By our choice of $\varphi$ and the boundary condition satisfied by each $\varphi_j$, we therefore have that
\begin{align*}
H^r_f[u](x)+2\delta&\geq\lim_{j\to\infty}\frac{1}{\omega_{d-1}}\int_{\bB^d}f^\infty(x_j+rz,\varphi_j(z),\nabla\varphi_j(z))\;\dd z\geq\limsup_{j\to\infty}H^r_f[u](x_j).
\end{align*}
It follows from the arbitrariness of $x\in K_\varepsilon$ and $\delta>0$ that $H^r_f[u]$ is upper semicontinuous when restricted to $K_\varepsilon$. Finally, define $F^r_\varepsilon\colon\Jcal_u\to[0,\infty]$ by
\[
F^r_\varepsilon(x):=\begin{cases}
H^r_f[u](x)&\text{ if }x\in K_\varepsilon,\\
\infty&\text{ otherwise},
\end{cases}
\]
and note that $F^r:=\inf_{\varepsilon>0}F^r_\varepsilon$ is equal to $H^r_f[u]$ at $|D^ju|$-almost every $x\in\Jcal_u$ and hence $\Hcal^{d-1}\restrict\Jcal_u$-almost every $x\in\Jcal_u$. The conclusion now follows from the fact that the pointwise infimum of a collection of upper semicontinuous functions is upper semicontinuous and hence measurable.
\end{proof}

\begin{corollary}\label{correctifiabledensities}
If $f\in\Rbf(\Omega\times\R^m)$ and $u\in\BV(\Omega;\R^m)$, then $H_f[u]$ is $\Hcal^{d-1}\restrict\Jcal_u$-measurable. If in addition $f^\infty\geq 0$ and
\[
\int_{\Jcal_u}H_f[u](x)\;\dd\Hcal^{d-1}(x)<\infty,
\]
then $H_f[u]\Hcal^{d-1}\restrict\Jcal_u\in\mbfM^+(\Omega)$ and is a $(d-1)$-rectifiable measure.
\end{corollary}

\begin{proof}
Lemma~\ref{lemenergyapprox} combined with Lemma~\ref{lemenergymeasurable} implies that $H_f[u]$ is the pointwise infimum of a countable collection of $\Hcal^{d-1}\restrict\Jcal_u$-measurable functions and is therefore also $\Hcal^{d-1}\restrict\Jcal_u$-measurable. Since $\Hcal^{d-1}\restrict\Jcal_u$-measurable functions coincide $\Hcal^{d-1}$-almost everywhere with Borel functions (see Exercise~1.3 in~\cite{AmFuPa00FBVF}) the second part of the Corollary follows from the discussion on rectifiability in Section~\ref{subsecmeasuretheory}.
\end{proof}

The following theorem is the main result of~\cite{RindlerShaw19} and will be used in Section~\ref{chapl1lsc}.
\begin{theorem}\label{thmwsclsc}
Let $f\colon\overline{\Omega}\times\R^m\times\R^{m\times d}\to\R$ where $d\geq 2$ and $m\geq 1$ be such that
\begin{enumerate}[(i)]
\item\label{wsthmhypoth1} $f$ is a Carath\'eodory function whose recession function $f^\infty$ exists in the sense of Definition~\ref{defrecessionfunction} and satisfies $f^\infty\geq 0$;
\item\label{wsthmhypoth2} $f$ satisfies a growth bound of the form 
\begin{equation*}
-C(1+|y|^p+h(A))\leq f(x,y,A)\leq C(1+|y|^{d/(d-1)}+|A|)
\end{equation*}
for some $C>0$, $p\in[1,d/(d-1))$, $h\in\C(\R^{m\times d})$ satisfying $h^\infty\equiv 0$, and for all $(x,y,A)\in\overline{\Omega}\times\R^m\times\R^{m\times d}$;
\item $f(x,y,\frarg)$ is quasiconvex for every $(x,y)\in\overline{\Omega}\times\R^m$.
\end{enumerate}
Then the sequential weak* relaxation of $\Fcal$ to $\BV(\Omega;\R^m)$ is given by
\begin{align*}
\Fcalrw[u]&=\int_\Omega f(x,u(x),\nabla u(x))\;\dd x+\int_\Omega f^\infty\left(x,u(x),\frac{\dd D^c u}{\dd|D^c u|}(x)\right)\;\dd|D^cu|(x)\\
&\qquad+\int_{\Jcal_u} K_f[u](x)\;\dd\Hcal^{d-1}(x).
\end{align*}
\end{theorem}

\section{$\Lrm^1$-relaxation versus $(\BV,\mathrm{weak*})$-relaxation} \label{secrelaxationcomparisons}
In this section, we clarify the relationship between our $\Lp^1$ relaxation result (Theorem~\ref{L1lscthm}), the previous $\Lp^1$ relaxation results found in~\cite{FonMul93RQFB,FonLeo01LSC}, and Theorem~\ref{thmwsclsc}. 

Our work shows that, when the hypotheses of Theorem~\ref{L1lscthm} are met, the strong $\Lp^1(\Omega;\R^m)$-relaxation of the functional $\Fcal$ defined by~\eqref{eqfunctional} to $\BV(\Omega;\R^m)$ is given by
\begin{align}
\begin{split}\label{eqourL1relaxation}
\Fcalro[u]&=\int_\Omega f(x,u(x),\nabla u(x))\;\dd x+\int_\Omega f^\infty\left(x,u(x),\frac{\dd D^c u}{\dd|D^c u|}(x)\right)\;\dd|D^cu|(x)\\
&\qquad+\int_{\Jcal_u} H_f[u](x)\;\dd\Hcal^{d-1}(x).
\end{split}
\end{align}
By contrast, the authors of~\cite{FonMul93RQFB,FonLeo01LSC}, assuming that $f$ is partially coercive as well as satisfying some additional technical hypotheses (see below), obtain the formula
\begin{align}
\begin{split}\label{eqFonMulL1relaxation}
\Fcalro[u]&=\int_\Omega f(x,u(x),\nabla u(x))\;\dd x+\int_\Omega f^\infty\left(x,u(x),\frac{\dd D^c u}{\dd|D^c u|}(x)\right)\;\dd|D^cu|(x)\\
&\qquad+\int_{\Jcal_u} K_f[u](x)\;\dd\Hcal^{d-1}(x).
\end{split}
\end{align}
In fact it can be shown that, under reasonable hypotheses on $f$,
\[
K_f[u](x)=\inf\left\{\frac{1}{\omega_{d-1}}\int_{\bB^d}f^\infty\left(x,\varphi(z),\nabla\varphi(z)\right)\;\dd z\colon\;\varphi\in\Acal_u(x)\text{ and }\norm{\varphi}_{\Lp^{1}}\leq 2\norm{u_x^\pm}_{\Lp^{1}}\right\},
\]
which makes it clear that $K_f[u](x)=H^0_f[u](x)$.

Since
\[
  \lim_{r\to 0} \int_{\bB^d}f^\infty\left(x+rz,\varphi(z),\nabla\varphi(z)\right)\;\dd z=\int_{\bB^d}f^\infty\left(x,\varphi(z),\nabla\varphi(z)\right)\;\dd z
\]
for any fixed $\varphi\in\Acal_u(x)$, we have $H_f[u]\leq K_f[u]$ whenever $f^\infty$ is continuous. Example~\ref{exdensitiesnotthesame} below shows that strict inequality $H_f[u]<K_f[u]$ can sometimes occur and hence that~\eqref{eqourL1relaxation} and~\eqref{eqFonMulL1relaxation} cannot simultaneously be true for arbitrary integrands $f$. That our results are compatible with those of~\cite{FonMul93RQFB,FonLeo01LSC} is due to the lower semicontinuity properties of $f$ which are assumed therein: in particular, these include the statement that for every $x\in\Omega$ and $\varepsilon>0$ there exists $r>0$ such that
\begin{equation}\label{eqfonassumption}
f^\infty(x,y,A) \leq (1+\varepsilon)f^\infty(x+rz,y,A) + \varepsilon\quad\text{ for all }(z,y,A)\in \bB^d\times\R^m\times\R^{m\times d}.
\end{equation}
This clearly implies that $K_f[u](x)\leq (1+\varepsilon)H^r_f[u](x)+\varepsilon$ for all $r$ sufficiently small and hence that $K_f[u]\equiv H_f[u]$ whenever the hypotheses of~\cite{FonMul93RQFB} apply.

To date, the results obtained in~\cite{FonMul93RQFB,FonLeo01LSC} are the only ones in this area which are valid for the full $u$-dependent, vector-valued $m>1$ case. The authors obtain~\eqref{eqFonMulL1relaxation} (or at least the lower bound ``$\geq$'') by assuming either 
\begin{itemize}
\item strong lower semicontinuity properties of $f$ in $(x,y)$ (but without requiring any kind of coercivity from $f$), or 
\item a combination of partial coercivity together with a weaker set of lower semicontinuity properties (including~\eqref{eqfonassumption}).
\end{itemize}
Our aim for this work was to obtain a sensible formula for $\Fcalro$ by requiring just partial coercivity from $f$ and without making any lower semicontinuity assumptions in the $(x,y)$-variables. The point is that, whilst integrands $f$ which arise from problems in the theory of phase transitions must always be partially coercive, it is not clear that they should satisfy the $(x,y)$ lower semicontinuity assumptions which have been made so far. The integrand $f(x,y,A)=|y|^{1-|x|}|A|$ defined on $\bB^d\times\R^m\times\R^{m\times d}$ arises naturally from the family of perturbed functionals given by~\eqref{eqperturbfamily} when $g(x,y)=|y|^{1-|x|}$ and $h(A)=|A|$, for instance, but does not satisfy~\eqref{eqfonassumption}. Theorem~\ref{L1lscthm} shows that a reasonable formula for $\Fcalro$ can indeed be obtained under partial coercivity without needing to make specific lower semicontinuity requirements from $f$. 

It is interesting to note that the formula~\eqref{eqFonMulL1relaxation} obtained in~\cite{FonMul93RQFB,FonLeo01LSC} is the same as the one appearing in Theorem~\ref{thmwsclsc}, which implies that $\Fcalrw=\Fcalro$ whenever the hypotheses in~\cite{FonMul93RQFB,FonLeo01LSC} are met. Whilst our Theorem~\ref{L1lscthm} shows that an integral representation for $\Fcalro$ can be obtained without requiring lower semicontinuity properties from $f$, it is an open question as to whether the equality $H_f=K_f$ can be obtained without assuming these and, if so, what the most natural alternative hypotheses guaranteeing $H_f=K_f$ are.

The following example shows that the relaxation formulae appearing in~\cite{FonMul93RQFB,RindlerShaw19} and Theorem~\ref{L1lscthm} are not the same in general. 
\RED{The integrand $f$ given below does not satisfy the hypotheses of Theorem~\ref{L1lscthm} (in particular, hypothesis~\eqref{eql1intro1}), although it does satisfy the hypotheses of Theorem~\ref{thmwsclsc}. It is unclear if an integrand exists which jointly satisfies the assumptions of Theorems~\ref{L1lscthm} and~\ref{thmwsclsc} and which also exhibits $K_f\neq H_f$.}
\begin{example}\label{exdensitiesnotthesame}
Let $u\in\BV((-1,1);\R^2)$ be given by 
\[u(x):=
\begin{cases}
\begin{pmatrix}
0\\
0\end{pmatrix} &\text{ if }x\leq 0,\\
\begin{pmatrix}
0\\
1\end{pmatrix}&\text{ if }x>0,
\end{cases}
\]
and define $f\in\Rbf((-1,1)\times\R^2)$ by
\[
f(x,y,A):=\Phi(x,y)|A|,\qquad\Phi\left(x,\begin{pmatrix}y_1\\y_2\end{pmatrix}\right):=\frac{|y_2(1-y_2)|}{1+|y_2|^2}\exp\left(-\sqrt{|x|}\cdot|y_1|\right).
\]
For $j\geq 4$, define the sequence $u_j\in\W^{1,1}((-1,1);\R^2)$ by
\begin{align*}
u_j(x):=&\begin{pmatrix}
j^2x\\
0
\end{pmatrix}\mathbbm{1}_{(0,1/j)}(x)+\begin{pmatrix}
j\\
j(x-1/j)
\end{pmatrix}\mathbbm{1}_{[1/j,2/j]}(x)+\begin{pmatrix}
j-j^2(x-2/j)\\
1
\end{pmatrix}\mathbbm{1}_{(2/j,3/j)}(x)\\
&\qquad+\begin{pmatrix}
0\\
1
\end{pmatrix}\mathbbm{1}_{[3/j,1)}(x).
\end{align*}
Note that $|Du_j|((-1,1))=2j+1$, $u_j\to u$ in $\Lp^1((-1,1);\R^2)$ as $j\to\infty$, and, for any $r>0$,
\begin{align*}
\int_{-1}^1\Phi(rz,u_j(z))|\nabla u_j(z)|\;\dd z & =\int_{1/j}^{2/j}\frac{(jz-1)(2-jz)}{1+|jz-1|^2}\exp\left(-r^{1/2}j\sqrt{|z|}\right)j\;\dd z\\
&\leq\exp\left(-(rj)^{1/2}\right)\int_1^2\frac{(t-1)(2-t)}{1+|t-1|^2}\;\dd t\\
&\to 0\text{ as }j\to\infty,
\end{align*}
from which it follows that $H_f^r[u](0)=0$ for every $r>0$ and hence that $H_f[u](0)=0$.

Next, we identify $K_f[u](0)$: since $f(0,y,A)=\frac{|y_2(1-y_2)|}{1+|y_2|^2}|A|$, we have for any \RED{$v=(v_1,v_2)\in(\C^\infty\cap\W^{1,1})((-1,1);\R^2)$},
\begin{align*}
\int_{-1}^1f\left(0,v(z),\nabla v(z)\right)\;\dd z&=\int_{-1}^1\frac{|v_2(z)(1-v_2(z))|}{1+|v_2(z)|^2}|\nabla v(z)|\;\dd z\\
&\geq \int_{-1}^1\frac{|v_2(z)(1-v_2(z))|}{1+|v_2(z)|^2}|\nabla v_2(z)|\;\dd z.
\end{align*}
If $v$ satisfies $v(-1)=\begin{pmatrix}0\\0\end{pmatrix}$ and $v(1)=\begin{pmatrix}0\\1\end{pmatrix}$, then we can use the substitution $\theta=v_2(z)$ to estimate
\begin{align*}
\int_{-1}^1f\left(0,v(z),\nabla v(z)\right)\;\dd z& \geq \int_{-1}^1\frac{|v_2(z)(1-v_2(z))|}{1+|v_2(z)|^2}|\nabla v_2(z)|\;\dd z\\
& \geq \int_{-1}^1\frac{|v_2(z)(1-v_2(z))|}{1+|v_2(z)|^2}\nabla v_2(z)\;\dd z \\
&= \int_0^1\frac{\theta(1-\theta)}{1+\theta^2}\;\dd\theta.
\end{align*}
Thus,
\[
K_f[u](0)\geq\frac{1}{2}\int_0^1\frac{\theta(1-\theta)}{1+\theta^2}\;\dd\theta>0,
\]
which shows that $K_f[u](0)>H_f[u](0)$.
\end{example}

Nevertheless, Proposition~\ref{proprecessionrecoveryseq} will show that essentially the same procedure can be used to simultaneously construct $\Lp^1$-recovery sequences for $\Fcalrw$ and $\Fcalro$ in the simple case where $f$ is positively $1$-homogeneous and satisfies a bound of the form $0\leq f(x,y,A)\leq C|A|$ for some $C>0$. 
	
\section{Localisation under partial coercivity}\label{chapl1lsc}
This section and the next are concerned with obtaining the lower semicontinuity component of Theorem~\ref{L1lscthm}, that is, the derivation of a lower bound for the $\Lp^1$-relaxation $\Fcalro$ of the functional $\Fcal$. We argue by applying the truncated blow-up method due to Fonseca \& M\"uller~\cite{FonMul92QCIL,FonMul93RQFB}, separately at points $x_0\in\Dcal_u$, $x_0\in\Ccal_u$ and $x_0\in\Jcal_u$. For integrands $f$ which satisfy
\[
g(x,y)|A|\leq f(x,y,A)\leq Cg(x,y)(1+|A|)\quad\text{ for all }(x,y,A)\in\Omega\times\R^m\times\R^{m\times d}
\]
for some $g\in\C(\Omega\times\R^m;[0,\infty))$, we can achieve this lower bound at points $x_0\in\Dcal_u\cup\Ccal_u$ by working in small cylinders $B(x_0,r)\times B(u(x_0),R)$ where $x_0$ is such that $g(x_0,u(x_0))>0$. At points where $g(x_0,u(x_0))=0$, we have that $f(x_0,u(x_0),A)=0$ for all $A\in\R^{m\times d}$ and so the lower bound is trivial. In this setting, the partial coercivity estimate $f(x,y,A)\geq g(x,y)|A|$ becomes a full coercivity estimate of the form $f(x,y,A)\geq c|A|$ and we can upgrade strong $\Lp^1$-convergence to weak* convergence in $\BV$ in order to make use of the results in~\cite{RindlerShaw19}.

As a consequence, the monotonicity and compatibility conditions~\eqref{prelimgbound} and~\eqref{prelimgcompatible} of Definition~\ref{defrepresentationfl1} are not needed to compute a lower bound at points $x_0\in\Dcal_u\cup\Ccal_u$: Proposition~\ref{propl1lebesgueineq} holds for any quasiconvex $f\in\mbfR(\Omega\times\R^m)$ (that is, $f$ Carath\'eodory and such that $f^\infty$ exists) \textcolor{black}{ that satisfies} $g(x,y)|A|\leq f(x,y,A)\leq Cg(x,y)(1+|A|)$ for \textcolor{black}{some} $g\in\C(\overline{\Omega}\times\R^m;[0,\infty))$\textcolor{black}{; no further conditions are required}. 

For points $x_0\in\Jcal_u$, the value $u(x_0)$ is not well-defined and minimising sequences must approximate both the inner and outer jump limits $u^+(x_0)$ and $u^-(x_0)$ near $x_0$. As a consequence, we cannot assume that minimising sequences only take values in regions where $f$ is non-degenerate, and an upgrade from partial to full coercivity is not possible. We can only work in the setting of $\Lp^1$-convergence which necessitates the use of stronger hypotheses on $g$ and $f$.

Let $f\in\Rbf(\Omega\times\R^m)$ satisfy $f\geq 0$ and \RED{consider a sequence $(u_j)_j\subset(\C^\infty\cap\W^{1,1})(\Omega;\R^m)$ such that $u_j\to u$ in $\Lp^1(\Omega;\R^m)$ for some $u\in\BV(\Omega;\R^m)$} and
\[
\liminf_{j\to\infty}\Fcal[u_j]<\infty.
\]
\RED{The restriction that $(u_j)_j\subset\C^\infty(\Omega;\R^m)$ rather than $(u_j)_j\subset\BV(\Omega;\R^m)$ will be lifted later for the purposes of Theorem~\ref{L1lscthm} using Corollary~\ref{corl1approx}.} Passing to a non-relabelled subsequence, we can assume that $\liminf_j\Fcal[u_j]=\lim_{j\to\infty}\Fcal[u_j]$ and hence (upon passing to another non-relabelled subsequence) that there exists a Radon measure $\mu\in\mbfM(\overline{\Omega})$ such that
\begin{equation*}
\wstarlim_{j\to\infty}f(x,u_j(x),\nabla u_j(x))\lL\restrict\Omega=\mu\quad\text{ in }\mbfM(\overline{\Omega}).
\end{equation*}
Using the Radon--Nikod\'{y}m Differentiation Theorem, we can write $\mu$ as the sum of mutually singular measures,
\[
\mu=\frac{\dd\mu}{\dd\lL}\lL\restrict\Omega+\frac{\dd\mu}{\dd|D^c u|}|D^c u|+\frac{\dd\mu}{\dd\Hcal^{d-1}\restrict\Jcal_u}\Hcal^{d-1}\restrict\Jcal_u+\mu^s.
\]
To obtain the lower semicontinuity statement
\begin{align*}
\liminf_{j\to\infty}\Fcal[u_j]&\geq\int_\Omega f(x,u(x),\nabla u(x))\;\dd x+\int_\Omega f^\infty\left(x,u(x),\frac{\dd D^cu}{\dd|D^cu|}(x)\right)\;\dd|D^cu|(x)\\
&\qquad+\int_{\Jcal_u}H_f[u](x)\;\dd\Hcal^{d-1}(x),
\end{align*}
it therefore suffices to prove the three pointwise inequalities
\begin{align}
\begin{split}\label{eqradonnikodymineqs}
&\frac{\dd\mu}{\dd\lL}(x)\geq f(x,u(x),\nabla u(x))\quad\text{ for }\lL\text{-almost every }x\in\Omega,\\
&\frac{\dd\mu}{\dd|D^c u|}(x)\geq f^\infty\left(x,u(x),\frac{\dd D^cu}{\dd|D^cu|}(x)\right)\quad\text{ for }|D^cu|\text{-almost every }x\in\Omega,\\
&\frac{\dd\mu}{\dd\Hcal^{d-1}\restrict\Jcal_u}(x)\geq H_f[u](x)\quad\text{ for }\Hcal^{d-1}\text{-almost every }x\in\Jcal_u.
\end{split}
\end{align}

Proposition~\ref{propl1lebesgueineq} below establishes the first two inequalities in~\eqref{eqradonnikodymineqs} via a unified argument and can be seen as the primary technical contribution of this paper.

\begin{proposition}\label{propl1lebesgueineq}
Let $f\in\Rbf(\Omega\times\R^m)$ satisfy a bound of the form 
\[
g(x,y)|A|\leq f(x,y,A)\leq Cg(x,y)(1+|A|),\qquad g\in\C(\overline{\Omega}\times\R^m;[0,\infty)),\quad C>0,
\]
and assume that $f(x,y,\frarg)$ is quasiconvex for all $(x,y)\in\overline{\Omega}\times\R^m$. Then, if $u\in\BV(\Omega;\R^m)$ and $(u_j)_j\subset(\C^\infty\cap\W^{1,1})({\Omega};\R^m)$ are such that
\[
\RED{u_j\to u\text{ in }\Lp^1(\Omega;\R^m)},\qquad f(x,u_j(x),\nabla u_j(x))\lL\restrict\Omega\wsc\mu\text{ in }\mbfM^+(\overline{\Omega}),
\]
it holds that
\[
\frac{\dd\mu}{\dd\lL}(x)\geq f(x,u(x),\nabla u(x))\text{ for }\lL\text{-almost every $x\in\Omega$,}
\]
and
\[
\frac{\dd\mu}{\dd|D^c u|}(x)\geq f^\infty\left(x,u(x),\frac{\dd D^cu}{\dd|D^cu|}(x)\right)\text{ for }|D^cu|\text{-almost every $x\in\Omega$.}
\]
\end{proposition}

\RED{The bulk of the proof of Proposition~\ref{propl1lebesgueineq} is contained in the following lemma, which tackles the special case where we assume that $(x,u(x))$ is always close to a fixed point of coercivity for $g$.}

\RED{
\begin{lemma} \label{lem:special}
Assume that $|\nabla u(x)|>0$ for $\lL$-almost every $x\in\Omega$ and assume also that there exist $y_0\in\R^m$ and $R$, $\delta>0$ such that
\[
g(x,y)\geq\delta\quad\text{for all }x\in\Omega\text{ and }y\in\R^m\text{ with }|y-y_0|\leq R.
\]
Then, the conclusions of Proposition~\ref{propl1lebesgueineq} hold for $(\lL+|D^cu|)$-almost every $x\in\Omega$ such that $|u(x)-y_0|<R$.
\end{lemma}}

\RED{Before proving Lemma~\ref{lem:special}, we first give an informal and simplified description of the most important points in order to ensure that the main ideas are clear. In particular, we ignore the need to ensure that various measure theoretic `almost everywhere' properties are satisfied, some measure theoretic approximations, and the matter of passing to subsequences in order to ensure convergence.

In order to be able to apply Theorem~\ref{thmwsclsc} to $\lim_{j\to\infty}\int_{B(x,r)}f(\overline{x},u_j(\overline{x}),\nabla u_j(\overline{x}))\;\dd \overline{x}$ for $r$ sufficiently small, we must modify the sequence $(u_j)_j$ to obtain a new weakly* convergent sequence. In Step~1 of the proof, we construct a truncated sequence $(\widetilde{u}_j)_j$ given by
\[
\widetilde{u}_j:=\mathbbm{1}_{A_j}u_j + \mathbbm{1}_{\Omega\setminus A_j}y_1
\]
for appropriately chosen sets of finite perimeter $A_j$. This truncated sequence converges weakly* in $\BV(\Omega;\R^m)$ to a limit $\widetilde{u}$ satisfying $\widetilde{u}=u$ in $\{\overline{x} : |u(\overline{x})-y_0|<R\}$.

For $x\in\{\overline{x} : |u(\overline{x})-y_0|<R\}$, the natural strategy might now seem to be to bound $\lim_{j\to\infty}\int_{B(x,r)}f(\overline{x},u_j(\overline{x}),\nabla u_j(\overline{x}))\;\dd \overline{x}$ from below by $\lim_{j\to\infty}\int_{B(x,r)}f(\overline{x},\widetilde{u}_j(\overline{x}),\nabla \widetilde{u}_j(\overline{x}))\;\dd \overline{x}$ up to an error term that goes to $0$ as $r\downarrow 0$ and then to use the fact that $\widetilde{u}(x)=u(x)$. However since the $\widetilde{u}_j$ are discontinuous, Theorem~\ref{thmwsclsc} can only be applied to the full functional
\[
\widetilde{u}_j\mapsto \int_{B(x,r)}f(\overline{x},\widetilde{u}_j(\overline{x}),\nabla \widetilde{u}_j(\overline{x}))\;\dd \overline{x} +  \int_{B(x,r)}\int_0^1 f^\infty\left(\overline{x},\widetilde{u}^\theta_j(\overline{x}),\frac{\dd D^s\widetilde{u}_j}{\dd|D^s\widetilde{u}_j|}(\overline{x})\right)\;\dd |D^s\widetilde{u}_j|(\overline{x}),
\]
which precludes this angle of attack, since the second term in the sum above cannot in general be expected to converge to $0$ as $j\to\infty$ and $r\to 0$. Instead, in Step~2 we obtain
\begin{align}\label{eqjustanotherbound}
\begin{split}
&\lim_{r\to 0}\lim_{j\to\infty}\frac{c_r}{r^d}\int_{B(x,r)}f(\overline{x},u_j(\overline{x}),\nabla u_j(\overline{x}))\;\dd \overline{x}\\
\geq&\lim_{r\to 0}\lim_{j\to\infty}\int_{\bB^d}\chi(\widetilde{w}^r_j(z))f_r(z,\Phi[\widetilde{w}^r_j](z),\nabla\Phi[\widetilde{w}^r_j](z))\;\dd z 
\end{split}
\end{align}
where $f_r$ is a localisation of $f$ defined below in~\eqref{defl1f_r} and $\Phi[\widetilde{w}^r_j]:=\varphi\circ\widetilde{w}^r_j$ with cut-off functions $\varphi$, $\chi$, and $\widetilde{w}^r_j=r^{-1}c_r(\widetilde{u}_j(x+r\frarg)-(\widetilde{u})_{x,r})$. 

In Step~3, we use the theory of liftings (Theorem~\ref{thmdiffusetangentlifting}) together with the graphical Besicovitch Derivation Theorem (more precisely, Lemma~\ref{lemgraphicallebesgueconverg}) to show that we can replace 
\[ 
 \lim_{r\to 0}\lim_{j\to\infty}\int_{\bB^d}\chi(\widetilde{w}^r_j(z))f_r(z,\Phi[\widetilde{w}^r_j](z),\nabla\Phi[\widetilde{w}^r_j](z))\;\dd z
\]
by
\[
\lim_{r\to 0}\lim_{j\to\infty}\int_{\bB^d}f_r(z,\Phi[\widetilde{w}^r_j](z),\nabla\Phi[\widetilde{w}^r_j](z))\;\dd z
\]
in~\eqref{eqjustanotherbound}. This is crucial, since it turns out that for $(\lL+|D^c\widetilde{u}|)$-almost every $x$, $\Phi[\widetilde{w}^r_j]$ is a $\C^\infty(\Omega;\R^m)$-function and so in Step~4 we are able to apply Theorem~\ref{thmwsclsc} as $j\to\infty$ to the functional
\[
\Fcal_r\colon\varphi\circ\widetilde{w}^r_j\mapsto \int_{\bB^d}f_r(z,\Phi[\widetilde{w}^r_j](z),\nabla\Phi[\widetilde{w}^r_j](z))\;\dd z
\]
before then letting $r\downarrow 0$ to complete the blow-up procedure. 

Finally, in Step~5, we show that the results of the previous steps obtained with $\widetilde{u}$ as the limit of $(\widetilde{u}_j)_j$ suffice to imply the conclusion of the lemma for $u$.
}

We will \RED{also} need the following lemma, which is due to Fonseca \& M\"uller and was first proved in~\cite{FonMul92QCIL} as Lemma~2.8:
\begin{lemma}\label{lemfonmul}
Let $v\in\C^\infty(\overline{\Omega};\R^m)$ and $0<a<b$. Then it follows that, for some constant $C>0$ independent of $v$,
\[
\operatorname{essinf}_{t\in[a,b]}\left(t\cdot\Hcal^{d-1}\left\{x\in\overline{\Omega}\colon |v(x)|=t\right\}\right)\leq\frac{C}{\ln (b/a)}\int_{\{|v(z)|\leq b\}}|\nabla v(z)|\;\dd z.
\]
\end{lemma}

\begin{proof}[Proof of Lemma~\ref{lem:special}]

\emph{Step 1:}
First note that
\begin{align*}
\limsup_{j\to\infty} \int_{\setsmall{ x \in \Omega}{|u_j(x)-y_0|< R}}|\nabla u_j(x)|\;\dd x & < \limsup_{j\to\infty} \frac{1}{\delta}\int_\Omega g(x,u_j(x))|\nabla u_j(x)|\;\dd x\\
&\leq\lim_{j\to\infty} \frac{1}{\delta} \Fcal[u_j]<\infty
\end{align*}
and fix $\tau\in(0,1)$. By Lemma~\ref{lemfonmul}, there exists a sequence $t_j\in[\RED{\frac{1}{2}(1+\tau)R,\frac{1}{4}(3+\tau)R}]$ such that
\begin{align*}
&\sup_j t_j\Hcal^{d-1}\Big(\Big\{x\in \Omega\colon|u_j(x)-y_0|=t_j\Big\}\Big)\\
&\qquad\leq\frac{C}{\log{\RED{\frac{3+\tau}{2(1+\tau)}}}}\sup_j\int_{\setsmall{x\in\Omega}{|u_j(x)-y_0|\leq \RED{\frac{1}{4}(3+\tau)R}}}|\nabla u_j(x)|\;\dd x\\
&\qquad\leq\frac{C}{\log{\RED{\frac{3+\tau}{2(1+\tau)}}}}\sup_j\int_{\setsmall{x\in\Omega}{|u_j(x)-y_0|< R}}|\nabla u_j(x)|\;\dd x\\
&\qquad<\infty.
\end{align*}
Letting $y_1\in\R^m$ be a fixed vector satisfying $|y_1-y_0|> R$, it follows that the sequence $(\widetilde{u}_j)_j\subset\BV(\Omega;\R^m)$ defined by
\begin{equation}\label{eq:defujtilde}
\widetilde{u}_j:=\mathbbm{1}_{\{|u_j(x)-y_0|< t_j\}}u_j + \mathbbm{1}_{\{|u_j(x)-y_0|\geq t_j\}}y_1
\end{equation}
is uniformly bounded in $\BV(\Omega;\R^m)$. Passing to a (non-relabelled) subsequence and using the fact that $(t_j)_j$ is bounded, we can assume that $t_j$ converges to a limit $t\in[\RED{\frac{1}{2}(1+\tau)R,\frac{1}{4}(3+\tau)R}]$. It follows then that $\widetilde{u}_j\wsc\widetilde{u}$ for some $\widetilde{u}\in\BV(\Omega;\R^m)$ of the form 
\begin{equation}\label{eq:defutilde}
\widetilde{u}=\mathbbm{1}_A u+\mathbbm{1}_{\Omega\setminus A}y_1,
\end{equation}
where $A\subset\Omega$ is a set of finite perimeter satisfying
\[
\{x\in \Omega\colon|u(x)-y_0|< \RED{\tau R}\}\subset A\subset\{x\in \Omega\colon|u(x)-y_0| \RED{<R}\}.
\]
This follows from the fact that, if $x$ is such that $|u(x)-y_0|<\RED{\tau R}$ and $u_j(x)\to u(x)$ then (since \RED{each $t_j >\tau R$}) eventually $|u_j(x)-y_0|<t_j$ (and hence $\widetilde{u}_j(x)=u_j(x)\to u(x)$) and similarly that if $x$ is such that $|u(x)-y_0|\RED{\geq R}$ and $u_j(x)\to u(x)$ then eventually $|u_j(x)-y_0|>t_j$ (at which point $\widetilde{u}_j(x)=y_1$).

We will show below that
\begin{equation}\label{equtildelebesguebound}
\frac{\dd\mu}{\dd\lL}(x)\geq f(x,\widetilde{u}(x),\nabla \widetilde{u}(x))\quad\text{ for }\lL\text{-a.e.\ }x\in \{|\widetilde{u}-y_0|<\tau R\}
\end{equation}
and
\begin{equation}\label{equtildecantorbound}
\frac{\dd\mu}{\dd|D^c\widetilde{u}|}(x)\geq f^\infty\left(x,\widetilde{u}(x),\frac{\dd D^c\widetilde{u}}{\dd|D^c\widetilde{u}|}(x)\right)\quad\text{ for }|D^c\widetilde{u}|\text{-a.e.\ }x\in\{|\widetilde{u}-y_0|<\tau R\}.
\end{equation}
\RED{The fact that $\widetilde{u}=u$ in $A$ will then allow us to replace $\widetilde{u}$ with $u$  in these inequalities to finish the proof.}

\emph{Step 2:} In the following, let 
\begin{equation}\label{defrho}
\rho = \omega_d^{-1}\lL\restrict\Omega  \qquad\text{or}\qquad
\rho = |D\widetilde{u}|.
\end{equation}
Pass to a (non-relabelled) subsequence so that
\begin{equation} \label{eq:gruj_eta}
  \gr^{\widetilde{u}_j}_\#(|\nabla\widetilde{u}_j|\lL\restrict\Omega)\wsc\eta  \qquad\text{as $j\to\infty$}
\end{equation}
for some $\eta\in\mbfM^+(\Omega\times\R^m)$. We abbreviate
\begin{equation} \label{eq:lambda_def}
\lambda:=|\gamma\asc{\widetilde{u}}|\restrict((\Omega\setminus\Jcal_{\widetilde{u}})\times\R^m)
\RED{=\gr^{\widetilde{u}}_\#(|\nabla \widetilde{u}|\lL\restrict\Omega + |D^c \widetilde{u}|)},
\end{equation}
\RED{where $\gamma[\widetilde{u}]$ is the elementary lifting associated to $\widetilde{u}$, introduced in Definition~\ref{defelementaryliftings},} and let
\[
\eta=\frac{\dd\eta}{\dd\lambda}\lambda+\eta^s\quad\text{ where }\quad\lambda\perp\eta^s
\]
be the Radon--Nikodym decomposition of $\eta$ with respect to $\lambda$.

Now define $\Bcal\subset\{x\in A\colon|\widetilde{u}(x)-y_0|<\tau R\}$ to be the set of points which additionally satisfy

\begin{equation}\label{eqlocalcond0} 
x\in\Dcal_{\widetilde{u}}\quad \text{ if }\quad\rho=\omega_d^{-1}\lL,\quad\text{ or }\quad x\in\Ccal_{\widetilde{u}} \quad\text{ if }\quad\rho=|D\widetilde{u}|;
\end{equation}
\begin{equation}\label{eqlocalcond0.5}
\lim_{r\to 0}\frac{r^d}{\rho(B(x,r))}=\omega_d^{-1}\frac{\dd\lL}{\dd\rho}(x)\quad\text{ and }\quad\displaystyle\lim_{r\to 0}\frac{\pi_\#\lambda(B(x,r))}{\rho(B(x,r))}=\frac{\dd\pi_\#\lambda}{\dd\rho}(x);
\end{equation}
\begin{align} \label{eqlocalcond2} 
\left\{\begin{aligned}
&\text{there exists a sequence $r_n\downarrow 0$ and a limit function ${\widetilde{u}^0}\in(\BV\cap\Lp^\infty)(\bB^d;\R^m)$}\\
&\text{of the form described in Theorem~\ref{thmbvblowup} such that}\\
&\RED{\widetilde{u}^{n}}\to \widetilde{u}^0 \quad\text{ strictly in }\quad (\BV\cap\Lp^\infty)(\bB^d;\R^m),\quad\text{where}\\
&\RED{\widetilde{u}^{n}:= c_nr_n^{-1}(\widetilde{u}(x+r_n\frarg)-(\widetilde{u})_{x,r_n})}\quad\text{and}\quad\RED{c_n^{-1} r_n\to 0};
\end{aligned}\right.
\end{align}
\begin{equation}\label{eqlocalcond1}
\text{$x$ is a Lebesgue point for $\widetilde{u}$}\quad\text{and}\quad\lim_{r\to 0}\frac{|D\widetilde{u}|(\Jcal_{\widetilde{u}}\cap B(x,r))}{|D\widetilde{u}|(B(x,r))} =0;
\end{equation}
\begin{align}\label{eqlocalcond1.5}
\left\{\begin{aligned}
\text{$x$ is a cylindrical $\lambda$-Lebesgue point for $\displaystyle \frac{\dd\eta}{\dd\lambda}$ in the sense of Theorem~\ref{thmgeneralisedbesicovitch} and}\\
\lim_{n\to\infty}\frac{\eta^s\left(\overline{B(x,r_n)}\times \overline{B((\widetilde{u})_{x,r_n},(1+\norm{\widetilde{u}^0}_\infty)c_n^{-1} r_n)}\right)}{\pi_\#\lambda({B(x,r_n)})}=0;
\end{aligned}\right.
\end{align}
\begin{align}\label{eqlocalcond2.75}
\left\{\begin{aligned}
&\text{if $\rho=\omega_d^{-1}\lL$, then $x$ is a $\lL$-Lebesgue point of the function $f(\frarg,\widetilde{u}(\frarg),\nabla\widetilde{u}(\frarg))$, or}\\
&\text{if  $\rho=|D\widetilde{u}|$, then}\quad\lim_{n\to\infty}\frac{(|\nabla \widetilde{u}|\lL+|D^j\widetilde{u}|)(B(x,r_n))}{|D\widetilde{u}|(B(x,r_n))}=0;
\end{aligned}\right.
\end{align}

Condition~\eqref{eqlocalcond0} holds for $\rho$-almost every $x\in A$ thanks to the definitions of $\Dcal_{\widetilde{u}}$ and $\Ccal_{\widetilde{u}}$, whilst Condition~\eqref{eqlocalcond0.5} holds for $\rho$-almost every $x\in A$ thanks to the Besicovitch Differentiation Theorem. Condition~\eqref{eqlocalcond2} is satisfied at $\rho$-almost every $x\in\Dcal_{\widetilde{u}}\cup\Ccal_{\widetilde{u}}$ as a consequence of Theorem~\ref{thmbvblowup} and Proposition~3.92 in~\cite{AmFuPa00FBVF}. That Condition~\eqref{eqlocalcond1} holds $\rho$-almost everywhere in $A$ stems from the fact that $\Hcal^{d-1}$-almost every $x\in\Omega\setminus\Jcal_{\widetilde{u}}$ (and hence $\rho$-almost every $x\in\Dcal_{\widetilde{u}}\cup\Ccal_{\widetilde{u}}$ since $\Jcal_{\widetilde{u}}\cap (\Dcal_{\widetilde{u}}\cup\Ccal_{\widetilde{u}})=\emptyset$) is a Lebesgue point of $\widetilde{u}$ and the fact that $(\Dcal_{\widetilde{u}}\cup\Ccal_{\widetilde{u}})\cap\Jcal_{\widetilde{u}}=\emptyset$ together with the Besicovitch Differentiation Theorem. Condition~\eqref{eqlocalcond1.5} holds for $|\lambda|=|\nabla \widetilde{u}|\lL+|D^c\widetilde{u}|$-almost every $x\in\Omega$ and therefore for $|D\widetilde{u}|$-almost every $x\in\Dcal_{\widetilde{u}}\cup\Ccal_{\widetilde{u}}$. Since $\nabla u=\nabla\widetilde{u}$ $\lL$-almost everywhere in $A$ and we have assumed that $|\nabla u(x)|>0$ for $\lL$-almost every $x\in\Omega$, it holds that $\lL\restrict A \ll |D\widetilde{u}|$ and so Condition~\eqref{eqlocalcond1.5} holds for $\rho$-almost every $x$ in $A$. That Condition~\eqref{eqlocalcond2.75} is satisfied for $\rho$-almost every $x\in\Dcal_{\widetilde{u}}\cup\Ccal_{\widetilde{u}}$ follows from the Lebesgue Differentiation Theorem and the Besicovitch Differentiation Theorem combined with the fact that $(|\nabla\widetilde{u}|\lL+|D^j\widetilde{u}|)(\Ccal_{\widetilde{u}})=0$. It follows then that $\rho(\{|\widetilde{u}-y_0|<\tau R\}\setminus\Bcal)=0$.

For a fixed $x\in\Bcal$, we can now consider
\begin{equation*}
\frac{\dd\mu}{\dd\rho}(x)=\lim_{n\to\infty}\frac{\mu({B(x,r_n)})}{\rho(B(x,r_n))}\geq\lim_{n\to\infty}\frac{1}{\rho(B(x,r_n))}\lim_{j\to\infty}\int_{B(x,r_n)}f(y,u_j(y),\nabla u_j(y))\;\dd y.
\end{equation*}
Rewrite $u_j(x+r_nz)=(\widetilde{u})_{x,r_n}+c_n^{-1} r_nw^{n}_j(z)$, where 
\[
w^n_j(z):=c_n\frac{u_j(x+r_nz)-(\widetilde{u})_{x,r_n}}{r_n},\quad c_n:=\frac{r_n^d}{\rho(B(x,r_n))},
\]
so that
\begin{align*}
\frac{1}{\rho(B(x,r_n))}&\int_{B(x,r_n)}f(y,u_j(y),\nabla u_j(y))\;\dd y\\
=&\int_{\bB^d}c_nf(x+r_nz,(\widetilde{u})_{x,r_n}+c_n^{-1} r_nw^n_j(z),c_n^{-1}\nabla w^n_j(z))\;\dd z\\
=&\int_{\bB^d}f_n(z,w^n_j(z),\nabla w^n_j(z))\;\dd z,
\end{align*}
where $f_n:=f_{r_n}$ and $f_r$ is given for each $r>0$ by
\begin{equation}\label{defl1f_r}
f_r(z,y,A):=c_r f\left(x+rz,(\widetilde{u})_{x,r}+c_r^{-1}ry,c_r^{-1}A\right).
\end{equation}
We can therefore write
\begin{equation}\label{eqrhobound}
\frac{\dd\mu}{\dd\rho}(x)\geq\liminf_{n\to\infty}\liminf_{j\to\infty}\int_{\bB^d}f_n(z,w^n_j(z),\nabla w^n_j(z))\;\dd z.
\end{equation}

Defining $(\widetilde{w}_j^n)_{n\in\mbN,j\in\mbN}\subset\BV(\bB^d;\R^m)$ by 
\[
\widetilde{w}_j^n(z)=c_n\frac{\widetilde{u}_j(x+r_nz)-(\widetilde{u})_{x,r_n}}{r_n},
\]
\RED{with $\widetilde{u}$ defined in~\eqref{eq:defutilde}}, we claim that, once $n$ is sufficiently large,
\begin{equation}\label{eql1blowupa}
|\widetilde{w}_j^n(z)|\leq 2+\norm{\widetilde{u}^0}_\infty\quad\text{ implies }\quad\widetilde{w}_j^n(z)=w_j^n(z)\quad\text{ for all }z\in\bB^d,\text{ }j\in\mbN,
\end{equation}
and
\begin{equation}\label{eql1blowupb}
|\widetilde{w}_j^n(z)|> 2+\norm{\widetilde{u}^0}_\infty\quad\text{ implies }\quad|{w}_j^n(z)|> 2+\norm{\widetilde{u}^0}_\infty\quad\text{ for all }z\in\bB^d,\text{ }j\in\mbN.
\end{equation}

To see this, note first that, for any $M>0$,
\[
\left|c_r\frac{\widetilde{u}_j(x+rz)-(\widetilde{u})_{x,r}}{r}\right|\leq M
\]
implies
\begin{equation*}
|\widetilde{u}_j(x+rz)-y_0|  \leq |y_0-\widetilde{u}(x)|+|\widetilde{u}(x)-(\widetilde{u})_{x,r}|+Mc_r^{-1}r.
\end{equation*}
By Conditions~\eqref{eqlocalcond2} and~\eqref{eqlocalcond1} combined with the fact that $|y_0-\widetilde{u}(x)|<\tau R$, we therefore have that, once $n$ is sufficiently large,
\begin{equation}\label{l1blowupeq3}
|\widetilde{w}_j^n(z)|\leq M \quad\text{ implies }\quad|\widetilde{u}_j(x+r_nz)-y_0|<\tau R\leq t_j.
\end{equation}
Since $|y_1-y_0|> R$, this is only possible if $\widetilde{u}_j(x+r_nz)=u_j(x+r_nz)$ which implies $\widetilde{w}_j^n(z)=w_j^n(z)$ and hence that~\eqref{eql1blowupa} holds as required. \RED{This implies, in particular, that $z \notin \Jcal_{\widetilde{w}_j^n}$.}

On the other hand, assume that $|{w}_j^n(z)|\leq M$. Repeating the preceding calculation with $\widetilde{u}_j(x+r_nz)$ replaced by $u_j(x+r_nz)$, we deduce that $|{u}_j(x+r_nz)-y_0|<\tau R\leq t_j$. This again implies $u_j(x+r_nz)=\widetilde{u}_j(x+r_nz)$ and hence that $\widetilde{w}_j^n(z)=w_j^n(z)$. In particular, we have that $|{w}_j^n(z)|\leq M$ implies $|\widetilde{w}_j^n(z)|\leq M$ which, upon taking $M=2+\norm{\widetilde{u}^0}_\infty$, is the contrapositive of~\eqref{eql1blowupb}. By passing to a tail of the sequence $(r_n)_{n}$ if necessary, we can therefore assume that~\eqref{eql1blowupa} and~\eqref{eql1blowupb} always hold.

Now let $\varphi\in\C_c^\infty(\R^m;\R^m)$ be such that
\[
\varphi(y)=\begin{cases}
y \quad\text{ for } &|y|\leq \norm{\widetilde{u}^0}_\infty+1,\\
0 \quad \text{ for }&|y|\geq \norm{\widetilde{u}^0}_\infty+2,
\end{cases}
\]
and define the nonlinear truncation operator $\Phi\colon\BV(\bB^d;\R^m)\to(\BV\cap\Lp^\infty)(\bB^d;\R^m)$ by
\begin{equation*}
\Phi[w]:=\varphi\circ w,\qquad w\in\BV(\bB^d;\R^m).
\end{equation*}
Let also $\chi\in\C_c(\R^m;[0,1])$ be such that $\chi(y)=1$ for $|y|\leq\norm{\widetilde{u}^0}_\infty$ and $\chi(y)=0$ for $|y|\geq\norm{\widetilde{u}^0}_\infty+1$. By~\eqref{eql1blowupa}, we have that $w_j^n=\widetilde{w}_j^n$ whenever $\chi(\widetilde{w}_j^n)>0$. Since our choice of $\varphi$ also implies that $\widetilde{w}_j^n=\Phi[\widetilde{w}_j^n]$ whenever $\chi(\widetilde{w}_j^n)> 0$, we can estimate
\begin{align}
\begin{split}\label{eqdisgustingintegralineq}
\int_{\bB^d}f_n(z,w^n_j(z),\nabla w^n_j(z))\;\dd z &\geq\int_{\bB^d}\chi(\widetilde{w}_j^n(z))f_n(z,w^n_j(z),\nabla w^n_j(z))\;\dd z\\
&=\int_{\bB^d}\chi(\widetilde{w}_j^n(z))f_n(z,\Phi[\widetilde{w}^n_j](z),\nabla \Phi[\widetilde{w}^n_j](z))\;\dd z
\end{split}
\end{align}
for all $j\in\mbN$. 

\emph{Step 3:} By~\eqref{eqrhobound} together with~\eqref{eqdisgustingintegralineq}, we have
\begin{align*}
\frac{\dd\mu}{\dd\rho}(x)&\geq\liminf_{n\to\infty}\liminf_{j\to\infty}\int_{\bB^d}f_n(z,\Phi[\widetilde{w}^n_j](z),\nabla \Phi[\widetilde{w}^n_j](z))\;\dd z\\
&\qquad-\limsup_{n\to\infty}\limsup_{j\to\infty}\int_{\bB^d}(1-\chi(\widetilde{w}_j^n(z)))f_n(z,\Phi[\widetilde{w}^n_j](z),\nabla \Phi[\widetilde{w}^n_j](z))\;\dd z.
\end{align*}
In what follows, we will first use Lemma~\ref{lemgraphicallebesgueconverg} to show that the second term on the right hand side of the inequality above is equal to zero before applying Theorem~\ref{thmwsclsc} to the sequence $(\Phi[\widetilde{w}_j^n])_j$ and the integrand $f_n$ to provide a lower bound on the first term.

By virtue of the definition~\eqref{defl1f_r} of $f_n$ together with the growth assumptions on $f$ in the statement of the proposition, there exists a constant $C>0$ depending only on $\norm{\widetilde{u}^0}_\infty$ such that
\[
|f_n(z,\varphi(y),\nabla\varphi(y)A)|\leq C \left(c_n+|\nabla\varphi(y)A|\right).
\]
Defining $\psi\in\C_0(\R^m)$ by
\begin{equation*}
\psi(w):=C(1-\chi(w))|\nabla\varphi(w)|,
\end{equation*}
we see that
\begin{align*}
&\limsup_{n\to\infty}\limsup_{j\to\infty}\int_{\bB^d} (1-\chi(\widetilde{w}_j^n(z)))f_n(z,\Phi[\widetilde{w}^n_j](z),\nabla \Phi[\widetilde{w}^n_j](z))\;\dd z\\
&\leq\lim_{n\to\infty}c_n\lim_{j\to\infty}C\int_{\bB^d}	1-\chi(\widetilde{w}_j^n(z))\;\dd z+\lim_{n\to\infty}\lim_{j\to\infty}\int_{\bB^d}\psi(\widetilde{w}_j^n(z))|\nabla\widetilde{w}_j^n(z)|\;\dd z.
\end{align*}
As $j\to\infty$, we have that $\widetilde{w}_j^n\to\widetilde{u}^n$ strongly in $\Lp^1(\bB^d;\R^m)$ (where $\widetilde{u}^n=c_nr_n^{-1}(\widetilde{u}(x+r_n\frarg)-(\widetilde{u})_{x,r_n})$ is as defined in Condition~\eqref{eqlocalcond2})  and as $n\to\infty$ we have that $\widetilde{u}^n\to\widetilde{u}^0$ in $\Lp^1(\bB^d;\R^m)$. Thus, since $\chi(y)=1$ whenever $|y|\leq\norm{\widetilde{u}^0}_\infty$ and Condition~\eqref{eqlocalcond0.5} holds,
\[
\lim_{n\to\infty}c_n\lim_{j\to\infty}C\int_{\bB^d} 1-\chi(\widetilde{w}_j^n(z))\;\dd z=\lim_{n\to\infty}c_n\int_{\bB^d} 1-\chi(\widetilde{u}^n(z))\;\dd z = \omega_d^{-1}\frac{\dd\lL}{\dd\rho}(x)\cdot 0 = 0,
\]
and so
\begin{align*}
&\limsup_{n\to\infty}\limsup_{j\to\infty}\int_{\bB^d}(1-\chi(\widetilde{w}_j^n(z)))f_n(z,\Phi[\widetilde{w}^n_j](z),\nabla \Phi[\widetilde{w}^n_j](z))\;\dd z\\\leq &\lim_{n\to\infty}\lim_{j\to\infty}\frac{1}{\rho(B(x,r_n))}\int_{B(x,r_n)\times\R^m}\psi\left(c_n\frac{y-(\widetilde{u})_{x,r_n}}{r_n}\right)\;\dd\gr^{\widetilde{u}_j}_\#(|\nabla\widetilde{u}_j|\lL)(\overline{x},y)\\
\leq &\lim_{n\to\infty}\frac{1}{\rho(B({x},r_n))}\int_{\overline{B(x,r_n)}\times\R^m}\psi\left(c_n\frac{y-(\widetilde{u})_{x,r_n}}{r_n}\right)\;\dd\eta(\overline{x},y),
\end{align*}
where we used~\eqref{eq:gruj_eta} in the last step. By Conditions~\RED{\eqref{eqlocalcond0.5} and}~\eqref{eqlocalcond1.5} together with the fact that $\supp\psi\subset\overline{B^m(0,1+\norm{\widetilde{u}^0}_\infty)} \subset A$, we see that
\begin{align*}
&\lim_{n\to\infty}\frac{1}{\rho(B(x,r_n))}\Bigg|\int_{\overline{B(x,r_n)}\times\R^m}\psi\left(c_n\frac{y-(\widetilde{u})_{x,r_n}}{r_n}\right)\;\dd\eta^s(\overline{x},y)\Bigg|\\
\leq &\lim_{n\to\infty}\norm{\psi}_\infty\frac{\pi_\#\lambda({B(x,r_n)})}{\rho(B(x,r_n))}\cdot\frac{\eta^s\left(\overline{B(x,r_n)}\times \overline{B((\widetilde{u})_{x,r_n},(1+\norm{\widetilde{u}^0}_\infty)c_n^{-1} r_n)}\right)}{\pi_\#\lambda({B(x,r_n)})}\\
= & \frac{\dd\pi_\#\lambda}{\dd\rho}(x)\cdot 0=0.
\end{align*}
Hence,
\begin{align}
\begin{split}\label{eqerrorbound}
&\lim_{n\to\infty}\frac{1}{\rho(B(x,r_n))}\int_{\overline{B(x,r_n)}\times\R^m}\psi\left(c_n\frac{y-(\widetilde{u})_{x,r_n}}{r_n}\right)\;\dd\eta(\overline{x},y)\\
=&\lim_{n\to\infty}\frac{1}{\rho(B(x,r_n))}\int_{\overline{B(x,r_n)}\times\R^m}\psi\left(c_n\frac{y-(\widetilde{u})_{x,r_n}}{r_n}\right)\frac{\dd\eta}{\dd\lambda}(\overline{x},y)\;\dd\lambda(\overline{x},y)\\
=&\lim_{n\to\infty}\int_{\bB^d\times\R^m}\psi\left(w\right)\;\dd\left[\frac{1}{\rho(B(x,r_n))}\left(T^{(x,r_n), ((\widetilde{u})_{x,r_n}, c_n^{-1}r_n)}\right)_\#\biggl(\frac{\dd\eta}{\dd\lambda}\lambda\biggr)\right](z,w),
\end{split}
\end{align}
\RED{where in the last line we used that $\lambda(\partial B(x,r_n)\times\R^m) = 0$.} Now note that, since (\RED{by~\eqref{eqrescaledliftings} together with~\eqref{defrho}})
\[
\frac{1}{\rho(B(x,r_n))}(T^{(x,r_n), ((\widetilde{u})_{x,r_n}, c_n^{-1}r_n)})_\# \bigl( \gamma\asc{\widetilde{u}}\restrict B(x,r_n)\times\R^m \bigr)=\gamma\asc{\widetilde{u}^{r_n}},
\]
and \RED{Conditions~\eqref{eq:lambda_def}} and~\eqref{eqlocalcond1} hold, we have that
\begin{align*}
&\lim_{n\to\infty}\left|\frac{1}{\rho(B(x,r_n))}\left(T^{(x,r_n), ((\widetilde{u})_{x,r_n}, c_n^{-1}r_n)}\right)_\#\lambda-|\gamma\asc{\widetilde{u}^{n}}|\right|(\bB^d\times\R^m)\\
=&\lim_{n\to\infty}\left|\frac{1}{\rho(B(x,r_n))}\left(T^{(x,r_n), ((\widetilde{u})_{x,r_n}, c_n^{-1}r_n)}\right)_\# \bigl(|\gamma\asc{\widetilde{u}}|\restrict(\Jcal_{\widetilde{u}}\times\R^m) \bigr)\right|(\bB^d\times\R^m)\\
\leq & \; C \cdot \lim_{n\to\infty}\frac{|D\widetilde{u}|(\Jcal_{\widetilde{u}}\cap B(x,r_n))}{|D\widetilde{u}|(B(x,r_n))}\\
= &\; 0.
\end{align*}
Since $|\gamma\asc{\widetilde{u}^{n}}|\to|\gamma\asc{\widetilde{u}^0}|$ strictly in $\mbfM^+(\bB^d\times\R^m)$ as $n\to\infty$ \RED{(by Reshetnyak's Continuity Theorem a sequence of strictly converging measures must have strictly converging total variations)}, we therefore have that
\[
\frac{1}{\rho(B(x,r_n))}\left(T^{(x,r_n), ((\widetilde{u})_{x,r_n}, c_n^{-1}r_n)}\right)_\#\lambda\to|\gamma\asc{\widetilde{u}^0}|\quad\text{strictly in }\mbfM(\bB^d\times\R^m)\text{ as }n\to\infty.
\]
Combining Lemma~\ref{lemgraphicallebesgueconverg} with Condition~\eqref{eqlocalcond1.5}, it follows that
\[
\frac{1}{\rho(B(x,r_n))}\left(T^{(x,r_n), ((\widetilde{u})_{x,r_n}, c_n^{-1}r_n)}\right)_\#\biggl(\frac{\dd\eta}{\dd\lambda}\lambda\biggr)\to\frac{\dd\eta}{\dd\lambda}(x,\widetilde{u}(x))|\gamma\asc{\widetilde{u}^0}|
\]
strictly in $\mbfM(\bB^d\times\R^m)$ as $n\to\infty$. By~\eqref{eqerrorbound} together with the fact that $\psi(w)=0$ for all $w\in\R^m$ satisfying $|w|\leq\norm{\widetilde{u}^0}_\infty$, we can then deduce
\begin{align*}
&\lim_{n\to\infty}\frac{1}{\rho(B(x,r_n))}\int_{\overline{B(x,r_n)}\times\R^m}\psi\left(c_n\frac{y-(\widetilde{u})_{x,r_n}}{r_n}\right)\;\dd\eta(\overline{x},y)\\
=&\int_{\bB^d\times\R^m}\psi(w)\frac{\dd\eta}{\dd\lambda}(x,\widetilde{u}(x))\;\dd|\gamma\asc{\widetilde{u}^0}|(z,w)\\
=& \; 0.
\end{align*}
Hence,
\begin{align*}
&\limsup_{n\to\infty}\limsup_{j\to\infty}\int_{\bB^d}(1-\chi(\widetilde{w}_j^n(z)))f_n(z,\Phi[\widetilde{w}^n_j](z),\nabla \Phi[\widetilde{w}^n_j](z))\;\dd z\\
\leq&\lim_{n\to\infty}\frac{1}{\rho(B(x,r_n))}\int_{\overline{B(x,r_n)}\times\R^m}\psi\left(c_n\frac{y-(\widetilde{u})_{x,r_n}}{r_n}\right)\;\dd\eta(\overline{x},y)\\
= &\; 0,
\end{align*}
which at last leaves us with
\begin{equation}\label{eql1blowuplb1}
\frac{\dd\mu}{\dd\rho}(x)\geq\liminf_{n\to\infty}\liminf_{j\to\infty}\int_{\bB^d}f_n(z,\Phi[\widetilde{w}^n_j](z),\nabla \Phi[\widetilde{w}^n_j](z))\;\dd z.
\end{equation}

\emph{Step 4:} Next, we claim that $\Phi[\widetilde{w}^n_j]=\Phi[w^n_j]$ in $\BV(\bB^d;\R^m)$ for all $n$, $j\in\mbN$ and in particular that we always have $\Phi[\widetilde{w}^n_j]\in\C^\infty(\bB^d;\R^m)$. If $z\in\bB^d$ is such that $|\widetilde{w}_j^n(z)|\leq 2+\norm{\widetilde{u}^0}_\infty$ then~\eqref{eql1blowupa} implies that $\widetilde{w}_j^n(z)=w_j^n(z)$ and hence that $\Phi[\widetilde{w}^n_j](z)=\Phi[w^n_j](z)$ as required. On the other hand, if $|\widetilde{w}_j^n(z)|> 2+\norm{\widetilde{u}^0}_\infty$ then~\eqref{eql1blowupb} combined with the definition of $\varphi$ implies that both $\Phi[\widetilde{w}_j^n](z)$ and $\Phi[w_j^n](z)$ are equal to zero, which also leaves us with the desired result.

We have that $(\Phi[\widetilde{w}_j^n])_j$ is a sequence of smooth functions converging weakly* in $\BV(\bB^d;\R^m)$ as $j\to\infty$ to the limit $\Phi[\widetilde{u}^n]$. Letting $\zeta_R\in\C_c(\R^m;[0,1])$ be a cut-off function such that $\zeta_R(y)=1$ for $|y|\leq R$ we have that \RED{$(z,y,A)\mapsto \zeta_R(y)f_n(z,y,A)=: (\zeta_Rf_n)(z,y,A)$} satisfies the assumptions of Theorem~\ref{thmwsclsc} and so,
\begin{align*}
&\liminf_{j\to\infty}\int_{\bB^d} f_n(z,\Phi[\widetilde{w}^n_j](z),\nabla \Phi[\widetilde{w}^n_j](z))\;\dd z \\
 \geq &\liminf_{j\to\infty}\int_{\bB^d}(\zeta_Rf_n)(z,\Phi[\widetilde{w}^n_j](z),\nabla \Phi[\widetilde{w}^n_j](z))\;\dd z\\
\geq &\int_{\bB^d}(\zeta_R f_n)(z,\Phi[\widetilde{u}^n](z),\nabla\Phi[\widetilde{u}^n](z))\;\dd z\\
&\quad+\int_{\bB^d}(\zeta_R f_n)^\infty\left(z,\Phi[\widetilde{u}^n](z),\frac{\dd D^c\Phi[\widetilde{u}^n]}{\dd|D^c\Phi[\widetilde{u}^n]|}(z)\right)\;\dd |D^c\Phi[\widetilde{u}^n]|(z),
\end{align*}
where we have neglected the positive $\Hcal^{d-1}$-absolutely continuous term since it can only contribute mass. Letting $R\uparrow\infty$ so that $\zeta_R\uparrow 1$, the Monotone Convergence Theorem therefore lets us deduce
\begin{align*}
&\liminf_{j\to\infty}\int_{\bB^d} f_n(z,\Phi[\widetilde{w}^n_j](z),\nabla \Phi[\widetilde{w}^n_j](z))\;\dd z \\
&\geq\int_{\bB^d} f_n(z,\Phi[\widetilde{u}^n](z),\nabla\Phi[\widetilde{u}^n](z))\;\dd z+\int f_n^\infty\left(z,\Phi[\widetilde{u}^n](z),\frac{\dd D^c\Phi[\widetilde{u}^n]}{\dd|D^c\Phi[\widetilde{u}^n]|}(z)\right)\;\dd |D^c\Phi[\widetilde{u}^n]|(z).
\end{align*}
The Chain Rule in $\BV$ (see Theorem~3.96 in~\cite{AmFuPa00FBVF}) implies that
\[
\nabla \Phi[\widetilde{u}^n]=\nabla\varphi(\widetilde{u}^n)\nabla\widetilde{u}^n\quad\text{and}\quad D^c\Phi[\widetilde{u}^n]=\nabla\varphi(\widetilde{u}^n)D^c\widetilde{u}^n,
\]
and so, recalling that $\varphi(y)=y$ whenever $\chi(y) > 0$, we have
\[
\Phi[\widetilde{u}^n]=\widetilde{u}^n,\quad\nabla\Phi[\widetilde{u}^n]=\nabla\widetilde{u}^n\quad\text{and}\quad D^c\Phi[\widetilde{u}^n]=D^c\widetilde{u}^n\quad\text{in}\quad\{z\in\bB^d\colon\chi(\widetilde{u}^n(z))>0\}.
\]
Since $0\leq\chi\leq 1$, it then follows that
\begin{align}
\begin{split}\label{eql1blowuplb2}
&\liminf_{j\to\infty}\int_{\bB^d} f_n(z,\Phi[\widetilde{w}^n_j](z),\nabla \Phi[\widetilde{w}^n_j](z))\;\dd z \\
 &\geq \int_{\bB^d} (\chi f_n)(z,\widetilde{u}^n(z),\nabla\widetilde{u}^n(z))\;\dd z+\int_{\bB^d}(\chi f_n)^\infty\left(z,\widetilde{u}^n(z),\frac{\dd D^c\widetilde{u}^n}{\dd|D^c\widetilde{u}^n|}(z)\right)\;\dd|D^c\widetilde{u}^n|(z).
 \end{split}
\end{align}
By virtue of the definitions of $f_n$ and $\widetilde{u}^n$, however,
\begin{align}
\begin{split}\label{eql1blowupeq}
f_n(z,\widetilde{u}^n(z),\nabla\widetilde{u}^n(z))&=c_n f(x+r_nz,\widetilde{u}(x+r_nz),\nabla\widetilde{u}(x+r_nz)),\\
f^\infty_n\left(z,\widetilde{u}^n(z),\frac{\dd D^c\widetilde{u}^n}{\dd|D^c\widetilde{u}^n|}(z)\right)&=f^\infty\left(x+r_nz,\widetilde{u}(x+r_nz),\frac{\dd D^c\widetilde{u}^n}{\dd|D^c\widetilde{u}^n|}(z)\right).
\end{split}
\end{align}

\emph{Case 1:} Assume $x\in\Dcal_{\widetilde{u}}$, so that $\rho=\omega_d^{-1}\lL\restrict\Omega$ and \RED{$c_{n}=1$}. Combining~\eqref{eql1blowuplb1},~\eqref{eql1blowuplb2}, and~\eqref{eql1blowupeq}, we have that
\begin{equation*}
\frac{\dd\mu}{\dd\lL}(x)\geq\liminf_{n\to\infty}\dashint_{\bB^d}\chi\left(\widetilde{u}^{n}(z)\right)f(x+r_nz,\widetilde{u}(x+r_nz),\nabla\widetilde{u}(x+r_nz))\;\dd z.
\end{equation*}
By Condition~\eqref{eqlocalcond2.75}, we have that
\begin{align*}
&\left|\lim_{n\to\infty}\int_{\bB^d}\chi\left(\widetilde{u}^{n}(z)\right)\left[f(x+r_nz,\widetilde{u}(x+r_nz),\nabla\widetilde{u}(x+r_nz))-f(x,\widetilde{u}(x),\nabla\widetilde{u}(x))\right]\;\dd z\right| \\
\leq&\lim_{n\to\infty}\int_{\bB^d}|f(x+r_nz,\widetilde{u}(x+r_nz),\nabla\widetilde{u}(x+r_nz))-f(x,\widetilde{u}(x),\nabla\widetilde{u}(x))|\;\dd z=0,
\end{align*}
and so
\[
\frac{\dd\mu}{\dd\lL}(x)\geq\liminf_{n\to\infty}\dashint_{\bB^d}\chi\left(\widetilde{u}^{n}(z)\right)f(x,\widetilde{u}(x),\nabla\widetilde{u}(x))\;\dd z.
\]
Since $\widetilde{u}^{n}\to\widetilde{u}^0$ strongly in $\Lp^1(\bB^d;\R^m)$, this further simplifies to
\[
\frac{\dd\mu}{\dd\lL}(x)\geq\dashint_{\bB^d}\chi\left(\widetilde{u}^{0}(z)\right)f(x,\widetilde{u}(x),\nabla\widetilde{u}(x))\;\dd z.
\]
Since $\chi(y)=1$ whenever $|y|\leq\norm{\widetilde{u}^0}_\infty$, however, we have that $\chi\left(\widetilde{u}^{0}(z)\right)\equiv 1$ and hence that
\[
\frac{\dd\mu}{\dd\lL}(x)\geq\dashint_{\bB^d}f(x,\widetilde{u}(x),\nabla\widetilde{u}(x))\;\dd z=f(x,\widetilde{u}(x),\nabla\widetilde{u}(x)).
\]
\RED{Since $\rho(\{|\widetilde{u}-y_0|<\tau R\}\setminus\Bcal)=0$, this implies~\eqref{equtildelebesguebound}.}

\emph{Case 2:} Assume $x\in\Ccal_{\widetilde{u}}$, so that $\rho=|D\widetilde{u}|$. Combining~\eqref{eql1blowuplb1},~\eqref{eql1blowuplb2}, and~\eqref{eql1blowupeq}, we have that
\[
\frac{\dd\mu}{\dd|D\widetilde{u}|}(x)\geq\liminf_{n\to\infty}\int_{\bB^d}\chi(\widetilde{u}^{n}(z))f^\infty\left(x+r_nz,\widetilde{u}(x+r_nz),\frac{\dd D^c\widetilde{u}^n}{\dd|D^c\widetilde{u}^n|}(z)\right)\;\dd |D^c\widetilde{u}^n|(z).
\]
If $z\in\bB^d$ is such that $\chi(\widetilde{u}^{n}(z)) > 0$, then $|\widetilde{u}^n(z)|\leq\norm{\widetilde{u}^0}_\infty + 1$ which implies
\[
|\widetilde{u}(x+r_nz)-\widetilde{u}(x)|\leq |(\widetilde{u})_{x,r_n}-\widetilde{u}(x)|+(\norm{\widetilde{u}^0}_\infty + 1) c_n^{-1} r_n.
\]
Since $(\widetilde{u})_{x,r_n}\to\widetilde{u}(x)$ and $c_n^{-1} r_n\to 0$ as $n\to\infty$, the uniform continuity of $f^\infty$ on compact sets implies \RED{(perhaps after modifying $\dd D^c\widetilde{u}^n/\dd |D^c\widetilde{u}^n|$ on a $|D^c\widetilde{u}^n|$-null set to guarantee that $|\dd D^c\widetilde{u}^n/\dd |D^c\widetilde{u}^n||\equiv 1$ in $\bB^d$ for every $n$)} that
\[
\Bigg|\chi(\widetilde{u}^{n}(z))\Bigg(f^\infty\left(x+r_nz,\widetilde{u}(x+r_nz),\frac{\dd D^c\widetilde{u}^n}{\dd|D^c\widetilde{u}^n|}(z)\right)-f^\infty\left(x,\widetilde{u}(x),\frac{\dd D^c\widetilde{u}^n}{\dd|D^c\widetilde{u}^n|}(z)\right)\Bigg)\Bigg|\to 0
\]
uniformly on $\bB^d$ as $n\to\infty$.
Thus,
\[
\frac{\dd\mu}{\dd|D\widetilde{u}|}(x)\geq\liminf_{n\to\infty}\int_{\bB^d}\chi(\widetilde{u}^{n}(z))f^\infty\left(x,\widetilde{u}(x),\frac{\dd D^c\widetilde{u}^n}{\dd|D^c\widetilde{u}^n|}(z)\right)\;\dd |D^c\widetilde{u}^n|(z).
\]
Condition~\eqref{eqlocalcond2.75} now guarantees that (recall $c_n = r_n^d / |D\widetilde{u}|(B(x,r_n))$)
\[
\lim_{n\to\infty}(|\nabla \widetilde{u}^{n}|\lL+|D^j\widetilde{u}^n|)(\bB^d)=\lim_{n\to\infty}\frac{(|\nabla \widetilde{u}|\lL+|D^j\widetilde{u}|)(B(x,r_n))}{|D\widetilde{u}|(B(x,r_n))}=0,
\]
from which we deduce
\begin{align*}
&\liminf_{n\to\infty}\int_{\bB^d}\chi(\widetilde{u}^{n}(z))f^\infty\left(x,\widetilde{u}(x),\frac{\dd D^c\widetilde{u}^n}{\dd|D^c\widetilde{u}^n|}(z)\right)\;\dd |D^c\widetilde{u}^n|(z)\\
=&\liminf_{n\to\infty}\int_{\bB^d}\int_0^1\chi((\widetilde{u}^{n})^\theta(z))f^\infty\left(x,\widetilde{u}(x),\frac{\dd D\widetilde{u}^n}{\dd|D\widetilde{u}^n|}(z)\right)\;\dd\theta\;\dd |D\widetilde{u}^n|(z)\\
=&\liminf_{n\to\infty}\int_{\bB^d\times\R^m}\chi(w)f^\infty\left(x,\widetilde{u}(x),\frac{\dd\gamma\asc{\widetilde{u}^n}}{\dd|\gamma\asc{\widetilde{u}^n}|}(z,w)\right)\;\dd|\gamma\asc{\widetilde{u}^n}|(z,w).
\end{align*}
By Condition~\eqref{eqlocalcond2} in combination with Theorem~\ref{thmdiffusetangentlifting}, $\gamma\asc{\widetilde{u}^n}\to\gamma\asc{\widetilde{u}^0}$ strictly in $\mbfM(\bB^d\times\R^m;\R^{m\times d})$, and so we can use Reshetnyak's Continuity Theorem to obtain
\[
\frac{\dd\mu}{\dd|D\widetilde{u}|}(x)\geq\int_{\bB^d\times\R^m}\chi(w)f^\infty\left(x,\widetilde{u}(x),\frac{\dd\gamma\asc{\widetilde{u}^0}}{\dd|\gamma\asc{\widetilde{u}^0}|}(z,w)\right)\;\dd|\gamma\asc{\widetilde{u}^0}|(z,w).
\]
\RED{By Condition~\eqref{eqlocalcond2}, we have that $\dd D^c\widetilde{u}^0/\dd|D^c\widetilde{u}^0|(z)=\dd D^c\widetilde{u}/\dd|D^c\widetilde{u}|(x)$ for $|D^c\widetilde{u}|$-almost every $z\in\bB^d$. Definition~\eqref{defelementaryliftings} therefore implies that $\dd\gamma\asc{\widetilde{u}^0}/\dd|\gamma\asc{\widetilde{u}^0}|(z,w)=\dd D^c\widetilde{u}/\dd |D^c \widetilde{u}|(x)$ for $|\gamma\asc{\widetilde{u}}|$-almost every $(z,w)\in\bB^d\times\R^m$.}
Since $\chi(y)=1$ for $|y|\leq\norm{\widetilde{u}^0}_\infty$ and $\supp|\gamma\asc{\widetilde{u}^0}|\subset\overline{\bB^d\times B^m(0,\norm{\widetilde{u}^0}_\infty)}$, we therefore have that
\[
\frac{\dd\mu}{\dd|D\widetilde{u}|}(x)\geq f^\infty\left(x,\widetilde{u}(x),\frac{\dd D^c{\widetilde{u}}}{\dd|D^c{\widetilde{u}}|}(x)\right)\int_{\bB^d\times\R^m}\;\dd|\gamma\asc{\widetilde{u}^0}|(z,w)=f^\infty\left(x,\widetilde{u}(x),\frac{\dd D^c{\widetilde{u}}}{\dd|D^c{\widetilde{u}}|}(x)\right).
\]
Since Condition~\eqref{eqlocalcond2.75} guarantees that $\lim_{r\to 0}|D^c\widetilde{u}|(B(x,r))/|D\widetilde{u}|(B(x,r))=1$ and hence that $\dd\mu/\dd|D\widetilde{u}|(x)=\dd\mu/\dd|D^c\widetilde{u}|(x)$ \RED{and $\rho(\{|\widetilde{u}-y_0|<\tau R\}\setminus\Bcal)=0$, we have therefore obtained~\eqref{equtildecantorbound}.}

\RED{\emph{Step~5:} Finally, we show that~\eqref{equtildelebesguebound} and~\eqref{equtildecantorbound} imply the conclusion of the lemma. Since (see Theorem~3.84 in~\cite{AmFuPa00FBVF}), $D\widetilde{u}$ admits the decomposition
\[
D\widetilde{u}=Du\restrict A^1 + (u^{\Fcal A}-y_1)\otimes\nu_{\Fcal A}\Hcal^{d-1}\restrict\Fcal A,
\]
where $\Fcal A$ denotes the $\Hcal^{d-1}$-rectifiable reduced (or measure-theoretic) boundary of $A$, $u^{\Fcal A}$ is the inner trace, $\nu_{\Fcal A}$ is the (measure-theoretic) unit inner normal, and $A^1$ is the set of points of density $1$ of $A$ (in the sense that $\omega_d^{-1}\lim_{r\to 0}r^{-d}\lL(B(x,r)\cap A)=1$). As $\Hcal^{d-1}(A \setminus A^1) < \infty$ (see Theorem~3.61 in~\cite{AmFuPa00FBVF}), we have  $(\lL+|D^cu|)(A\setminus A^1)=0$, whereby it follows that
\[
\nabla \widetilde{u}\lL\restrict\Omega=\nabla u\lL\restrict A\quad\text{ and }D^c\widetilde{u}=D^c u\restrict A.
\]
This implies
\[
\nabla\widetilde{u}(x)=\nabla u(x)\quad\text{ for $\lL$-almost every $x\in A$}
\]
and
\[
\frac{\dd D^c\widetilde{u}}{\dd|D^c\widetilde{u}|}(x)=\frac{\dd D^c u}{\dd|D^cu|}(x),\quad\frac{\dd\mu}{\dd|D^c\widetilde{u}|}(x)=\frac{\dd\mu}{\dd|D^c{u}|}(x)\quad\text{ for $|D^c u|$-almost every $x\in A$.}
\]
Since we also have that $u=\widetilde{u}$ for $\Hcal^{d-1}$-almost every (and hence $(\lL+|D^cu|)$-almost every) $x\in A$ and $\{|u-y_0|<\tau R\}\subset A$, we therefore see that~\eqref{equtildelebesguebound} and~\eqref{equtildecantorbound} imply
\begin{equation*}
\frac{\dd\mu}{\dd\lL}(x)\geq f(x,{u}(x),\nabla {u}(x))\quad\text{ for }\lL\text{-a.e.\ }x\in \{|u-y_0|<\tau R\}
\end{equation*}
and
\begin{equation*}
\frac{\dd\mu}{\dd|D^c {u}|}(x)\geq f^\infty\left(x,{u}(x),\frac{\dd D^c{u}}{\dd|D^c{u}|}(x)\right)\quad\text{ for }|D^c{u}|\text{-a.e.\ }x\in \{|u-y_0|<\tau R\}.
\end{equation*}
Letting $\tau\uparrow 1$ then provides us with the desired conclusion  at $(\lL+|D^cu|)$-almost every $x\in\Omega$ such that $|u(x)-y_0|<R$.}
\end{proof}

\begin{proof}[Proof of Proposition~\ref{propl1lebesgueineq}]
\RED{In the following we will show that the conclusion of the proposition holds at $(\lL + |D^cu|)$-almost every point $x\in\Omega$ for which $g(x,u(x))>0$. If $x\in\Omega$ is such that $g(x,u(x))=0$, then
\[
f(x,u(x),A)=f^\infty(x,u(x),A)=0\quad\text{ for all $A\in\R^{m\times d}$}
\]
and the conclusion is immediate in this case as well.}

\emph{Step 1:} Assume first that $|\nabla u(x)|>0$ for $\lL$-almost every $x\in\Omega$ and let $g\in\C(\Omega\times\R^m;[0,\infty))$ be arbitrary. Let 
\[
\{B(x_k,r_k)\times B(y_k,R_k)\colon k\in\mbN\}
\]
be a countable collection of cylinders in $\Omega\times\R^m$ such that
\[
\bigcup_{k\in\mbN}B(x_k,r_k)\times B(y_k,R_k)=g^{-1}((0,\infty))
\]
and $\inf\{g(x,y)\colon (x,y)\in B(x_k,r_k)\times B(y_k,R_k)\}>0$ for each $k\in\mbN$. Applying \RED{Lemma~\ref{lem:special}} to each cylinder by taking $\Omega= B(x_k,r_k)$, $y_0=y_k$, $R=R_k$ we obtain that
\[
\frac{\dd\mu}{\dd\lL}(x)\geq f(x,u(x),\nabla u(x))\quad\text{ for }\lL\text{-a.e.\ }x\in \{x\in B(x_k,r_k)\colon u(x)\in B(y_k,R_k)\}
\]
and
\[
\frac{\dd\mu}{\dd|D^cu|}(x)\geq f^\infty\left(x,u(x),\frac{\dd D^cu}{\dd|D^cu|}(x)\right)
\]
for $|D^cu|\text{-a.e.\ }x\in\{x\in B(x_k,r_k)\colon u(x)\in B(y_k,R_k)\}$.
The conclusion now follows from the countable additivity of $\lL\restrict\Omega$ and $|D^cu|$.

\emph{Step 2:} Finally, let $u\in\BV(\Omega;\R^m)$ be arbitrary and define $U\in\BV(\Omega;\R^{d+m})$ by $U(x):=\gr^u(x)=(x,u(x))$ for $\lL$-almost every $x\in\Omega$ so that $\nabla U=(\id_{\R^m},\nabla u)$ and $D^s U=(0,D^su)$. Define $F\in\Rbf(\Omega\times\R^{d+m})$ by 
\[
F(x,(z,y),(B,A)):=f(x,y,A)\quad\text{ for }x,\;z\in\Omega,\,B\in\R^{d\times d},\,A\in\R^{m\times d},
\]
so that, for $\lL$-almost every $x\in\Omega$ and $(|D^cU|=|D^cu|)$-almost every $x\in\Omega$ respectively,
\[
F(x,U(x),\nabla U(x))=f(x,u(x),\nabla u(x))
\]
and
\[
F^\infty\left(x,U(x),\frac{\dd D^c U}{\dd|D^cU|}(x)\right)=f^\infty\left(x,u(x),\frac{\dd D^c u}{\dd|D^cu|}(x)\right).
\]
Since $|\nabla U(x)|=\sqrt{|\id_{\R^d}|^2+|\nabla u(x)|^2}=\sqrt{d^2+|\nabla u(x)|^2}$ for $\lL$-almost every $x\in\Omega$ and it is clear that the quasiconvexity of $f$ over $\R^{m\times d}$ implies that $F$ is quasiconvex over $\R^{(d+m)\times d}$, we see that Step~1 implies the desired conclusion.
\end{proof}

\RED{\begin{remark}
Note that in the proof of Proposition~\ref{propl1lebesgueineq} it would not be possible to bound
\[
\liminf_{j\to\infty}\int_{\bB^d}\chi(\widetilde{w}_j^n(z))f_n(z,\Phi[\widetilde{w}^n_j](z),\nabla \Phi[\widetilde{w}^n_j](z))\;\dd z
\]
 from below for each $n\in\mbN$ by simply applying Theorem~\ref{thmwsclsc} to the weakly* convergent sequence $(\widetilde{w}^n_j)_j\subset\BV(\bB^d;\R^m)$ and the integrand $(x,y,A)\mapsto \chi(y)f_n(z,\varphi(y),\nabla\varphi(y)A)$. This is so because, although the contrapositive of~\eqref{l1blowupeq3} guarantees that eventually 
\[
\varphi((\widetilde{w}_j^n)^+(z))=\varphi((\widetilde{w}_j^n)^-(z))=0\quad\text{ for all }z\in\Jcal_{\widetilde{w}_j^n}\quad\text{ and all }j\in\mbN,
\]
it need not be the case that $\varphi(\theta(\widetilde{w}_j^n)^+(z) +(1-\theta)(\widetilde{w}_j^n)^-)(z)=0$ for all $\theta\in(0,1)$. Consequently, it might occur that
\[
\int_{\Jcal_{\widetilde{w}_j^n}}\int_0^1 \chi((\widetilde{w}_j^n)^\theta(z)) f_n^\infty\left(z,\varphi((\widetilde{w}_j^n)^\theta(z)),\nabla\varphi((\widetilde{w}_j^n)^\theta(z))\frac{\dd D^j\widetilde{w}_j^n}{\dd|D^j\widetilde{w}_j^n|}(z)\right)\;\dd\theta\;\dd|D^j\widetilde{w}_j^n|(z)>0
\]
uniformly in $j$ and $n$ and hence that, for some $\delta>0$,
\begin{equation}\label{eqwrongestimate}
\liminf_{j\to\infty}\Fcal^n_\varphi[\widetilde{w}^n_j]>\liminf_{j\to\infty}\int_{\bB^d}\chi(\widetilde{w}_j^n(z))f_n(z,\Phi[\widetilde{w}^n_j](z),\nabla \Phi[\widetilde{w}^n_j](z))\;\dd z + \delta
\end{equation}
for every $n$, where $\Fcal^n_\varphi[v]$ is defined for every $v\in\BV(\bB^d;\R^m)$ by
\begin{align*}
\Fcal^n_\varphi[v]:=&\int_{\bB^d} f_n(z,\varphi(v(z)),\nabla\varphi(v(z))\nabla v(z))\;\dd z\\
 &\qquad+ \int_{\bB^d}\int_0^1 f_n^\infty\left(z,\varphi(v^\theta(z)),\nabla\varphi(v^\theta(z))\frac{\dd D^s v}{\dd|D^s v|}(z)\right)\;\dd\theta\;\dd|D^sv|(z).
\end{align*}
This precludes any application of Theorem~\ref{thmwsclsc} to the sequence $(\widetilde{w}_j^n)_j$, since this could only be used to provide a potentially non-optimal lower bound for the left-hand side of~\eqref{eqwrongestimate}.
\end{remark}}

\section{Localisation over \texorpdfstring{$\Jcal_u$}{J_u}}\label{secl1jumplsc}

The following proposition shows that any integrand $f\in\RBVL(\Omega\times\R^m)$ (see Definition~\ref{defrepresentationfl1}) can be well approximated by a sequence of `good' integrands \RED{(i.e., those which satisfy the hypotheses of Theorems~\ref{thmareastrictcontinuity} and~\ref{thmwsclsc})} in a manner which is stable under $\Lp^1(\Omega;\R^m)$ convergence of sequences in $\BV(\Omega;\R^m)$.

\begin{lemma}\label{propfphi}
For $\varphi\in\C_b^1(\R^m;\R^m)$ and $f\in\RBVL(\Omega\times\R^m)$, define $f_\varphi\colon\Omega\times\R^m\times\R^{m\times d}\to[0,\infty)$ by
\[
f_\varphi(x,y,A):=f(x,\varphi(y),\nabla\varphi(y)A).
\]
There exists a sequence $(\varphi_M)_M\subset\C^\infty_b(\R^m;\R^m)$ such that for any $f\in\RBVL(\Omega\times\R^m)$, if we abbreviate $f_M:=f_{\varphi_M}$, the following conditions hold:
\begin{enumerate}[(i)]
\item\label{propfcond1} $\sup_M\norm{\nabla\varphi_M}_\infty<\infty$, $\norm{\varphi_M}_\infty\leq 2M^2$, and $\varphi_M(y)\to y$ as $M\to\infty$ for all $y\in\R^m$.
\item\label{propfcond2} For any sequences $r_j\downarrow 0$ and $(u_j)_j$ with $u_j\in(\C^\infty\cap\W^{1,1})(B(x_0,r_j);\R^m)$,
it holds that
\begin{align*}
\limsup_{M\to\infty}\limsup_{j\to\infty}r_j^{1-d}&\int_{B(x_0,r_j)}f_M(x,u_j(x),\nabla u_j(x))\;\dd x\\
\leq\limsup_{j\to\infty}r_j^{1-d}&\int_{B(x_0,r_j)}f(x,u_j(x),\nabla u_j(x))\;\dd x.
\end{align*}
\item\label{propfcond3} If $u\in\BV(\Omega;\R^m)$ with $\Fcal[u]<\infty$, then
\[
\int_\Omega f^\infty_M\left(x,u(x),\frac{\dd D^cu}{\dd|D^cu|}(x)\right)\;\dd|D^cu|(x)\to\int_\Omega f^\infty\left(x,u(x),\frac{\dd D^cu}{\dd|D^cu|}(x)\right)\;\dd|D^cu|(x),
\]
\begin{align*}
&\int_\Omega\int_0^1 f^\infty_M\left(x,u^\theta(x),\frac{\dd D^ju}{\dd|D^ju|}(x)\right)\;\dd\theta\;\dd|D^ju|(x)\\
\to&\int_\Omega\int_0^1 f^\infty\left(x,u^\theta(x),\frac{\dd D^ju}{\dd|D^ju|}(x)\right)\;\dd\theta\;\dd|D^ju|(x),
\end{align*}
as $M\to\infty$, and, if $\int_{\Jcal_u}H_f[u](x)\;\dd\Hcal^{d-1}(x)<\infty$,
\[
\int_{\Jcal_u} H_{f_M}[u]\;\dd\Hcal^{d-1}(x)\to\int_{\Jcal_u} H_{f}[u]\;\dd\Hcal^{d-1}(x)\quad\text{ as }M\to\infty,
\]
where $H_f$ is the surface energy introduced \RED{in~\eqref{eq:Hdef}} at the end of Section~\ref{secpreliminaries}.

If $\Fcal[u]<\infty$ and also $\int_\Omega g(x,u(x))\;\dd x<\infty$, then
\[
\int_\Omega f_M(x,u(x),\nabla u(x))\;\dd x\to\int_\Omega f(x,u(x),\nabla u(x))\;\dd x.
\]

\item\label{propfcond4} Define $h_M:=f_M-f^\infty_M$. \RED{For every $\varepsilon >0$, there exists $R>0$ such that
\[
|h_M(x,y,A)|\leq\varepsilon g(x,\varphi_M(y))(1+|\nabla\varphi_M(y)A|)
\]
for all $(x,y)\in\Omega\times\R^m$ and $A\in\R^{m\times d}$ with $|A|\geq R$.}

\end{enumerate}
\end{lemma}

\begin{proof}
Let $(\varphi_M)_{M\in\mbN}\subset\C^\infty(\R^m;\R^m)$ be a family of functions such that each $\varphi_M$ satisfies
\begin{equation}\label{eqphiproperties}
\begin{cases}
\varphi_M(y)=y &\text{ if }|y|\leq M,\\
M\leq|\varphi_M(y)|\leq\RED{|y|}\text{ and }|\nabla\varphi_M(y)|\leq C &\text{ if }|y|\in[M,M	^2],\\
|\varphi_M(y)|= \RED{M}\text{ and }|\nabla\varphi_M(y)|\leq \frac{C}{M} &\text{ if }|y|\geq M^2,\\
\varphi_M(y)\text{ is positive scalar multiple of } y &\text{ for all }y\in\R^m,
\end{cases}
\end{equation}
where $C>0$ depends only on $m$. We can construct such a family as follows: let $\kappa\in\C^\infty([0,1];[0,1])$ be such that $\kappa(0)=1$, $\kappa(1)=0$, $\frac{\dd^n}{\dd^n x} \kappa \big|_{x=0,1}=0$ for all $n\in\mbN$, $\norm{\kappa'}_\infty\leq 2$, and define
\[
\varphi_M(y):=\begin{cases}
y &\text{ if }|y|\leq M,\\
\kappa\left(\frac{|y|-M}{M^2-M}\right)y+\left(1-\kappa\left(\frac{|y|-M}{M^2-M}\right)\right)M\frac{y}{|y|}&\text{ if }|y|\in[M,M^2],\\
M\frac{y}{|y|}&\text{ if }|y|\geq M^2.
\end{cases}
\]
It follows immediately that $(\varphi_M)_M$ satisfies~\eqref{propfcond1}.

To obtain~\eqref{propfcond2}, note first that we can assume without loss of generality that
\[
\limsup_{j\to\infty}r_j^{1-d}\int_{B(x_0,r_j)}f(x,u_j(x),\nabla u_j(x))\;\dd x<\infty.
\]
Now observe that, since $\varphi_M(y)=y$ for $|y|\leq M$,
\begin{align*}
&\limsup_{j\to\infty}r_j^{1-d}\int_{B(x_0,r_j)}f_M(x,u_j(x),\nabla u_j(x))\;\dd x\\
\leq&\limsup_{j\to\infty}r_j^{1-d}\int_{B(x_0,r_j)\cap\{|u_j(x)|< M\}}f(x,u_j(x),\nabla u_j(x))\;\dd x\\
&\qquad +\limsup_{j\to\infty} Cr_j^{1-d}\int_{B(x_0,r_j)\cap\{|u_j(x)|\geq M\}}g(x,(\varphi_M\circ u_j)(x))(1+|\nabla(\varphi_M\circ u_j)(x)|)\;\dd x,
\end{align*}
where $g\in\C(\overline{\Omega}\times\R^m;[0,\infty))$ \RED{is a witness to the conditions of Definition~\ref{defrepresentationfl1} for $f$.} For fixed $M$, the sequence $(g(\frarg,(\varphi_M\circ u_j)(\frarg))_j$ is uniformly bounded in $\Lp^\infty(\Omega;\R^m)$, from which it follows that
\[
\limsup_{j\to\infty} Cr_j^{1-d}\int_{B(x_0,r_j)}g(x,(\varphi_M\circ u_j)(x))\;\dd x=0.
\]
Now define the sequence of measures $(\eta_j)_j\subset\mbfM^+(\overline{\bB^d}\times\R^m)$ via their action
\[
\ip{\psi}{\eta_j}:=r_j^{d-1}\int_{B(x_0,r_j)}\psi\left(\frac{x-x_0}{r_j},u_j(x)\right)g(x,u_j(x))|\nabla u_j(x)|\;\dd x,\quad\psi\in\C_0\big(\overline{\bB^d}\times\R^m\big).
\]
Since
\begin{align*}
\norm{\psi}_\infty &\cdot \limsup_{j\to\infty}r_j^{1-d}\int_{B(x_0,r_j)}f(x,u_j(x),\nabla u_j(x))\;\dd x\\
 \geq & \limsup_{j\to\infty}r_j^{1-d}\int_{B(x_0,r_j)}\psi\left(\frac{x-x_0}{r_j},u_j(x)\right)f(x,u_j(x),\nabla u_j(x))\;\dd x\\
\geq&\limsup_{j\to\infty}\RED{r_j^{1-d}}\int_{B(x_0,r_j)}\psi\left(\frac{x-x_0}{r_j},u_j(x)\right)g(x,u_j(x))|\nabla u_j(x)|\;\dd x\\
=&\limsup_{j\to\infty}\ip{\psi}{\eta_j},
\end{align*}
it follows that $(\eta_j)_j$ is a norm-bounded sequence in $\mbfM^+(\overline{\bB^d}\times\R^m)$ and so we can pass to a non-relabelled subsequence in order to assume that $\eta_j\wsc\eta$ for some $\eta\in\mbfM^+(\overline{\bB^d}\times\R^m)$.

Now let $M\in\mbN$ be so large that there exists $C>0$ such that $g(x,y)\leq Cg(x,ty)$ whenever $t>1$ for all $x\in\overline{\Omega}$ and $y\in\R^m$ with $|y|\geq M$ and hence that 
\[
g(x,(\varphi_M\circ u_j)(x))\leq C g(x,u_j(x))\qquad\text{ for all $x\in\Omega$ such that $|u_j(x)|\geq M$.}
\]
Combined with the fact that 
\[
|\nabla(\varphi_M\circ u_j)(x)|=|\nabla\varphi_M(u_j(x))\nabla u_j(x)|\leq|\nabla\varphi_M(u_j(x))| \cdot |\nabla u_j(x)|,
\]
we can then use the bounds on $|\nabla\varphi_M|$ collected in~\eqref{eqphiproperties} to compute
\begin{align*}
&Cr_j^{1-d}\int_{B(x_0,r_j)\cap\{|u_j(x)|\geq M\}}g(x,(\varphi_M\circ u_j)(x))|\nabla(\varphi_M\circ u_j)(x)|\;\dd x \\
\leq &\; Cr_j^{1-d}\int_{B(x_0,r_j)\cap\{|u_j(x)|\geq M\}}g(x,u_j(x))|\nabla\varphi_M(u_j(x))| \cdot |\nabla u_j(x)|\;\dd x\\
\leq &\; Cr_j^{1-d}\int_{B(x_0,r_j)\cap\{M^2\geq|u_j(x)|\geq M\}}g(x,u_j(x))|\nabla u_j(x)|\;\dd x\\
&\qquad +\frac{C}{M}r_j^{1-d}\int_{B(x_0,r_j)\cap\{|u_j(x)|\geq M^2\}}g(x,u_j(x))|\nabla u_j(x)|\;\dd x\\
= &\; C\left(\eta_j(\overline{\bB^d}\times\{y\colon|y|\in[M,M^2]\})+\frac{1}{M}\eta_j(\overline{\bB^d}\times\{y\colon|y|\geq M^2\})\right).
\end{align*}
Using the convergence $\eta_j\wsc\eta$ and the upper semicontinuity of the total variation on compact sets, we can therefore deduce
\begin{align*}
&\limsup_{j\to\infty}r_j^{1-d}\int_{B(x_0,r_j)}f_M(x,u_j(x),\nabla u_j(x))\;\dd x\\
\leq &\limsup_{j\to\infty}r_j^{1-d}\int_{B(x_0,r_j)}f(x,u_j(x),\nabla u_j(x))\;\dd x\\
&\qquad+C\left(\eta(\overline{\bB^d}\times\{y\colon|y|\in[M,M^2]\})+\frac{1}{M}\sup_j\eta_j(\overline{\bB^d}\times\R^m)\right).
\end{align*}
Since $\eta$ is a finite measure on $\overline{\bB^d}\times\R^m$ and $(\eta_j)_j$ is a bounded sequence in $\mbfM(\overline{\bB^d}\times\R^m)$, it follows that
\begin{align}
\begin{split}\label{eqcutoffineq}
\limsup_{M\to\infty}\limsup_{j\to\infty}r_j^{1-d}&\int_{B(x_0,r_j)}f_M(x,u_j(x),\nabla u_j(x))\;\dd x\\
\leq\limsup_{j\to\infty}r_j^{1-d}&\int_{B(x_0,r_j)}f(x,u_j(x),\nabla u_j(x))\;\dd x,
\end{split}
\end{align}
as required.

Next, we deduce~\eqref{propfcond3}: to see that
\[
\lim_{M\to\infty}\int_\Omega f_M(x,u(x),\nabla u(x))\;\dd x=\int_\Omega f(x,u(x),\nabla u(x))\;\dd x,
\]
note
\begin{align*}
\int_\Omega f_M(x,u(x),\nabla u(x))\;\dd x &=\int_{\{|u(x)|< M\}} f(x,u(x),\nabla u(x))\;\dd x\\
&\qquad+\int_{\{|u(x)|\geq M\}}f(x,(\varphi_M\circ u)(x),\nabla(\varphi_M\circ u)(x))\;\dd x,
\end{align*}
and that, since $\norm{\nabla\varphi_M}_\infty$ is bounded independently of $M$, $\int_\Omega g(x,u(x))|\nabla u(x)|\;\dd x<\infty$ (since $f(x,y,A)\geq g(x,y)|A|$), and $\int_\Omega g(x,u(x))\;\dd x<\infty$, we can use the fact that, for $M$ sufficiently large $g(x,(\varphi_M\circ u)(x))\leq Cg(x,u(x))$ whenever $|u(x)|\geq M$ in order to bound
\begin{align}
\begin{split}\label{eqpropfphicomputation}
\int_{\{|u(x)|\geq M\}}f(x,(\varphi_M\circ u)&(x),\nabla\varphi_M(u(x))\nabla u(x))\;\dd x \\
&\leq C\int_{\{|u(x)|\geq M\}}g(x,\varphi_M(u(x)))(1+|\nabla\varphi_M(u(x))| \cdot |\nabla u(x)|)\;\dd x\\
&\leq C\int_{\{|u(x)|\geq M\}} g(x,u(x))(1+|\nabla u(x)|)\;\dd x\\
&\to 0\text{ as }M\to\infty.
\end{split}
\end{align}

Similarly, the estimate
\begin{align*}
&\int_\Omega\int_0^1\mathbbm{1}_{\{|y|\geq M\}}(\RED{u^\theta(x)})g(x,\varphi_M(u^\theta(x)))\;\dd|D^su|(x)\\
&\qquad \leq C\int_\Omega\int_0^1\mathbbm{1}_{\{|y|\geq M\}}(u^\theta(x))g(x,u^\theta(x))\;\dd|D^su|(x)\\
&\qquad \to 0\quad\text{ as }M\to\infty
\end{align*}
implies that
\[
\lim_{M\to\infty}\int_\Omega f^\infty_M\left(x,u(x),\frac{\dd D^c u}{\dd|D^cu|}(x)\right)\;\dd|D^cu|(x)=\int_\Omega f^\infty\left(x,u(x),\frac{\dd D^c u}{\dd|D^cu|}(x)\right)\;\dd|D^cu|(x),
\]
and
\begin{align*}
\lim_{M\to\infty}&\int_\Omega\int_0^1 f^\infty_M\left(x,u^\theta(x),\frac{\dd D^ju}{\dd|D^ju|}(x)\right)\;\dd\theta\;\dd|D^ju|(x)\\
=&\int_\Omega\int_0^1 f^\infty \left(x,u^\theta(x),\frac{\dd D^ju}{\dd|D^ju|}(x)\right)\;\dd\theta\;\dd|D^ju|(x),
\end{align*}
\RED{where $(\theta, x)\mapsto u^\theta(x)$ denotes the jump interpolant of $u$ at $x$ defined in~\eqref{eqdefjumpinterpolant}.}

To see
\[
\lim_{M\to\infty}\int_{\Jcal_u} H_{f_M}[u](x)\;\dd\Hcal^{d-1}(x)=\int_{\Jcal_u} H_{f}[u](x)\;\dd\Hcal^{d-1}(x),
\]
note first that the estimate
\[
0\leq f^\infty_M(x,y,A)\leq Cg(x,\varphi_M(y))|\nabla\varphi_M(y)| \cdot |A|\leq C g(x,y)|A|\leq C f^\infty(x,y,A)
\]
for large enough $M$ implies that $H_{f_M}[u]\leq C H_{f}[u]$ pointwise on $\Jcal_u$. If we can show that $H_{f_M}[u](x)\to H_f[u](x)$ for $\Hcal^{d-1}$-almost every $x\in\Jcal_u$, then the conclusion will follow from the Dominated Convergence Theorem. That $\liminf_{M}H_{f_M}[u](x)\geq H_f[u](x)$ follows from the fact that $\varphi_M\circ v\in\Acal_u(x)$ \RED{with $\norm{\varphi_M\circ v}_{\Lp^1}\leq 2\norm{u_x^\pm}_{\Lp^1}$ (since $|\varphi_M(y)|\leq |y|$ for all $y$)} whenever $M\geq\max\{u^+(x),u^-(x)\}$ and $v\in\Acal_u(x)$ with \RED{$\norm{v}_{\Lp^1}\leq 2\norm{u_x^\pm}_{\Lp^1}$} combined with the definition of $H_f$ and the fact that $\nabla(\varphi_M\circ v)=\nabla\varphi_M(v)\nabla v$. After a change of variables and using $f^\infty$ in place of $f$ (note that $f\in\RBVL(\Omega\times\R^m)$ implies $f^\infty\in\RBVL(\Omega\times\R^m)$), \eqref{eqcutoffineq} reads as
\begin{align*}
\limsup_{M\to\infty}\limsup_{j\to\infty}&\int_{\bB^d}f^\infty_M(x+r_jz,u_j(z),\nabla u_j(z))\;\dd z\\
\leq\limsup_{j\to\infty}&\int_{\bB^d}f^\infty(x+r_jz,u_j(z),\nabla u_j(z))\;\dd z
\end{align*}
for any sequence $(u_j)_j\subset\C^\infty(\bB^d;\R^m)$ and any sequence $r_j\downarrow 0$. Letting $(u_j)_j\subset\Acal_u(x)$, $r_j\downarrow 0$ be such that
\[
H_f[u](x)=\lim_{j\to\infty}\int_{\bB^d}f^\infty(x+r_jz,u_j(z),\nabla u_j(z))\;\dd z,
\]
we then deduce
\[
\limsup_{M\to\infty}H_{f_M}[u](x)\leq\limsup_{M\to\infty}\limsup_{j\to\infty}\int_{\bB^d}f^\infty_M(x+r_jz,u_j(z),\nabla u_j(z))\;\dd z\leq H_f[u](x),
\]
as required.

To prove~\eqref{propfcond4}, fix $\varepsilon>0$ and use the fact that $f\in\RBVL(\Omega\times\R^m)$ to obtain $R_\varepsilon>0$ such that $|f(x,y,A)-f^\infty(x,y,A)|\leq\varepsilon g(x,y)(1+|A|)$ for all $A\in\R^{m\times d}$ with $|A|\geq R_\varepsilon$ and all $(x,y)\in\overline{\Omega}\times 2M^2\overline{\bB^m}$. Using~\eqref{propfcond1}, it follows that
\[
|\mathbbm{1}_{\{|A|\geq R_\varepsilon\}}h_M(x,y,A)|\leq\varepsilon g(x,\varphi_M(y))(1+|\nabla\varphi_M(y)A|).
\]
This finishes the proof.
\end{proof}

\begin{corollary}\label{corl1approx}
Let $f\in\RBVL(\Omega\times\R^m)$ and $u\in\BV(\Omega;\R^m)$ be such that $\Fcal[u]<\infty$ and $\int_\Omega g(x,u(x))\;\dd x<\infty$. Then there exists a sequence $(u_j)_j\subset(\C^\infty\cap\W^{1,1}\cap\Lp^\infty)(\Omega;\R^m)$ such that $u_j\to u$ area-strictly in $\BV(\Omega;\R^m)$ and
\[
\lim_{j\to\infty}\int_\Omega f(x,u_j(x),\nabla u_j(x))\;\dd x= \Fcal[u].
\]
\end{corollary}

\begin{proof}
Given $f\in\RBVL(\Omega\times\R^m)$, let $f_M$ be as in Lemma~\ref{propfphi}. Defining $\Fcal_M\colon\BV(\Omega;\R^m)\to[0,\infty)$ by
\[
\Fcal_M[u]:=\int_\Omega f_M(x,u(x),\nabla u(x))\;\dd x+\int_\Omega\int_0^1 f_M^\infty\left(x,u^\theta(x),\frac{\dd D^su}{\dd|D^su|}(x)\right)\;\dd\theta\;\dd|D^su|(x),
\]
Lemma~\ref{propfphi} implies
\[
\lim_{M\to\infty}\Fcal_M[u]=\Fcal[u].
\]
By Theorem~\ref{thmareastrictcontinuity}, for any sequence $(u_j)_j\subset(\C^\infty\cap\W^{1,1})(\Omega;\R^m)$ converging area-strictly to $u$ (see Proposition~\ref{propfixedbdary}), it holds that
\[
\lim_{j\to\infty}\int_\Omega f_M(x,u_j(x),\nabla u_j(x))\;\dd x=\lim_{j\to\infty}\int_\Omega f(x,\varphi_M\circ u_j(x),\nabla\varphi_M\circ u_j(x))\;\dd x=\Fcal_M[u]
\]
for each $M\in\mbN$. Define $u_j^M\in\RED{(\C^\infty\cap\W^{1,1})}(\Omega;\R^m)$ by $u_j^M:=\varphi_M\circ u_j$. It follows from the area-strict convergence of $u_j$ to $u$ that (see again Theorem~\ref{thmareastrictcontinuity} \RED{applied to the integrand $(x,y,A)\mapsto\sqrt{1+|\nabla\varphi_M(y)A|^2}$})
\begin{align*}
\lim_{j\to\infty}\int_\Omega\sqrt{1+|\nabla u_j^M(x)|^2}\;\dd x & =\lim_{j\to\infty}\int_\Omega\sqrt{1+|\nabla\varphi_M(u_j(x))\nabla u_j(x)|^2}\;\dd x\\
&=\int_\Omega\sqrt{1+|\nabla\varphi_M(u(x))\nabla u(x)|^2}\;\dd x\\
&\qquad+\int_\Omega\int_0^1\left|\nabla\varphi_M(u^\theta(x))\frac{\dd D^su}{\dd|D^su|}(x)\right|\;\dd\theta\;\dd|D^su|(x).
\end{align*}
Applying Lemma~\ref{propfphi} to the integrand $\sqrt{1+|A|^2}$, we see that
\begin{align*}
\lim_{M\to\infty}&\int_\Omega\sqrt{1+|\nabla\varphi_M(u(x))\nabla u(x)|^2}\;\dd x+\int_\Omega\int_0^1\left|\nabla\varphi_M(u^\theta(x))\frac{\dd D^su}{\dd|D^su|}(x)\right|\;\dd\theta\;\dd|D^su|(x)\\
&=\int_\Omega \sqrt{1+|\nabla u(x)|^2}\;\dd x+|D^su|(\Omega).	
\end{align*}
It follows then that we can use a diagonal argument to extract a sequence $w_M:=u_{j_M}^M$ such that $w_M\to u$ area-strictly in $\BV(\Omega;\R^m)$ and 
\[
\int_\Omega f(x,w_M(x),\nabla w_M(x))\;\dd x\to\Fcal[u],
\]
as required.
\end{proof}

\begin{remark} \label{rem:G}
If $f\in\RBVL(\Omega\times\R^m)$ satisfies a stronger bound of the form $g(x,y)|A|\leq f(x,y,A)\leq Cg(x,y)|A|$ (rather than just $g(x,y)|A|\leq f(x,y,A)\leq Cg(x,y)(1+|A|)$) for some $g\in\C(\overline{\Omega}\times\R^m;[0,\infty))$ which verifies the conditions present in Definition~\ref{defrepresentationfl1}, it can be seen that the computation~\eqref{eqpropfphicomputation} in the proof of Lemma~\ref{propfphi} can be carried out without needing to assume that $\int_\Omega g(x,u(x))\;\dd x<\infty$. Consequently, Corollary~\ref{corl1approx} (and, later, Theorems~\ref{thml1recovery} and~\ref{thml1relaxation}) is also valid for integrands satisfying $g(x,y)|A|\leq f(x,y,A)\leq Cg(x,y)|A|$ without needing to assume that $\int_\Omega g(x,u(x))\;\dd x<\infty$.
\end{remark}

\begin{proposition}\label{propl1jumpineq}
If $f\in\RBVL(\Omega\times\R^m)$, and $u\in\BV(\Omega;\R^m)$, $(u_j)_j\subset(\C^\infty\cap\W^{1,1})({\Omega};\R^m)$ are such that
\[
\RED{u_j\to u\text{ in }\Lp^1(\Omega;\R^m)},\qquad f(x,u_j(x),\nabla u_j(x))\lL\restrict\Omega\wsc\mu\text{ in }\mbfM^+(\overline{\Omega}),
\]
it holds that
\[
\frac{\dd\mu}{\dd\Hcal^{d-1}\restrict\Jcal_u}(x_0)\geq H_f[u](x_0)\text{ for }\Hcal^{d-1}\text{-almost every }x_0\in\Jcal_u.
\]	
\end{proposition}

\begin{proof}
We start with the estimate
\[
\frac{\dd\mu}{\dd\Hcal^{d-1}\restrict\Jcal_u}(x_0)=\lim_{r\to 0}\frac{\mu(\overline{B(x_0,r)})}{\omega_{d-1}r^{d-1}}\geq\limsup_{r\to 0}\limsup_{j\to\infty}\frac{r^{1-d}}{\omega_{d-1}}\int_{B(x_0,r)}f(x,u_j(x),\nabla u_j(x))\;\dd x.
\]
It suffices to show that
\[
\limsup_{r\to 0}\limsup_{j\to\infty}r^{1-d}\int_{B(x_0,r)}f_M(x,u_j(x),\nabla u_j(x))\;\dd x\geq \omega_{d-1}H_{f_M}(x_0) \RED{ +\varepsilon_M}
\]
for any $M>0$ where $f_M=f_{\varphi_M}$ is as defined in Lemma~\ref{propfphi} \RED{and $(\varepsilon_M)_{M>0}$ is a sequence converging to zero as $M \to \infty$}, since the conclusion will then follow from statements~\eqref{propfcond2} and~\eqref{propfcond3} of Lemma~\ref{propfphi} \RED{upon taking a specific sequence of $r_j$'s} and letting $M\to\infty$.

We claim first that
\begin{align*}
&\limsup_{r\to 0}\limsup_{j\to\infty}r^{1-d}\int_{B(x_0,r)}f_M(x,u_j(x),\nabla u_j(x))\;\dd x\\
=&\limsup_{r\to 0}\limsup_{j\to\infty}r^{1-d}\int_{B(x_0,r)}f_M^\infty(x,u_j(x),\nabla u_j(x))\;\dd x + \RED{\varepsilon_M}.
\end{align*}
Subtracting $f_M^\infty$ from $f_M$, this is equivalent to showing
\[
\limsup_{r\to 0}\limsup_{j\to\infty}r^{1-d}\int_{B(x_0,r)}h_M(x,u_j(x),\nabla u_j(x))\;\dd x=\RED{\varepsilon_M}
\]
where $h_M:=f_M-f_M^\infty$ is as in Lemma~\ref{propfphi}.

Let $r_j\downarrow 0$ be a sequence such that, taking a subsequence if necessary,
\begin{align*}
\limsup_{r\to 0}\limsup_{j\to\infty}r^{1-d}&\left|\int_{B(x_0,r)}h_M(x,u_j(x),\nabla u_j(x))\;\dd x\right|\\
=\limsup_{j\to\infty}r_j^{1-d}&\left|\int_{B(x_0,r_j)}h_M(x,u_j(x),\nabla u_j(x))\;\dd x\right|.
\end{align*}
 By statement~\eqref{propfcond4} of Lemma~\ref{propfphi} \RED{combined with the $\Lp^1$-convergence of $u_j$ to $u$}, we have that, for $R>0$ large enough,
\begin{align*}
&\limsup_{j\to\infty}\Big|r_j^{1-d}\int_{B(x_0,r_j)}\mathbbm{1}_{\{|A|\geq R\}}\left(\nabla(\varphi_M\circ u_j)(x)\right) h_M(x,u_j(x),\nabla u_j(x))\;\dd x\Big|\\
\leq &\; \varepsilon\limsup_{j\to\infty} r_j^{1-d}\left(\int_{B(x_0,r_j)} g(x,\varphi_M\circ u_j(x))\;\dd x+\int_{B(x_0,r_j)} g(x,(\varphi_M\circ u_j(x)))|\nabla(\varphi_M\circ u_j)(x)|\;\dd x\right)\\
\leq &\; \varepsilon\limsup_{j\to\infty} r_j^{1-d}\left(\int_{B(x_0,r_j)} g(x,\varphi_M\circ u(x))\;\dd x+\int_{B(x_0,r_j)} g(x,(\varphi_M\circ u_j(x)))|\nabla(\varphi_M\circ u_j)(x)|\;\dd x\right)\\
= &\; \varepsilon\limsup_{j\to\infty}r_j^{1-d}\int_{B(x_0,r_j)} g(x,(\varphi_M\circ u_j(x)))|\nabla(\varphi_M\circ u_j)(x)|\;\dd x\\
\leq &\; \varepsilon \limsup_{j\to\infty}r_j^{1-d}\int_{B(x_0,r_j)}f_M(x,u_j(x),\nabla u_j(x))\;\dd x.
\end{align*}
\RED{Using statement~\eqref{propfcond2} from Lemma~\ref{propfphi}, we therefore obtain
\begin{align*}
&\limsup_{j\to\infty}\Big|r_j^{1-d}\int_{B(x_0,r_j)}\mathbbm{1}_{\{|A|\geq R\}}\left(\nabla(\varphi_M\circ u_j)(x)\right) h_M(x,u_j(x),\nabla u_j(x))\;\dd x\Big|\\
\leq &\; \varepsilon \limsup_{j\to\infty}r_j^{1-d}\int_{B(x_0,r_j)}f_M(x,u_j(x),\nabla u_j(x))\;\dd x\\
\leq &\; \varepsilon\limsup_{j\to\infty}r_j^{1-d}\int_{B(x_0,r_j)}f(x,u_j(x),\nabla u_j(x))\;\dd x +\varepsilon_M \\
\leq &\; \varepsilon\frac{\dd\mu}{\dd\Hcal^{d-1}\restrict\Jcal_u}(x_0) +\varepsilon_M.
\end{align*}}
On the other hand, since \RED{$|h_M(x,y,A)|\leq Cg(x,\varphi_M(y))(1+|\nabla\varphi_M(y)A|)$ implies that} the sequence $\mathbbm{1}_{\{|A|< R\}}\left(\nabla(\varphi_M\circ u_j)(x)\right) h_M(x,u_j(x),\nabla u_j(x))$ is bounded from above by $C\norm{g(\frarg,\varphi_M(\frarg))}_\infty(1+ R)$, it also holds that
\begin{align*}
&\limsup_{r\to 0}\limsup_{j\to\infty}\Big|r^{1-d}\int_{B(x_0,r)}\mathbbm{1}_{\{|A|< R\}}\left(\nabla(\varphi_M\circ u_j)(x)\right) h_M(x,u_j(x),\nabla u_j(x))\;\dd x\Big|\\
\leq&\limsup_{r\to 0}r^{1-d}\int_{B(x_0,r)}C\norm{g(\frarg,\varphi_M(\frarg))}_\infty(1+ R)\;\dd x\\
=&\; 0.
\end{align*}
Since $\varepsilon>0$ is arbitrary, we conclude
\[
\lim_{r\to 0}\lim_{j\to\infty}r^{1-d}\int_{B(x_0,r)}h_M(x,u_j(x),\nabla u_j(x))\;\dd x=\RED{\varepsilon_M},
\]
as desired. By a diagonal argument then,
\begin{align*}
\lim_{r\to 0}\lim_{j\to\infty}r^{1-d}\int_{B(x_0,r)}f_M(x,u_j(x),\nabla u_j(x))\;\dd x&=\lim_{n\to\infty}\int_{\bB^d}f_M^\infty(x_0+r_nz,w_n(z),\nabla w_n(z))\;\dd z\\
&\qquad\RED{+\varepsilon_M}
\end{align*}
for some sequence $r_n\downarrow 0$ where $w_n\in\RED{(\C^\infty\cap\W^{1,1})}(\bB^d;\R^m)$ is given by $w_n(z):=u_{j_n}(x_0+r_nz)$ for some subsequence $(j_n)_n\subset\mbN$. Since it can readily be checked that \RED{$w_n\to u^\pm_{x_0}$} in $\Lp^1(\bB^d;\R^m)$, the result now follows from Lemma~\ref{lemslicinglemma} below \RED{(note that $\|v_n\|_{\Lp^1} \leq 2\|u_{x_0}^\pm\|_{\Lp^1}$ eventually for sufficiently large $n$)}.
\end{proof}

The following lemma and its proof \RED{(given in full for clarity)} are minor adaptations of Lemma~3.1 from~\cite{FonMul93RQFB}:
\begin{lemma}\label{lemslicinglemma}
Let $x_0\in\Omega$, $f\in\Rbf(\Omega\times\R^m)$ satisfy \textcolor{black}{$0 \leq f(x,y,A)\leq C(1+|A|)$} and let $(v_n)_n\subset(\C^\infty\cap\W^{1,1})(\bB^d;\R^m)$ be a sequence such that \RED{$v_n\to u^\pm_{x_0}$} in $\Lp^1(\bB^d;\R^m)$. It follows that there exists a sequence $(u_n)_n\subset\Acal_u(x_0)$ such that \RED{$u_n\to u^\pm_{x_0}$} in $\Lp^1(\bB^d;\R^m)$ and, for any sequence $r_n\downarrow 0$,
\begin{align*}
\limsup_{n\to\infty}&\int_{\bB^d}f\left(x_0+r_nz,u_n(z),\nabla u_n(z)\right)\;\dd z\\
&\leq\liminf_{n\to\infty}\int_{\bB^d}f\left(x_0+r_nz,v_n(z),\nabla v_n(z)\right)\;\dd z.
\end{align*}
\end{lemma}
\RED{
\begin{proof}
We can assume that
\[
\liminf_{n\to\infty}\int_{\bB^d}f\left(x_0+r_nz,v_n(z),\nabla v_n(z)\right)\;\dd z<\infty
\]
since the conclusion is an immediate consequence of Proposition~\ref{propfixedbdary} otherwise.

We begin by using Proposition~\ref{propfixedbdary} to obtain a sequence $(w_n)_n\subset\Acal_u(x_0)$ with $w_n\to u^\pm$ strictly in $\BV(\bB^d;\R^m)$ as $n\to\infty$. Now define
\[
\alpha_n:=\sqrt{\norm{v_n-w_n}_{\Lp^1}},\quad k_n\:=n \floorb{1+\norm{w_n}_{\BV(\bB^d;\R^m)}+\norm{v_n}_{\BV(\bB^d;\R^m)}},\quad s_n:=\frac{\alpha_n}{k_n},
\]
where $\floor{t}$ denotes the largest integer less than or equal to $t \in \R$. Since $\alpha_n\downarrow 0$, we can assume that $0\leq\alpha_n<1$ and define
\[
B_{0,n} := (1-\alpha_n)\bB^d,\quad B_{i,n} := (1-\alpha_n + is_n)\bB^d,\quad i=1,\ldots, k_n.
\]
For each $i\in\{1,\ldots,k_n\}$, let $\varphi_i$ be a cut-off function satisfying
\[
\varphi_i\in\C^\infty_0(B_{i,n}),\quad \varphi_i \equiv 1\text{ in }B_{{i-1,n}}, \quad \norm{\nabla\varphi_i}_\infty\leq \frac{A}{s_n}
\]
for some constant $A>0$ independent of $i$, and define
\[
v_n^i:= \varphi_i v_n + (1-\varphi_i)w_n,
\]
so that $(v_n^i)_{i,n}\subset\Acal_u(x_0)$ and $v^i_n\to u^\pm$ as $n\to\infty$ in $\Lp^1(\bB^d;\R^m)$ for each $i$.

We can then use the fact that $|f(x,y,A)|\leq C(1+|A|)$ together with the properties of $\varphi_i$ to estimate
\begin{align*}
&\int_{\bB^d}f(x_0+r_nz,v^i_n(z), \nabla v^i_n(z))\;\dd z \\
&\qquad \leq \int_{\bB^d}f(x_0+r_nz,v_n(z), \nabla v_n(z))\;\dd z \\
&\qquad + C\int_{B_{i,n}\setminus B_{i-1,n}}1 + |\nabla v_n(z)| + |\nabla w_n(z)| + |v_n(z)-w_n(z)| \cdot \frac{A}{s_n}\;\dd z\\
&\qquad + C\int_{\bB^d\setminus B_{0,n}}1 + |\nabla w_n(z)|\;\dd z.
\end{align*}
Taking the average of the above inequality over every layer $B_{i,n}\setminus B_{i-1,n}$, we deduce
\begin{align*}
&\frac{1}{k_n}\sum_{i=1}^{k_n}\int_{\bB^d}f(x_0+r_nz,v^i_n(z), \nabla v^i_n(z))\;\dd z \\
&\qquad \leq \int_{\bB^d}f(x_0+r_nz,v_n(z), \nabla v_n(z))\;\dd z \\
&\qquad + \frac{C}{k_n}\int_{\bB^d}1 + |\nabla v_n(z)| + |\nabla w_n(z)| + |v_n(z)-w_n(z)| \cdot \frac{A}{s_n}\;\dd z\\
&\qquad + C\int_{\bB^d\setminus B_{0,n}}1 + |\nabla w_n(z)|\;\dd z.
\end{align*}
Thus, by our definition of $k_n$ and $\alpha_n$,
\begin{align}\label{eq:averagedbound}
\begin{split}
\frac{1}{k_n}\sum_{i=1}^{k_n}\int_{\bB^d}f(x_0+r_nz,v^i_n(z), \nabla v^i_n(z))\;\dd z & \leq \int_{\bB^d}f(x_0+r_nz,v_n(z), \nabla v_n(z))\;\dd z \\
&\quad + \frac{C}{n} +C\sqrt{\norm{v_n(z)-w_n(z)}_{\Lp^1(\bB^d;\R^m)}}\\
&\quad + C\int_{\bB^d\setminus B_{0,n}}1 + |\nabla w_n(z)|\;\dd z,
\end{split}
\end{align}
where the constant $C$ may have increased in a manner independent of $i$ and $n$.

Next, we argue that
\begin{equation}\label{eqremainderto0}
\limsup_{n\to\infty}\int_{\bB^d\setminus B_{0,n}}1 + |\nabla w_n(z)|\;\dd z=0.
\end{equation}
To see this, let $\varepsilon>0$ and let $n_0$ be large enough that $|Du^\pm|(\bB^d\setminus B_{0,n_0})<\varepsilon$. Since $|Dw_n|\to |Du^\pm|$ strictly in $\mbfM^+(\bB^d)$ (recall again that this follows from the strict convergence of $Dw_n$ to $Du^\pm$ combined with Reshetnyak's Continuity Theorem), we therefore have that $\limsup_{n\to\infty}|Dw_n|(\bB^d\setminus B_{0,n_0})\leq\varepsilon$. Since $\lL(\bB^d\setminus B_{0,n})\to 0$ as $n\to\infty$ and $\bB^d\setminus B_{0,n}\subset\bB^d\setminus B_{0,n_0}$ whenever $n\geq n_0$, it follows that
\[
\limsup_{n\to\infty}\int_{\bB^d\setminus B_{0,n}}1 + |\nabla w_n(z)|\;\dd z\leq \limsup_{n\to\infty}(\lL(\bB^d\setminus B_{0,n})+|Dw_n|(\bB^d\setminus B_{0, n_0}))\leq\varepsilon.
\]
Since $\varepsilon>0$ was arbitrary, we have that $\int_{\bB^d\setminus B_{0,n}}1 + |\nabla w_n(z)|\;\dd z\to 0$ as required.

From~\eqref{eq:averagedbound} and~\eqref{eqremainderto0}, then, we have that, there exists a sequence $(\varepsilon_n)_n$ with $\varepsilon_n\downarrow 0$ such that,
\[
\frac{1}{k_n}\sum_{i=1}^{k_n}\int_{\bB^d}f(x_0+r_nz,v^i_n(z), \nabla v^i_n(z))\;\dd z \leq \int_{\bB^d}f(x_0+r_nz,v_n(z), \nabla v_n(z))\;\dd z  + \varepsilon_n.
\]
It follows that for every $n \in \N$ there must exist $i_n\in\{1,\ldots, k_n\}$ satisfying
\[
\int_{\bB^d}f(x_0+r_nz,v^{i_n}_n(z), \nabla v^{i_n}_n(z))\;\dd z\leq  \int_{\bB^d}f(x_0+r_nz,v_n(z), \nabla v_n(z))\;\dd z  + \varepsilon_n,
\]
and so we may simply define $(u_n)_n\subset\Acal_u(x_0)$ by $u_n:=v^{i_n}_n$ and note that $u_n\to u^\pm$ in $\Lp^1(\bB^d;\R^m)$ as $n\to\infty$ by a diagonal argument to complete the proof.
\end{proof}
}

The results of this section now culminate in the following $\Lp^1$-lower semicontinuity theorem:
\begin{theorem}\label{thml1lsc}
Let $f\in\RBVL(\Omega\times\R^m)$ be such that $f(x,y,\frarg)$ is quasiconvex for every $(x,y)\in\overline{\Omega}\times\R^m$. If $(u_j)_j\subset(\C^\infty\cap\W^{1,1})(\Omega;\R^m)$ and $u\in\BV(\Omega;\R^m)$ are such that $u_j\to u$ in $\Lp^1(\Omega;\R^m)$, then
\begin{align*}
\liminf_{j\to\infty}\Fcal[u_j]\geq &\int_\Omega f(x,u(x),\nabla u(x))\;\dd x+\int_\Omega f^\infty\left(x,u(x),\frac{\dd D^cu}{\dd |D^cu|}(x)\right)\;\dd|D^cu|(x)\\
&\qquad+\int_{\Jcal_u}H_f[u](x)\;\dd\Hcal^{d-1}(x).
\end{align*}
\end{theorem}

\begin{proof}
We can assume without loss of generality that $\sup_j\Fcal[u_j]<\infty$ since the conclusion is trivial otherwise. Combining the discussion at the start of Section~\ref{chapl1lsc} with the lower semicontinuity of the total variation and Propositions~\ref{propl1lebesgueineq} and~\ref{propl1jumpineq}, it follows that
\begin{align*}
\liminf_{j\to\infty}\Fcal[u_j] & =\liminf_{j\to\infty}\int_\Omega f(x,u_j(x),\nabla u_j(x))\;\dd x\\
&\geq \mu(\Omega)\\
&\RED{\geq} \int_\Omega \frac{\dd\mu}{\dd\lL}(x)\;\dd x+\int_\Omega \frac{\dd\mu}{\dd|D^cu|}(x)\;\dd|D^cu|(x)+\int_{\Jcal_u}\frac{\dd\mu}{\dd\Hcal^{d-1}\restrict\Jcal_u}(x)\;\dd\Hcal^{d-1}(x)\\
&\geq\int_\Omega f(x,u(x),\nabla u(x))\;\dd x+\int_\Omega f^\infty\left(x,u(x),\frac{\dd D^cu}{\dd |D^cu|}(x)\right)\;\dd|D^cu|(x)\\
&\qquad+\int_{\Jcal_u}H_f[u](x)\;\dd\Hcal^{d-1}(x),
\end{align*}
as required.	
\end{proof}

\section{Recovery sequences and relaxation}\label{chaprecoveryseqs}
This section is devoted to the construction of recovery sequences to show that the lower bound obtained in Theorem~\ref{thml1lsc} for $\Fcalro$ is attained. We then finish by combining Theorems~\ref{thml1lsc} and~\ref{thml1recovery} (proved below) to finally obtain Theorem~\ref{L1lscthm}.

\subsection{Primitive recovery sequences}
The following proposition explicitly constructs $\Lp^1$-recovery sequences in $\BV(\Omega;\R^m)$ for $\Fcalro$ in the case where $f=f^\infty$, following the same procedure introduced and used in~\cite{RindlerShaw19} to construct recovery sequences for Theorem~\ref{thmwsclsc}.
\begin{proposition}\label{proprecessionrecoveryseq}
Let $f\in\C(\overline{\Omega}\times\R^m\times\R^{m\times d})$ be a positively one-homogeneous integrand and assume that 
\[
0\leq f(x,y,A)\leq C|A|\qquad\text{ for all }(x,y,A)\in\overline{\Omega}\times\R^m\times\R^{m\times d}
\]
for some $C>0$. Let $u\in\BV(\Omega;\R^m)$ be such that
\begin{align*}
\int_\Omega f(x,u(x),\nabla u(x))\;\dd x&+\int_\Omega f\left(x,u(x),\frac{\dd D^cu}{\dd|D^cu|}(x)\right)\;\dd|D^cu|(x)\\
&\qquad\qquad\qquad + \int_{\Jcal_u} H_f[u](x)\;\dd\Hcal^{d-1}(x)<\infty.
\end{align*}
Then, abbreviating $\Hfrak_f[u]=H_f[u]\Hcal^{d-1}\restrict\Jcal_u$, there exists a sequence $(u_j)_j\subset\BV(\Omega;\R^m)$ such that $u_j\to u$ in $\Lp^1(\Omega;\R^m)$, $u_j$ and $\nabla u_j$ converge pointwise $\lL$-almost everywhere to $u$ and $\nabla u$ respectively, and
\begin{align*}
\lim_{j\to\infty}&\left(\int_\Omega f(x,u_j(x),\nabla u_j(x))\;\dd x+\int_\Omega\int_0^1 f\left(x,u_j^\theta(x),\frac{\dd D^su_j}{\dd|D^s u_j|}(x)\right)\;\dd\theta\;\dd|D^su_j|(x)\right)\\
=&\int_\Omega f(x,u(x),\nabla u(x))\;\dd x+\int_\Omega f\left(x,u(x),\frac{\dd D^cu}{\dd|D^c u|}(x)\right)\;\dd|D^c u|(x) +\Hfrak_{f}[u](\Jcal_u).
\end{align*}
\end{proposition}

\begin{proof}
Let $L_0$ be the set of points $x\in\mathcal{J}_u$ which are such that
\begin{enumerate}[(1)]
\item\label{eqrecoverytangentcond} $\Hfrak_f[u]$ and $|Du|$ possess approximate tangent planes at $x$, \RED{and
\[
  H_f[u](x)=\frac{\dd\Hfrak_f[u]}{\dd\Hcal^{d-1}\restrict\Jcal_u}(x);
\]}
\item\label{eqrecoveryblowupcond} $u(x+r\frarg)$ converges strictly  in $\BV(\bB^d;\R^m)$ to $u^\pm_x$ as $r\to 0$ as discussed after Definition~\ref{defjumpblowup}.
\end{enumerate}
Corollary~\ref{correctifiabledensities} implies that $H_f[u]$ is $\Hcal^{d-1}\restrict\Jcal_u$-measurable and so, since $\Jcal_u$ is countably $\Hcal^{d-1}$-rectifiable, we have that Condition~\eqref{eqrecoverytangentcond} is satisfied for $\Hcal^{d-1}$-almost every $x\in\Jcal_u$. Theorem~\ref{thmbvblowup} combined with the definition of $\Jcal_u$ implies that Condition~\eqref{eqrecoveryblowupcond} holds at $\Hcal^{d-1}$-every $x\in\Jcal_u$ and so we have that $\mathcal{H}^{d-1}(\mathcal{J}_u\setminus L_0)=0$.

For $i=1,2$, let $F_i\colon L_0\times(0,1]\to\R$ be the functions defined by
\[
F_1(x,r):=\dashint_{B(x,r)}\left|u^\pm_x\left(\frac{\overline{x}-x}{r}\right)-u(\overline{x})\right|\;\dd \overline{x},
\]
\[
 F_2(x,r):=\frac{1}{\omega_{d-1}}r^{1-d}\frac{|Du|(B(x,r))}{|u^+(x)-u^-(x)|}.
\]
It follows from Conditions~\eqref{eqrecoverytangentcond} and~\eqref{eqrecoveryblowupcond} that $\lim_{r\downarrow 0}F_1(x,r)=0$ and $\lim_{r\downarrow 0}F_2(x,r)=1$ for each $x\in L_0$. Since the $F_i$ are $(\Hcal^{d-1}\restrict\Jcal_u)\times(\Lcal^1\restrict(0,1])$-measurable and hence $(\Hfrak_f[u])\times(\Lcal^1\restrict(0,1])$-measurable, for any $\varepsilon>0$ we can write $L_0$ as the following countable union of increasing $\Hfrak_f[u]$-measurable sets:
\[
L_0=\bigcup_{k\in\mbN}\left\{x\in L_0\colon F_1(x,r)\leq\varepsilon\frac{|Du|(B(x,r))}{r^{d-1}},\;\; F_2(x,r)\in(1-\varepsilon,1+\varepsilon)\text{ for all }r\leq \frac{1}{k} \right\}.
\]
Hence, for fixed $\varepsilon>0$ we can write $L_0=L_0^{\varepsilon}\cup E^{\varepsilon}$ where $\Hfrak_f[u](E^{\varepsilon})<\varepsilon$ and, for some $k_{\varepsilon}\in\mbN$,
\[
L_0^{\varepsilon}\subset\left\{x\in L_0\colon F_1(x,r)\leq\varepsilon\frac{|Du|(B(x,r))}{r^{d-1}},\quad F_2(x,r)\in(1-\varepsilon,1+\varepsilon)\text{ for all }r\leq \frac{1}{k_\varepsilon} \right\}.
\]	
By the outer regularity of Radon measures, there exists an open set $U_\varepsilon$ with $\Jcal_u\subset U_\varepsilon$ and $\lL(U_\varepsilon)<\varepsilon$. 

For a fixed $x\in L_0$, the fact that $\Hfrak_f[u]$ possesses an approximate tangent plane at $x$ implies that
\[
\lim_{r\to 0}\frac{\Hfrak_f[u]\left(B(x,r)\right)}{r^{d-1}}=\omega_{d-1}\frac{\dd\Hfrak_f[u]}{\dd\Hcal^{d-1}\restrict\Jcal_u}(x).
\]
Since $|Du|$ also possesses an approximate tangent plane at $x$ we have that
\[
\lim_{r\to 0}\frac{|Du|(B(x,r))}{r^{d-1}}=\omega_{d-1}\frac{\dd|Du|}{\dd\Hcal^{d-1}\restrict\Jcal_u}(x)=\omega_{d-1}|u^+(x)-u^-(x)|>0,
\]
and so we can deduce that, for all $r>0$ sufficiently small,
\[
\left|\frac{\Hfrak_f[u]\left(B(x,r)\right)}{r^{d-1}}-\omega_{d-1}\frac{\dd\Hfrak_f[u]}{\dd\Hcal^{d-1}\restrict\Jcal_u}(x)\right|\leq\frac{\varepsilon}{2}\frac{|Du|(B(x,r))}{r^{d-1}}.
\]
Similarly the definition of $H_f$  implies that, \RED{for any $\tau\in(0,1)$ fixed, there exist arbitrarily small $r>0$ for which}
\[
\inf_{\substack{v\in\Acal_u(x)\\ \norm{v}_{\Lp^1}\leq 2\norm{u^\pm_{x}}_{\Lp^1}}}\left|\int_{\bB^d}f\left(x+\RED{\tau} rz,v(z),\nabla v(z)\right)\dd z-\omega_{d-1}\frac{\dd\Hfrak_f[u]}{\dd\Hcal^{d-1}\restrict\Jcal_u}(x)\right|\leq\frac{\varepsilon}{2}\frac{|Du|(B(x,r))}{r^{d-1}}
\]
holds. Noting also that, for any $\mu\in\mbfM(\Omega)$, $\mu(\pd B(x,r))>0$ can only be true for countably many balls $B(x,r)\subset\Omega$,  it therefore follows that, for $\tau\in(0,1)$ fixed, the collection
\begin{align*}
F^{\RED{\tau}, \varepsilon}&:=\Bigg\{\overline{B(x,r)}\colon x\in L_0^{\varepsilon},\;r\leq\frac{1}{k_{\varepsilon}},\;\overline{B(x,r)}\Subset U_\varepsilon,\;\Hfrak_f[u](\pd B(x,r))=0,\text{ and }\\
&\inf_{\substack{v\in\Acal_u(x)\\ \norm{v}_{\Lp^1}\leq 2\norm{u^\pm_{x}}_{\Lp^1}}}\left|\int_{\bB^d}f\left(x+\RED{\tau} rz,v(z),\nabla v(z)\right)\dd z-\frac{\Hfrak_f[u]\left(B(x,r)\right)}{r^{d-1}}\right|<\varepsilon\frac{|Du|(B(x,r))}{r^{d-1}}\Bigg\}
\end{align*}
is a fine cover for $L_0^{\varepsilon}$. So, by the Vitali--Besicovitch Covering Theorem~\ref{thm:besi}, there exists a countable, disjoint subcover $\mathcal{F}^{\tau, \varepsilon}\subset F^{\tau, \varepsilon}$ of $L_0^{\varepsilon}$ with respect to $\Hfrak_f[u]$.

Let $\overline{B(x_1,r_1)},\,\overline{B(x_{2},r_{2})}\ldots$ be a sequence of elements from $\Fcal^{\tau, \varepsilon}$ such that there exists an increasing sequence $N_1,\,N_2\ldots$ in $\mbN$ with
\[
\RED{\Hfrak_f[u]}\left(L_0^{\varepsilon}\setminus\bigcup_{i=1}^{N_j}\overline{B(x_{i},r_{i})}\right)\leq\frac{1}{j}.
\]
Let $\eta_\tau\in\C_c^\infty(\bB^d;[0,1])$ be such that $\eta_\tau \equiv 1$ on $\tau\bB^d$. For $i=1\ldots N_j$, let $\RED{v_i^\tau}\in\Acal_u(x_i)$ be such that
\[
\left|\int_{\bB^d}f\left(x_i+\tau r_iz,v_i^\tau(z),\nabla v_i^\tau(z)\right)\;\dd z-\frac{\Hfrak_f[u]\left(B(x_i,{r_i})\right)}{r_i^{d-1}}\right|<\varepsilon\frac{|Du|(B(x_i,{r_i}))}{r_i^{d-1}}
\] 
and
\[
\norm{v_i^\tau}_{\Lp^1(\bB^d)}\leq 2\norm{u_{x_i}^\pm}_{\Lp^1(\bB^d)}.
\]
Define
$(\RED{w_i^\tau})_i\subset\BV(\bB^d;\R^m)$ by
\[
w_i^\tau(z):=\begin{cases}
v_i^\tau\left(\frac{z}{\tau}\right)&\text{ if }|z|<\tau,\\
u^\pm_{x_i}(z)&\text{ if }\tau\leq|z|<1.
\end{cases}
\]
We can now define
\begin{align*}
v_j^{\varepsilon,\tau}(x)&:=\sum_{i=1}^{N_j}w_i^\tau\left(\frac{x-x_i}{r_i}\right)\eta_\tau\left(\frac{x-x_i}{r_i}\right),\\
 u_j^{\varepsilon,\tau}(x)&:=u(x)\left(1-\sum_{i=1}^{N_j}\eta_\tau\left(\frac{x-x_i}{r_i}\right)\right) +v_j^{\varepsilon,\tau}(x).
\end{align*}
Invoking the criterion for membership of $L_0^\varepsilon$ involving $F_1$, that is,
\[
  \dashint_{\bB^d} \abs{u_{x_i}^\pm} \;\dd x \leq \varepsilon r_i^{1-d}|Du|(B(x_i,{r_i}))+\dashint_{B(x_i,{r_i})}|u(x)|\;\dd x,
\]
we have that
\begin{align*}
\norm{v_j^{\varepsilon,\tau}}_{\Lp^1(\Omega)}&\leq \sum_{i=1}^\infty r_i^{d}\norm{w_i^\tau}_{\Lp^1(\bB^d)} \\
&\leq C\sum_{i=1}^\infty r_i^{d}\norm{u_{x_i}^\pm}^1_{\Lp^1(\bB^d)}\\
&\leq C\sum_{i=1}^\infty r_i^{d}\left(\varepsilon r_i^{1-d}|Du|(B(x_i,{r_i}))+\dashint_{B(x_i,{r_i})}|u(x)|\;\dd x\right)\\
&\leq C\left(\frac{\varepsilon}{k_\varepsilon}|Du|\left(\RED{\bigcup_{\overline{B}\in\Fcal^{\tau,\varepsilon}}\overline{B}}\right)+\int_{\bigcup\Fcal^{\tau,\varepsilon}}|u(x)|\;\dd x\right)\\
&\leq C\left(\frac{\varepsilon}{k_\varepsilon}|Du|(U_\varepsilon)+\int_{U_\varepsilon}|u(x)|\;\dd x\right),
\end{align*}
which implies that $v_j^{\varepsilon,\tau}\to 0$ in $\Lp^1(\Omega;\R^m)$ as $\varepsilon\to 0$ uniformly in $\tau$ and $j$. \RED{Since
\[
\int_{\Omega}\left|u(x)\sum_{i=1}^{N_j}\eta_\tau\left(\frac{x-x_i}{r_i}\right)\right|\;\dd x\leq \int_{\bigcup_{\overline{B}\in\Fcal^{\tau,\varepsilon}}\overline{B}}|u(x)|\;\dd x\leq \int_{U_\varepsilon}|u(x)|\;\dd x,
\]
and $\lL(U_\varepsilon)<\varepsilon$, the absolute continuity of the Lebesgue integral implies  that $u\cdot\sum_{i=1}^{N_j}\eta_\tau\left(\frac{\frarg-x_i}{r_i}\right)$ also converges to $0$ in $\Lp^1(\Omega;\R^m)$  uniformly in $\tau$ and $j$.} We therefore have that  $u_j^{\varepsilon,\tau}\to u$ in $\Lp^1(\Omega;\R^m)$ as $\varepsilon\to 0$ uniformly in $\tau$ and $j$.

Now observe
\begin{align}
\begin{split}\label{eqrecoverdecomp1}
&\int_{\bigcup_{i=1}^{N_j} B(x_i,r_i)}\int_0^1 f\left(x,\left(u_j^{\varepsilon,\tau}\right)^\theta(x),\frac{\dd D u_j^{\varepsilon,\tau}}{\dd|Du_j^{\varepsilon,\tau}|}(x)\right)\;\dd\theta\;\dd|Du_j^{\varepsilon,\tau}|(x)\\
=&\sum_{i=1}^{N_j}\Bigg\{\int_{B(x_i,\tau r_i)}f\left(x,v_i^\tau\left(\frac{x-x_i}{\tau r_i}\right),\frac{1}{\tau r_i}\nabla v_i^\tau\left(\frac{x-x_i}{\tau r_i}\right)\right)\;\dd x\\
&\qquad+\int_{\left\{\tau\leq\frac{|x-x_i|}{r_i}<1\right\}}\int_0^1 f\left(x,\left(u_j^{\varepsilon,\tau}(x)\right)^\theta,\frac{\dd D u_j^{\varepsilon,\tau}}{\dd|D u_j^{\varepsilon,\tau}|}(x)\right)\;\dd\theta\;\dd|D u_j^{\varepsilon,\tau}|(x) \Bigg\}.
\end{split}
\end{align}
Changing coordinates, we can manipulate the first term in this expression as follows:
\begin{align}
\begin{split}\label{eqrecoverydecomp1.5}
&\sum_{i=1}^{N_j}\int_{B(x_i,\tau r_i)}f\left(x,v_i^\tau\left(\frac{x-x_i}{\tau r_i}\right),\frac{1}{\tau r_i}\nabla v_i^\tau\left(\frac{x-x_i}{\tau r_i}\right)\right)\;\dd x\\
=&\sum_{i=1}^{N_j}\tau^{d-1} r_i^{d-1}\int_{\bB^d}f\left(x_i+\tau r_i z,v_i^\tau(z),\nabla v_i^\tau(z)\right)\;\dd z.
\end{split}
\end{align}
Since $\Fcal^\varepsilon$ is a fine cover for $L_0^\varepsilon$ with respect to $\Hfrak_f[u]$ and $\Hfrak_f[u](\pd B(x_i,r_i))=0$ for each $\overline{B(x_i,r_i)}\in\Fcal^\varepsilon$, we can write
\begin{align*}
&\lim_{j\to\infty}\Bigg|\Bigg(\tau^{d-1}\sum_{i=1}^{N_j}r_i^{d-1}\int_{\bB^d}f\left(x_i+\tau r_i z,v_i^\tau(z),\nabla v_i^\tau(z)\right)\;\dd z\Bigg)-\Hfrak_f[u](\Jcal_u)\Bigg|\\
\leq\tau^{d-1}&\sum_{\overline{B(x_i,r_i)}\in\mathcal{F}^{\tau, \varepsilon}}r_i^{d-1}\Bigg|\int_{\bB^d}f\left(x_i+\tau r_i z,v_i^\tau(z),\nabla v_i^\tau(z)\right)\;\dd z-\frac{\Hfrak_f[u]\left({B(x_i,r_i)}\right)}{r_i^{d-1}}\Bigg|\\
&\qquad\qquad\qquad+(1-\tau^{d-1})\Hfrak_f[u](\Jcal_u)+\tau^{d-1}\Hfrak_f[u]\left(\Jcal_u\setminus \RED{\bigcup_{\overline{B}\in\Fcal^{\tau,\varepsilon}}\overline{B}}\right).
\end{align*}
By our choice of $v_i^\tau$ and the fact that \RED{$\Hfrak_f[u](\Jcal_u\setminus L_0^\varepsilon)<\varepsilon$ and $\Hfrak_f[u](\Jcal_u\setminus \bigcup_{\overline{B}\in\Fcal^{\tau,\varepsilon}}\overline{B})\leq \Hfrak_f[u](\Jcal_u\setminus L_0^\varepsilon)$ since $L_0^\varepsilon\subset\bigcup_{\overline{B}\in\Fcal^{\tau,\varepsilon}}\overline{B}$}, we therefore have that
\begin{align}
\begin{split}\label{eqrecoverydecomp2}
&\lim_{j\to\infty}\Bigg|\Bigg(\tau^{d-1}\sum_{i=1}^{N_j}r_i^{d-1}\int_{\bB^d}f\left(x_i+\tau r_i z,v_i^\tau(z),\nabla v_i^\tau(z)\right)\;\dd z\Bigg)-\Hfrak_f[u](\Jcal_u)\Bigg|\\
&\leq \tau^{d-1}\sum_{\overline{B(x_i,r_i)}\in\mathcal{F}^{\tau, \varepsilon}}\varepsilon|Du|(B(x_i,r_i))+(1-\tau^{d-1})\Hfrak_f[u](\Jcal_u)+\varepsilon\\
&\leq \varepsilon|Du|(\Omega)+(1-\tau^{d-1})\Hfrak_f[u](\Jcal_u)+\varepsilon.
\end{split}
\end{align}

For each $x$ satisfying $\tau\leq\frac{|x-x_i|}{r_i}<1$,
\[
u_j^{\varepsilon,\tau}(x)=u(x)+\eta_\tau\left(\frac{x-x_i}{r_i}\right)\left[u^\pm_{x_i}\left(\frac{x-x_i}{r_i}\right)-u(x)\right], 
\]
and so we can use the product rule to deduce
\begin{align*}
|Du_j^{\varepsilon,\tau}|\left(\left\{x\colon\tau\leq\frac{|x-x_i|}{r_i}<1\right\}\right)\leq &\frac{1}{r_i}\norm{\nabla\eta_\tau}_\infty\int_{\left\{\tau\leq\frac{|x-x_i|}{r_i}<1\right\}}\left|u^\pm_{x_i}\left(\frac{x-x_i}{r_i}\right)-u(x)\right|\;\dd x\\
&\qquad+r_i^{d-1}(1-\tau^{d-1})\omega_{d-1}|u^+(x_i)-u^-(x_i)|\\
&\qquad+|Du|(B(x_i,r_i))-|Du|(B(x_i,\tau r_i)).
\end{align*}
Since $f$ satisfies $0\leq f(x,y,A)\leq C|A|$ for some $C>0$, we can therefore estimate
\begin{align*}
&\int_{\left\{\tau\leq\frac{|x-x_i|}{r_i}<1\right\}}\int_0^1 f\left(x,\left(u_j^{\varepsilon,\tau}(x)\right)^\theta,\frac{\dd D u_j^{\varepsilon,\tau}}{\dd |D u_j^{\varepsilon,\tau}|}(x)\right)\;\dd\theta\;\dd|D u_j^{\varepsilon,\tau}|(x)\\
&\leq C\int_{\left\{\tau\leq\frac{|x-x_i|}{r_i}<1\right\}}\;\dd|Du_j^{\varepsilon,\tau}|(x)\\
&\leq C\Big(r_i^{d-1}\norm{\nabla\eta_\tau}_\infty\dashint_{B(x_i,r_i)}\left|u^\pm_{x_i}\left(\frac{x-x_i}{r_i}\right)-u(x)\right|\;\dd x\\
&\qquad+r_i^{d-1}(1-\tau^{d-1})\omega_{d-1}|u^+(x_i)-u^-(x_i)|+|Du|(B(x_i,r_i))-|Du|(B(x_i,\tau r_i))\Big).
\end{align*}
The membership criterion for $L_0^{\varepsilon}$ with respect to $F_2$ implies that
\begin{align*}
|Du|(B(x_i,\tau r_i))\geq(1-\varepsilon)\tau^{d-1}\omega_{d-1}r_i^{d-1}|u^+(x_i)-u^-(x_i)|\geq\tau^{d-1}\frac{1-\varepsilon}{1+\varepsilon}|Du|(B(x_i,r_i)),
\end{align*}
which gives
\[
|Du|(B(x_i,r_i))-|Du|(B(x_i,\tau r_i))\leq |Du|(B(x_i,r_i))\left(1-\tau^{d-1}\frac{1-\varepsilon}{1+\varepsilon}\right).
\]
From the same criterion we also deduce
\[
r_i^{d-1}(1-\tau^{d-1})\omega_{d-1}|u^+(x_i)-u^-(x_i)|\leq\frac{1-\tau^{d-1}}{1-\varepsilon}|Du|(B(x_i,r_i)).
\]
Thus, using also the membership criterion for $L_0^\varepsilon$ with respect to $F_1$ to bound 
\[
r_i^{d-1}\norm{\nabla\eta_\tau}_\infty\dashint_{B(x_i,r_i)}\left|u^\pm_{x_i}\left(\frac{x-x_i}{r_i}\right)-u(x)\right|\;\dd x\leq \varepsilon \norm{\nabla\eta_\tau}_\infty|Du|(B(x_i,r_i)),
\]
we obtain
\begin{align*}
&\int_{\left\{\tau\leq\frac{|x-x_i|}{r_i}<1\right\}}\int_0^1 f\left(x,\left(u_j^{\varepsilon,\tau}(x)\right)^\theta,\frac{\dd D u_j^{\varepsilon,\tau}}{\dd |D u_j^{\varepsilon,\tau}|}(x)\right)\;\dd\theta\;\dd|D u_j^{\varepsilon,\tau}|(x)\\
\leq&C\Big(\varepsilon\norm{\nabla\eta_\tau}_\infty|Du|(B(x_i,{r_i}))+\frac{1-\tau^{d-1}}{1-\varepsilon}|Du|(B(x_i,{r_i}))\Big)\\
&\qquad+\left(1-\tau^{d-1}\frac{1-\varepsilon}{1+\varepsilon}\right)|Du|(B(x_i,{r_i}))\\
\leq &\; C\left(\varepsilon\norm{\nabla\eta_\tau}_\infty+3(1-\tau^{d-1})+2\varepsilon\right)|Du|(B(x_i,{r_i}))
\end{align*}
for $\varepsilon > 0$ sufficiently small. Hence,
\begin{align}
\begin{split}\label{eqrecoverydecomp3}
&\left|\sum_{i=1}^{N_j}\int_{\left\{\tau\leq\frac{|x-x_i|}{r_i}<1\right\}}\int_0^1 f\left(x,\left(u_j^{\varepsilon,\tau}(x)\right)^\theta,\frac{\dd D u_j^{\varepsilon,\tau}}{\dd |D u_j^{\varepsilon,\tau}|}(x)\right)\;\dd\theta\;\dd|D u_j^{\varepsilon,\tau}|(x)\right|\\
&\leq C\left(\varepsilon\norm{\nabla\eta_\tau}_\infty+2(1-\tau^{d-1})+2\varepsilon\right)|Du|(U_\varepsilon).
\end{split}
\end{align}
Combining~\eqref{eqrecoverdecomp1},~\eqref{eqrecoverydecomp1.5},~\eqref{eqrecoverydecomp2}, and~\eqref{eqrecoverydecomp3}, we finally deduce
\[
\lim_{\tau\to 1}\lim_{\varepsilon\to 0}\lim_{j\to\infty}\int_{\bigcup_{i=1}^{N_j}B(x_i,r_i)}\int_0^1 f\left(x,\left(u_j^{\varepsilon,\tau}\right)^\theta(x),\frac{\dd D u_j^{\varepsilon,\tau}}{\dd|Du_j^{\varepsilon,\tau}|}(x)\right)\;\dd\theta\;\dd|Du_j^{\varepsilon,\tau}|(x)=\Hfrak_f[u](\Jcal_u).
\]
Since $u_j^{\varepsilon,\tau}\equiv u$ in $\Omega\setminus U_\varepsilon$, we can use a diagonal argument to obtain a sequence $(u_j)_j\subset\BV(\Omega;\R^m)$ satisfying $u_j\to u$ in $\Lp^1(\Omega;\R^m)$, $u_j(x)=u(x)$, $\nabla u_j(x)=\nabla u(x)$ in $\Omega_j:=\Omega\setminus U_{1/j}$ and which is such that
\begin{align*}
\lim_{j\to\infty}&\int_\Omega f(x,u_j(x),\nabla u_j(x))\;\dd x+\int_\Omega\int_0^1 f\left(x,u^\theta_j(x),\frac{\dd D^su_j}{\dd |D^s u_j|}(x)\right)\;\dd\theta\;\dd|D^s u_j|(x)\\
=&\int_\Omega f(x,u(x),\nabla u(x))\;\dd x+\int_\Omega f\left(x,u(x),\frac{\dd D^cu}{\dd |D^c u|}(x)\right)\;\dd|D^c u|(x)+\Hfrak_f[u](\Jcal_u),
\end{align*}
as required.
\end{proof}

Combining Proposition~\ref{proprecessionrecoveryseq} together with the approximation properties afforded to integrands in $\RBVL(\Omega\times\R^m)$ by Lemma~\ref{propfphi} now allows us to construct smooth recovery sequences for the general case:

\begin{theorem}\label{thml1recovery}
Let $f\in\RBVL(\Omega\times\R^m)$ and $u\in\BV(\Omega;\R^m)$ be such that $\int_\Omega g(x,u(x))\;\dd x<\infty$. Then there exists a sequence $(u_j)_j\subset(\C^\infty\cap\W^{1,1})(\Omega;\R^m)$ such that $u_j\to u$ in $\Lp^1(\Omega;\R^m)$ and
\begin{align*}
\lim_{j\to\infty}\Fcal[u_j]=&\int_\Omega f(x,u(x),\nabla u(x))\;\dd x+\int_\Omega f^\infty\left(x,u(x),\frac{\dd D^c u}{\dd|D^c u|}(x)\right)\;\dd|D^c u|(x)\\
&\qquad+\int_{\mathcal{J}_u}H_f[u](x)\;\dd\mathcal{H}^{d-1}(x).
\end{align*}
\end{theorem}

\begin{proof}
For $f\in\RBVL(\Omega\times\R^m)$, let $f_M\colon\Omega\times\R^m\times\R^{m\times d}\to[0,\infty)$ be defined as in Lemma~\ref{propfphi}. Since \RED{$f_M^\infty$} satisfies the hypotheses required for Proposition~\ref{proprecessionrecoveryseq}, we are guaranteed that there exists a sequence $(u_j)_j\subset\BV(\Omega;\R^m)$ such that $u_j\to u$ in $\Lp^1(\Omega;\R^m)$, $u_j\equiv u$ in $\Omega_j$, and
\begin{align}\label{eqrecessioneq1}
\begin{split}
&\int_\Omega f_M^\infty(x,u_j(x),\nabla u_j(x))\;\dd x+\int_\Omega\int_0^1 f_M^\infty\left(x,u_j^\theta(x),\frac{\dd D^su_j}{\dd|D^s u_j|}(x)\right)\;\dd\theta\dd|D^su_j|(x)\\
&\to\int_\Omega f_M^\infty(x,u(x),\nabla u(x))\;\dd x+\int_\Omega f_M^\infty\left(x,u(x),\frac{\dd D^cu}{\dd|D^c u|}(x)\right)\;\dd|D^c u|(x)\\
&\qquad +\int_{\Jcal_u} H_{f_M}[u](x)\;\dd\Hcal^{d-1}(x)
\end{split}
\end{align}
as $j\to\infty$. Now let $h_M:=f_M-f_M^\infty$ be as in Lemma~\ref{propfphi} and let $\varepsilon>0$.

By statement~\eqref{propfcond4} of Lemma~\ref{propfphi} \RED{combined with the $\Lp^1$-convergence of $u_j$ to $u$}, we have that, for $R>0$ sufficiently large,
\begin{align*}
&\left|\lim_{j\to\infty}\int_{\Omega}\mathbbm{1}_{\{|A|> R\}}(\nabla (\varphi_M\circ u_j)(x)) h_M(x,u_j(x),\nabla u_j(x))\;\dd x\right|\\
\leq &\; \varepsilon\left(\int_\Omega g(x,(\varphi_M\circ u)(x))\;\dd x+\lim_{j\to\infty}\int_\Omega g(x,(\varphi_M\circ u_j)(x))|\nabla(\varphi_M\circ u_j)(x)|\;\dd x\right)\\
\leq &\; \varepsilon\left(\int_\Omega g(x,(\varphi_M\circ u)(x))\;\dd x+\lim_{j\to\infty}\int_\Omega f_M^\infty(x,u_j(x),\nabla u_j(x))\;\dd x\right)\\
\leq &\; \varepsilon\Bigg(\int_\Omega g(x,(\varphi_M\circ u)(x))\;\dd x+\int_\Omega f_M^\infty(x,u(x),\nabla u(x))\;\dd x\\
&\qquad+\int_\Omega f^\infty_M\left(x,u(x),\frac{\dd D^cu}{\dd|D^c u|}(x)\right)\;\dd|D^c u|(x)+\int_{\Jcal_u} H_{f_M}[u](x)\;\dd\Hcal^{d-1}(x)\Bigg).
\end{align*}
Thus, since $\varepsilon>0$ was arbitrary,
\[
\lim_{R\to\infty}\lim_{j\to\infty}\left|\int_{\Omega}\mathbbm{1}_{\{|A|> R\}}(\nabla (\varphi_M\circ u_j)(x)) h_M(x,u_j(x),\nabla u_j(x))\;\dd x\right|=0.
\]
Hence, using the $\lL$-almost everywhere pointwise convergence of $u_j$ and $\nabla u_j$ to $u$ and $\nabla u$, we see that
\begin{align}
&\lim_{j\to\infty}\int_\Omega h_M(x,u_j(x),\nabla u_j(x))\;\dd x \notag\\
=&\lim_{R\to\infty}\lim_{j\to\infty}\int_\Omega\mathbbm{1}_{\{|A|\leq R\}}(\nabla (\varphi_M\circ u_j)(x))h_M(x,u_j(x),\nabla u_j(x))\;\dd x \notag\\
=&\lim_{R\to\infty}\int_\Omega\mathbbm{1}_{\{|A|\leq R\}}(\nabla (\varphi_M\circ u)(x))h_M(x,u(x),\nabla u(x))\;\dd x \notag\\
=&\int_\Omega h_M(x,u(x),\nabla u(x))\;\dd x.  \label{eqrecessioneq2}
\end{align}
Adding equations~\eqref{eqrecessioneq1} and~\eqref{eqrecessioneq2}, we therefore obtain
\begin{align*}
\lim_{j\to\infty}\Fcal_M[u_j]=&\int_\Omega f_M(x,u(x),\nabla u(x))\;\dd x+\int_\Omega f_M^\infty\left(x,u(x),\frac{\dd D^cu}{\dd|D^c u|}(x)\right)\;\dd|D^c u|(x)\\
&\qquad +\int_{\Jcal_u} H_{f_M}[u](x)\;\dd\Hcal^{d-1}(x),
\end{align*}
where $\Fcal_M\colon\BV(\Omega;\R^m)\to\R$ is defined for each $M$ by
\[
\Fcal_M[u]:=\int_\Omega f_M(x,u(x),\nabla u(x))\;\dd x+\int_\Omega\int_0^1 f_M^\infty\left(x,u^\theta(x),\frac{\dd D^su}{\dd|D^s u|}(x)\right)\;\dd\theta\dd|D^su|(x).
\]
Applying Theorem~\ref{thmareastrictcontinuity} to each $u_j$ and using a diagonal argument, we can find a new (non-relabelled) sequence $(u_j)_j\subset\RED{(\C^\infty\cap\W^{1,1})}(\Omega;\R^m)$ such that $u_j\to u$ in $\Lp^1(\Omega;\R^m)$ and
\begin{align*}
\lim_{j\to\infty}\Fcal_M[u_j]=&\int_\Omega f_M(x,u(x),\nabla u(x))\;\dd x+\int_\Omega f_M^\infty\left(x,u(x),\frac{\dd D^cu}{\dd|D^c u|}(x)\right)\;\dd|D^c u|(x)\\
&\qquad +\int_{\Jcal_u} H_{f_M}[u](x)\;\dd\Hcal^{d-1}(x).
\end{align*}
By Lemma~\ref{propfphi}, however,
\begin{align*}
\lim_{M\to\infty}&\int_\Omega f_M(x,u(x),\nabla u(x))\;\dd x+\int_\Omega f_M^\infty\left(x,u(x),\frac{\dd D^cu}{\dd|D^c u|}(x)\right)\;\dd|D^c u|(x)\\
&\qquad +\int_{\Jcal_u} H_{f_M}[u](x)\;\dd\Hcal^{d-1}(x)\\
=&\int_\Omega f(x,u(x),\nabla u(x))\;\dd x+\int_\Omega f^\infty\left(x,u(x),\frac{\dd D^cu}{\dd|D^c u|}(x)\right)\;\dd|D^c u|(x)\\
&\qquad +\int_{\Jcal_u} H_{f}[u](x)\;\dd\Hcal^{d-1}(x).
\end{align*}
The statement of the theorem now follows from a diagonal argument to obtain sequences $M_j\to\infty$ and $(v_j)_j\subset\RED{(\C^\infty\cap\W^{1,1})}(\Omega;\R^m)$ such that $v_j\to u$ in $\Lp^1(\Omega;\R^m)$ and
\begin{align*}
\lim_{j\to\infty}\Fcal_{M_j}[v_j]=&\int_\Omega f(x,u(x),\nabla u(x))\;\dd x+\int_\Omega f^\infty\left(x,u(x),\frac{\dd D^cu}{\dd|D^c u|}(x)\right)\;\dd|D^c u|(x)\\
&\qquad +\int_{\Jcal_u} H_{f}[u](x)\;\dd\Hcal^{d-1}(x).
\end{align*}
Recalling the definition of $f_M$ and defining a the new sequence $(u_j)_j\subset\RED{(\C^\infty\cap\W^{1,1})}(\Omega;\R^m)$ by $u_j:=\varphi_{M_j}\circ v_j$ so that $f_{M_j}(x,v_j(x),\nabla v_j(x))=f(x,u_j(x),\nabla u_j(x))$ and noting that $u_j\to u$ in $\Lp^1(\Omega;\R^m)$, we arrive at the desired conclusion.
\end{proof}

\subsection{Relaxation}
Combining the results of Sections~\ref{chapl1lsc},~\ref{secl1jumplsc}, and~\ref{chaprecoveryseqs}, we can finally complete the proof of Theorem~\ref{L1lscthm}, which we restate in a slightly different form as follows:

\begin{theorem}\label{thml1relaxation}
Let $f\in\RBVL(\Omega\times\R^m)$ be such that $f(x,y,\frarg)$ is quasiconvex for every $(x,y)\in\Omega\times\R^m$ and for a fixed $g\in\C(\overline{\Omega}\times\R^m;[0,\infty))$ as in Definition~\ref{defrepresentationfl1} define
\[
\Gcal:=\left\{u\in\Lp^1(\Omega;\R^m)\colon\int_\Omega g(x,u(x))\;\dd x<\infty\right\}.
\]
Then the $\Lp^1$-relaxation of the functional \RED{$\Fcal\colon\BV(\Omega;\R^m)\cap\Gcal\to[0,\infty]$},
\[
\Fcal[u]:=\int_\Omega f(x,u(x),\nabla u(x))\;\dd x+\int_\Omega\int_0^1 f^\infty\left(x,u^\theta(x),\frac{\dd D^su}{\dd|D^su|}(x)\right)\;\dd\theta\;\dd|D^su|(x),
\]
onto $\BV(\Omega;\R^m)\cap\Gcal$ is given by
\begin{align*}
\Fcalro[u]&=\int_\Omega f(x,u(x),\nabla u(x))\;\dd x+\int_\Omega f^\infty\left(x,u(x),\frac{\dd D^c u}{\dd|D^c u|}(x)\right)\;\dd|D^c u|(x)\\
&\qquad+\int_{\mathcal{J}_u}H_f[u](x)\;\dd\mathcal{H}^{d-1}(x).
\end{align*}
\end{theorem}

\begin{proof}
Corollary~\ref{corl1approx} implies that
\[
\Fcalro[u]=\inf\left\{\liminf_{j\to\infty}\Fcal[u_j]\colon(u_j)_j\subset(\C^\infty\cap\W^{1,1}\cap\Lp^\infty)(\Omega;\R^m)\text{ and }u_j\to u\text{ in }\Lp^1(\Omega;\R^m)\right\}.
\]
Theorem~\ref{thml1lsc} therefore implies
\begin{align*}
\Fcalro[u]&\geq\int_\Omega f(x,u(x),\nabla u(x))\;\dd x+\int_\Omega\int_0^1 f^\infty\left(x,u(x),\frac{\dd D^cu}{\dd|D^cu|}(x)\right)\;\dd|D^c u|(x)\\
&\qquad+\int_{\Jcal_u} H_f[u](x)\;\dd\Hcal^{d-1}(x)
\end{align*}
for any $u\in\BV(\Omega;\R^m)$. If $u$ satisfies $\int_\Omega g(x,u(x))\;\dd x<\infty$, then Theorem~\ref{thml1recovery} provides a sequence $(u_j)_j\subset(\C^\infty\cap\W^{1,1})(\Omega;\R^m)$  such that $u_j\to u$ in $\Lp^1(\Omega;\R^m)$ and
\begin{align*}
\Fcalro[u]\leq & \lim_{j\to\infty}\int_\Omega f(x,u_j(x),\nabla u_j(x))\;\dd x\\
=&\int_\Omega f(x,u(x),\nabla u(x))\;\dd x+\int_\Omega f^\infty\left(x,u(x),\frac{\dd D^c u}{\dd|D^c u|}(x)\right)\;\dd|D^c u|(x)\\
&\qquad+\int_{\mathcal{J}_u}H_f[u](x)\;\dd\mathcal{H}^{d-1}(x),
\end{align*}
from which the conclusion follows.
\end{proof}

\bibliography{math_bib}
\bibliographystyle{plain}
\end{document}